\documentclass[3p]{elsarticle}

\usepackage{amsmath,amssymb,amsthm,mathtools}
\usepackage{bm}
\usepackage{dsfont}
\usepackage{stmaryrd}
\usepackage{booktabs}
\usepackage{multirow}
\usepackage{siunitx}
\usepackage{xcolor}
\usepackage{graphicx}
\usepackage{epstopdf}
\usepackage{capt-of}
\usepackage{tikz-cd}
\usepackage{float,newfloat}
\usepackage{comment,ccicons}
\usepackage{hhline,rotating}
\usepackage{placeins}
\usepackage[expansion=false]{microtype}
\usepackage{tikz}
\usetikzlibrary{arrows.meta,positioning,calc,fit,shapes.geometric}
\usepackage{hyperref}
\hypersetup{hidelinks,pdfauthor={S. Cardonna, M. Lupi, L. Gude Vila, W. Boscheri}}
\usepackage[nameinlink,noabbrev]{cleveref}
\crefformat{appendix}{#2#1#3}\Crefformat{appendix}{#2#1#3}
\crefformat{subappendix}{#2#1#3}\Crefformat{subappendix}{#2#1#3}
\graphicspath{{./figures/}}
\DeclareFloatingEnvironment[fileext=lop]{Diagram}

\allowdisplaybreaks

\journal{Computers \& Fluids}

\newcommand{\R}{\mathbb{R}}
\newcommand{\Q}{\bm q}
\newcommand{\B}{\bm B}
\newcommand{\uvec}{\bm u}
\newcommand{\mvec}{\bm m}
\newcommand{\nvec}{\bm n}
\newcommand{\xvec}{\bm x}
\newcommand{\F}{\mathds F}
\newcommand{\I}{\mathds I}
\newcommand{\jump}[1]{\llbracket #1\rrbracket}
\newcommand{\avg}[1]{\{\!\{#1\}\!\}}
\newcommand{\dd}{\,\mathrm d}
\newcommand{\Th}{\mathcal T_h}
\newcommand{\Ehi}{\mathcal E_h^{\mathrm i}}
\newcommand{\Mh}{\mathds M}
\newcommand{\Ch}{\mathds C}
\newcommand{\Ah}{\mathds A}

\newcommand{\frakD}{\mathfrak D}
\newcommand{\frakJ}{\mathfrak J}
\newcommand{\frakL}{\mathfrak L}
\newcommand{\dt}{\Delta t}
\newcommand{\figurepair}[5]{%
	\begin{minipage}[t]{#1\textwidth}\centering
		\includegraphics[width=\linewidth]{#2}\\[-1mm]
		{\footnotesize #3}
	\end{minipage}\hspace{0.015\textwidth}
	\begin{minipage}[t]{#1\textwidth}\centering
		\includegraphics[width=\linewidth]{#4}\\[-1mm]
		{\footnotesize #5}
\end{minipage}}

\begin{document}
	
	\begin{frontmatter}
		
		\title{A locally and globally divergence-free semi-implicit FV/DG scheme for ideal magnetohydrodynamics at all Mach and Alfv\'en numbers}

		\author[WAT,IMAG]{S. Cardonna}
		\ead{sacha.cardonna@uwaterloo.ca}
		
		\author[UniPD]{M. Lupi}
		\ead{mattia.lupi@phd.unipd.it}
		
		\author[UniSC,CITMAga]{L. Gude Vila}
		\ead{lidia.gude.vila@usc.es}
		
		\author[UniSMB,UniFE]{W. Boscheri \corref{cor1}}
		\ead{walter.boscheri@univ-smb.fr}
		\cortext[cor1]{Corresponding author}

		\address[WAT]{Department of Applied Mathematics, University of Waterloo, N2L3G1 Waterloo, Canada}
		
		\address[IMAG]{Institut Montpelliérain Alexander Grothendieck UMR 5149, 34090 Montpellier, France}

		\address[UniPD]{Department of Civil and Environmental Engineering, University of Padova, 35131 Padova, Italy}
		
		\address[UniSC]{Department of Applied Mathematics, Universidade de Santiago de Compostela, 15782 Santiago de Compostela, Spain}
		
		\address[CITMAga]{Galician Center for Mathematical Research and Technology, CITMAga, 15782 Santiago de
		Compostela, Spain}
		
		\address[UniSMB]{Laboratoire de Math\'ematiques UMR 5127 CNRS
			Universit\'e Savoie Mont Blanc,
			73376 Le Bourget du Lac, France}
		
		\address[UniFE]{Department of Mathematics and Computer Science, University of
			Ferrara, 44121 Ferrara, Italy}

		\begin{abstract}
			We present a new semi-implicit hybrid finite volume/discontinuous Galerkin (FV/DG) scheme for the ideal magnetohydrodynamics (MHD) equations, which is locally and globally divergence-free and works uniformly at all Mach and Alfv\'en numbers. Following a flux splitting approach, the governing equations are divided into a convective, a magnetic and a pressure sub-system. The convective sub-system is discretized explicitly with a finite volume scheme, so that the time step is only restricted by the material speed in smooth flows, while the magnetic and pressure sub-systems are treated implicitly. The magnetic field is approximated in each cell with a five-dimensional basis of linear divergence-free polynomials, thus the local divergence vanishes identically by construction. The magnetic forces are evaluated on a subgrid within each cell, which yields a linear system for the magnetic field without introducing additional unknowns. Global solenoidality, i.e. the continuity of the normal magnetic field across the cell faces, is enforced by a constrained projection in the mass-matrix norm, combined with an energy correction that preserves the total energy. The pressure sub-system leads to a symmetric positive definite elliptic equation for the total energy. The resulting first order scheme is extended to second order of accuracy in time by Richardson extrapolation. The new method is validated in terms of accuracy and robustness against a set of benchmarks. Second order of convergence is retrieved for all Mach numbers, the local divergence is exactly zero and the normal jumps of the magnetic field remain at the level of the solver tolerance in all tests, hence ensuring global divergence-free magnetic fields.
		\end{abstract}
		
		
	\end{frontmatter} 
	
	
	\section{Introduction}\label{sec:intro}
	The equations of ideal magnetohydrodynamics (MHD) provide a mathematical description of electrically conducting fluids interacting with magnetic fields. They arise in a broad range of applications, including astrophysical plasmas, space physics, and magnetic and inertial confinement fusion. From the mathematical viewpoint, compressible ideal MHD constitutes a nonlinear hyperbolic system of conservation laws supplemented with the solenoidal constraint on the magnetic field, $\nabla\cdot\B=0$. The simultaneous presence of multiple wave families and of this differential involution poses two distinct challenges for numerical discretizations: accurately resolving solutions across disparate characteristic time scales and maintaining a physically admissible magnetic field by preserving its solenoidal property.
		
	The multiscale character of MHD is conveniently described by the Mach number $\mathcal M$, measuring the ratio between the characteristic fluid velocity and the sound speed, and the Alfv\'en number $\mathcal M_A$, defined analogously using the Alfv\'en speed. When $\mathcal M\ll1$, acoustic waves propagate much faster than material transport and the dynamics approach a low-Mach, nearly incompressible regime. Similarly, $\mathcal M_A\ll1$ characterizes strongly magnetized configurations in which magnetic waves can impose an additional severe restriction on the time step. In contrast, at moderate or large Mach numbers, the numerical method must remain capable of capturing compressive waves, shocks, and contact discontinuities. Hence, a numerical solver intended for all Mach and Alfv\'en numbers should remain robust in compressible regimes without becoming prohibitively expensive or excessively dissipative as acoustic or magnetic characteristic speeds grow.
	
	Classical Godunov-type finite volume schemes provide robust shock-capturing approximations for hyperbolic conservation laws \cite{Godunov1959,HartenLaxVanLeer1983,Balsara2004}. Yet their explicit stability restriction depends on the fastest characteristic wave speed, which can render simulations of low-Mach or low-Alfv\'en flows computationally demanding. In addition, standard upwind numerical dissipation may become overly large as the Mach number tends to zero, thereby compromising the resolution of the flow structures \cite{Dellacherie1,GH}. Several approaches have been proposed to address these difficulties, including suitable modifications of numerical viscosity, fully implicit methods, and semi-implicit or implicit-explicit (IMEX) strategies. The latter offer an attractive compromise: nonlinear convective transport can be handled explicitly by robust shock-capturing finite volume fluxes, while the stiff acoustic and magnetic contributions are integrated implicitly \cite{AscRuuSpi,BosFil2016,BoscheriPareschi2021,Chen2022,BirkeBoscheriKlingenberg2024}.
	
	Semi-implicit methods have also stimulated the development of hybrid spatial discretizations, which assign different numerical techniques to different terms of the PDE system. Staggered discontinuous Galerkin (DG) approaches have been developed for all-Mach compressible flow problems, combining explicit treatment of advection with an implicit pressure solve \cite{TavelliDumbser2017AllMach}. Likewise, hybrid finite volume/finite element (FV/FE) schemes employ finite volume discretizations for nonlinear convective transport and continuous finite elements for the elliptic-type equations associated with the pressure \cite{Busto2021}. Extensions of such ideas to MHD include divergence-preserving semi-implicit finite volume algorithms \cite{DumbserBalsaraTavelliFambri2019,Fambri2021} and hybrid FV/FE formulations for incompressible MHD on unstructured meshes \cite{FambriEtAl2023,ZampaBustoDumbser2024}. Of particular relevance to the present work are all-Mach and all-Alfv\'en semi-implicit MHD schemes in which operator splitting and a suitable time linearization lead to separate implicit solves for the magnetic and pressure sub-systems \cite{BoscheriThomann2024}. These developments motivate the combination of conservative finite volume transport with a DG approximation of the fields entering the implicit sub-system.
	
	Independently of time step restrictions, the magnetic divergence constraint represents a fundamental structural requirement. At the continuous level, the induction equation propagates $\nabla\cdot\B=0$ from divergence-free initial data. In a fully discrete scheme, however, this property is not automatic, and spurious magnetic divergence can compromise the robustness and physical fidelity of the computation. Widely used remedies include constrained transport on staggered meshes \cite{yee1966,BalsaraSpicer1999}, hyperbolic/parabolic divergence cleaning \cite{dedner2002}, and compatible or mimetic discrete differential operators \cite{HymanShashkov1997,BoscheriDimarcoPareschi2023}. An interesting alternative is to construct the magnetic approximation directly in a polynomial space whose basis functions have vanishing divergence inside each element. Such locally divergence-free DG approximations have been investigated for ideal MHD and provide a natural way to eliminate the element-interior part of the divergence error \cite{LiXu2013LocalDG,LiXu2012GlobalDG}.
	
	It is essential, nevertheless, to distinguish local solenoidality from its global counterpart. For a piecewise smooth magnetic field $\B_h$ that is divergence-free in each mesh element, the distributional divergence can still contain face-supported contributions proportional to the jump of the normal magnetic field, $\llbracket\B_h\cdot\nvec\rrbracket$. Consequently, an element-wise divergence-free polynomial basis alone does not yield a globally divergence-free DG solution. Normal continuity across interfaces must additionally be enforced, whether through compatible degrees of freedom, a constrained reconstruction, or an appropriate coupling of numerical traces. This observation also clarifies the role of divergence-free and curl-free approximation spaces in structure-preserving numerical methods: such spaces encode differential identities locally, while inter-element traces and fluxes determine whether the corresponding involutions hold globally \cite{BoscheriDimarcoPareschi2023,ZampaBustoDumbser2024}.
	
	In this paper, we introduce a locally and globally divergence-free semi-implicit hybrid FV/DG method for the compressible ideal MHD equations at all Mach and Alfv\'en numbers. The temporal splitting separates the nonlinear convective terms, treated explicitly by a shock-capturing finite volume scheme, from the acoustic and magnetic contributions, treated implicitly to alleviate restrictions associated with fast waves. The distinguishing feature of the spatial discretization is that the magnetic field itself is the primary unknown: unlike vector-potential formulations, no auxiliary magnetic potential is evolved and no curl reconstruction is required to obtain $\B_h$. Instead, the magnetic field is approximated directly in a genuinely divergence-free local DG basis. Inter-element consistency is enforced through a constrained reconstruction procedure formulated as a minimization problem, in which jumps of the normal magnetic field across element interfaces are constrained to vanish. This procedure complements the locally divergence-free polynomial approximation space by controlling the inter-element contributions to the discrete divergence, thereby ensuring global solenoidality when the normal-flux continuity constraints are satisfied exactly. Finally, the pressure sub-system is solved implicitly following the flux splitting originally proposed in \cite{Toro_Vazquez12}, hence making the new method suitable to deal with low-Mach flows.
	
	The proposed approach combines the complementary advantages of the two spatial discretizations: robust conservative finite volume transport for the nonlinear hyperbolic part and flexible high-order DG approximation for the semi-implicit magnetic sub-system. In this framework, local divergence preservation is guaranteed by the construction of the polynomial basis, while global divergence preservation is achieved through the constrained reconstruction of the magnetic field across element interfaces. This strategy avoids introducing a magnetic vector potential and enables the implicit treatment of fast magnetic waves. The resulting method provides a unified framework for the accurate and efficient simulation of compressible MHD flows across a broad range of Mach and Alfv\'en numbers. 
	
	The remainder of this paper is organized as follows. In \Cref{sec:model} we recall the compressible ideal MHD equations and the adopted flux splitting technique. We then describe the semi-implicit splitting and the hybrid FV/DG discretization, with particular attention to the locally divergence-free magnetic basis and to the constrained projection that enforces the continuity of the normal magnetic field across the cell faces. Finally, numerical experiments assess the accuracy, robustness, and divergence control of the proposed scheme for a range of Mach and Alfv\'en numbers, followed by concluding remarks.

	\section{Governing equations}
	\label{sec:model}
	
	We consider the ideal magnetohydrodynamics (MHD) equations on a domain $\Omega\subset\R^2$. The velocity is $\uvec=(u,v)^\top$, the magnetic field is $\B=(B_x,B_y)^\top$, and the momentum is $\mvec=\rho\uvec$. The vector of conservative variables is
	\begin{equation}
		\Q=(\rho,\mvec,E,\B)^\top .
		\label{eq:state}
	\end{equation}
	Here $\rho$ and $E$ denote density and total energy per unit volume. For a perfect gas with ratio of specific heats $\gamma>1$, the pressure is given by
	\begin{equation}
		p=(\gamma-1)\left(E-\frac{|\mvec|^2}{2\rho}
		-\frac{|\B|^2}{8\pi}\right).
		\label{eq:eos}
	\end{equation}
	The magnetic pressure is $m_B=|\B|^2/(8\pi)$, and $\I$ is the identity tensor. We also introduce the specific kinetic energy $k=|\uvec|^2/2$, the specific internal energy $e=p/((\gamma-1)\rho)$ and the specific enthalpy $h=e+p/\rho$. The mathematical model, in Gaussian units, reads
	\begin{equation}
		\partial_t\Q+\nabla\cdot\F(\Q)=0,
		\qquad \nabla\cdot\B=0,
		\label{eq:mhd}
	\end{equation}
	with the nonlinear flux tensor
	\begin{equation}
		\F(\Q)=
		\begin{pmatrix}
			\rho\uvec\\
			\rho\uvec\otimes\uvec+
			(p+m_B)\I-\B\otimes\B/(4\pi)\\
			(E+p+m_B)\uvec-(\uvec\cdot\B)\B/(4\pi)\\
			\B\otimes\uvec-\uvec\otimes\B
		\end{pmatrix}.
		\label{eq:full-flux}
	\end{equation} 
	Following \cite{BoscheriThomann2024}, the flux can be split into three parts, that is $\F=\F^{\mathrm c}+\F^{\mathrm m}+\F^{\mathrm p}$, where 
	\begin{subequations}\label{eq:split-components}
		\begin{align}
			\F^{\mathrm c}(\Q)
			=\begin{pmatrix}
				\rho\uvec\\ \rho\uvec\otimes\uvec\\ \rho k \uvec \\ \mathds O
			\end{pmatrix},
			\qquad
			\F^{\mathrm m}(\Q)
			=\begin{pmatrix}
				\bm0\\ m_B\I-\B\otimes\B/(4\pi)\\
				2m_B\uvec-(\uvec\cdot\B)\B/(4\pi)\\
				\B\otimes\uvec-\uvec\otimes\B
			\end{pmatrix},\qquad
			\F^{\mathrm p}(\Q)
			=\begin{pmatrix}
				\bm0\\ p\I\\ h \rho \uvec\\ \mathds O
			\end{pmatrix}.
		\end{align}
	\end{subequations}
	Zero blocks have the dimensions of the corresponding flux components. Let us remark that the hydrodynamic energy flux has been reformulated as $(\rho e+\rho k+p)\uvec = \rho k \uvec + h \rho \uvec$, hence allowing the pressure contribution to be split from the pure transport of the kinetic energy. Consequently, the magnetic flux $\F^{\mathrm m}$ carries the whole electromagnetic energy flux $2m_B\uvec-(\uvec\cdot\B)\B/(4\pi)$. The convective sub-system has material and stationary characteristic speeds; its pressure-less momentum block is not strongly hyperbolic in general. For a unit direction $\nvec$, let us define the following speeds: 
	\begin{equation}
		c_s^2=\frac{\gamma p}{\rho}, \qquad c_A^2=\frac{|\B|^2}{4\pi\rho}, \qquad  c_{A,\nvec}^2=\frac{(\B\cdot\nvec)^2}{4\pi\rho}.
	\end{equation}
	The fast magnetosonic speed is
	\begin{equation}
		c_f^2=\frac12\left[
		c_s^2+c_A^2+
		\sqrt{(c_s^2+c_A^2)^2-4c_s^2c_{A,\nvec}^2}\right].
		\label{eq:fast-magnetosonic-speed}
	\end{equation}
	The acoustic Mach number is $\mathcal{M}=|\uvec|/c_s$, and the so-called Alfv\'en Mach number can be defined as $\mathcal{M}_A=|\uvec|/c_A$.

	\section{Numerical scheme}
	\label{sec:num_scheme}
	
	As done in \cite{BoscheriThomann2024,Fambri2021,Dematte2024MHD}, the convective sub-system will be discretized explicitly, while the magnetic and  the pressure sub-systems will undergo an implicit time marching scheme to obtain a solver which uniformly works at all Mach and Alfv\'en numbers.

	\subsection{Semi-discrete scheme in time}
	\label{sec:sdscheme}
	We detail the first order semi-discrete time marching scheme over the time interval $ [t^{n};t^{n+1}]$. The time step $\dt = t^{n+1}-t^n$ is computed under a CFL-type stability condition which is only restricted by the material speed rather than by the acoustic or Alfv\'en speeds, thus
	\begin{equation}
		\dt = \textnormal{CFL} \, \frac{\min_{\Omega} \varrho}{\max_{\Omega}|\uvec|},
		\label{eqn.timestep}
	\end{equation}
	with $\varrho$ being a characteristic mesh size and $\textnormal{CFL}=1/2$.
	
	\begin{enumerate}
		\item Applying a forward Euler method to the convective sub-system, we obtain the intermediate states $\Q^{*}$:
		\begin{equation}
			\Q^{*} = \Q^n - \dt \, \nabla \cdot \F^{\mathrm c}(\Q^n).
			\label{eqn.sd-conv}
		\end{equation}
		
		\item The implicit magnetic sub-system is then solved, which couples the evolution equations for $\mvec$ and $\B$:
		\begin{subequations}
			\begin{align}
				\mvec^{**} & = \mvec^{*} - \dt \, \nabla \cdot \left( \frac{\B^n \cdot \B^{n+1}}{8\pi}\I- \frac{\B^{n+1} \otimes \B^{n}}{4\pi}  \right), \label{eqn.sd-magnetic-rhou} \\
				\B^{n+1} & = \B^n - \dt \, \nabla \cdot \left( \B^n \otimes \uvec^{**}-\uvec^{**}\otimes\B^n \right). \label{eqn.sd-magnetic-B}
			\end{align}	
			\label{eqn.sd-magnetic}
		\end{subequations}
		Let us notice that $\rho^{n+1}=\rho^*$, and that a semi-implicit discretization has been adopted in order to linearize the quadratic terms involving the magnetic forces. Plugging the expression for the new velocity from \eqref{eqn.sd-magnetic-rhou} into \eqref{eqn.sd-magnetic-B} leads to a vector wave equation for the magnetic field $\B^{n+1}$:
		\begin{align}
			\B^{n+1} & - \dt^2 \, \nabla \cdot \left( \B^n \otimes \frac{1}{\rho^*} \nabla \cdot  \left( \frac{\B^n \cdot \B^{n+1}}{8\pi}\I- \frac{\B^{n+1} \otimes \B^{n}}{4\pi}  \right)  \right) \nonumber \\
			&+ \dt^2 \, \nabla \cdot \left( \left( \frac{1}{\rho^*} \nabla \cdot  \left( \frac{\B^n \cdot \B^{n+1}}{8\pi}\I- \frac{\B^{n+1} \otimes \B^{n}}{4\pi}  \right)  \right) \otimes \B^n \right) = \bm b^n,
			\label{eqn.sd-waveEqn_B}
		\end{align}
		with the known right hand side
		\begin{equation}
			\bm b^n = \B^{n}  - \dt \, \nabla \cdot \left( \B^n \otimes \frac{\mvec^{*}}{\rho^*} - \frac{\mvec^{*}}{\rho^*} \otimes \B^n \right).
			\label{eqn.sd-rhs_B}
		\end{equation}
		Once the new magnetic field is known, $\mvec^{**}$ is directly calculated from \eqref{eqn.sd-magnetic-rhou}, and the work of magnetic forces can be added to the total energy, hence computing
		\begin{equation}
			E^{**} = E^* - \dt \, \nabla \cdot \left( 2\,m_B^{n+1} \uvec^{n} - \frac{(\uvec^{n}\cdot\B^{n+1})\B^{n+1}}{4\pi} \right).
			\label{eqn.sd-magnetic-E}
		\end{equation}
		
		\item Lastly, the pressure sub-system is considered, which links the momentum equation with the energy equation:
		\begin{subequations}
			\begin{align}
				\mvec^{n+1} & = \mvec^{**} - \dt \, \nabla \cdot \left( p^{n+1}\I \right), \label{eqn.sd-pressure-rhou} \\
				E^{n+1} & = E^{**} - \dt \, \nabla \cdot \left( h^n \mvec^{n+1} \right). \label{eqn.sd-pressure-E}
			\end{align}	
			\label{eqn.sd-pressure}
		\end{subequations}
		Following \cite{BoscheriThomann2024}, formal substitution of \eqref{eqn.sd-pressure-rhou} into \eqref{eqn.sd-pressure-E} yields a scalar wave equation for the unknown total energy:
		\begin{equation}
			E^{n+1} - \dt^2 \, \nabla \cdot \left( h^n \nabla p^{n+1} \right) = E^{**} - (\rho k )^{n} - m_B^{n+1} - \dt \, \nabla \cdot \left( h^n \mvec^{**} \right),
			\label{eqn.sd-waveEq-p}
		\end{equation}
		where $p^{n+1}=p(E^{n+1})$ according to the equation of state \eqref{eq:eos}.	The new momentum is then evaluated by the conservative update \eqref{eqn.sd-pressure-rhou}.
		
	\end{enumerate}

	\subsection{Mesh and discrete approximations}
	\label{sec:mesh-approx}
	
	Let $\Th=\{\omega_c\}_{c=1}^{N_c}$ be a uniform Cartesian mesh of the rectangular
	domain $\Omega \subset\R^2$, and let $L_{\max}$ be its largest side length. Each cell $\omega_c$
	has center $\xvec_c$, side lengths $\Delta x$ and $\Delta y$, and area
	$|\omega_c|=\Delta x\,\Delta y$; we set $\varrho_{\min}=\min(\Delta x,\Delta y)$.
	We denote by $\mathcal F(c)$ the faces of $\omega_c$, by $\mathcal N(c)$ its
	face neighbors, and by $\Ehi$ the set of interior faces. For $f\in\mathcal F(c)$,
	$\nvec_{cf}$ is the outward unit normal. An interior face $f$ separates two cells
	$c^-$ and $c^+$, its normal $\nvec_f$ points from $c^-$ to $c^+$, and for a quantity
	$w$ with traces $w^\mp$ we write
	\begin{equation}
		\jump{w}=w^+-w^-,\qquad \avg{w}=\tfrac12\,(w^-+w^+).
		\label{eq:jump-mean}
	\end{equation}
	The Cartesian unit vectors are $e_x,e_y$, and $\mathbb P^r$ denotes polynomials
	of total degree at most $r$. 
	
	\paragraph{Finite volume (FV) approximation}
	The conservative variables $\Q=(\rho,\mvec,E,\B)$ are approximated by their
	cell averages
	\begin{equation}
		\Q_c(t)=\frac{1}{|\omega_c|}\int_{\omega_c}\Q(\xvec,t)\dd\xvec ,
		\qquad c=1,\dots,N_c,
		\label{eq:cell-average}
	\end{equation}
	that is, by a piecewise constant function $\Q_h$ with $\Q_h|_{\omega_c}=\Q_c$.
	Face values are obtained by a conservative reconstruction: from the averages on
	$\omega_c$ and $\mathcal N(c)$ we build a polynomial $\mathcal R_c(\xvec)$ on
	$\omega_c$ with mean $\Q_c$, and its traces on the faces of $\omega_c$ give the
	fluid states $\Q_f^\mp$ entering the numerical fluxes. In particular, we use a compact CWENO reconstruction \cite{dumbser2017central,cravero2018cweno}, which is briefly summarized in \ref{app:reconstruction-coefficients}. The pressure of the finite
	volume state follows from the equation of state \eqref{eq:eos} evaluated on $\Q_c$:
	\begin{equation}
		p_{\mathrm{FV},c}=(\gamma-1)
		\left(E_c-\frac{|\mvec_c|^2}{2\rho_c}-\frac{|\overline\B_c|^2}{8\pi}\right).
		\label{eq:cell-fv-pressure}
	\end{equation}
	This pressure is used by the reconstruction and its indicators. The fluid
	components of $\Q_c$ are the finite volume unknowns. The magnetic component
	$\overline\B_c$ is not an independent unknown: it is the mean of the
	divergence-free polynomial field introduced next.
	
	\paragraph{Divergence-free Discontinuous Galerkin (DG) approximation of the magnetic field}
	On the reference cell $\widehat\omega=[0,1]^2$, with coordinates $(\xi,\eta)$,
	consider the local solenoidal space
	\begin{equation}
		\widehat V^{\mathrm{div}}_1
		=\bigl\{\bm v\in(\mathbb P^1(\widehat\omega))^{2}
		\;\big|\;\widehat\nabla\cdot\bm v=0\bigr\}.
		\label{eq:complete-local-space}
	\end{equation}
	The divergence maps the six-dimensional space $(\mathbb P^1)^2$ onto $\mathbb P^0$,
	so $\dim\widehat V^{\mathrm{div}}_1=5$. A basis is given by the rotated gradients
	$\widehat\nabla^\perp\psi=(\partial_\eta\psi,-\partial_\xi\psi)^\top$ of
	$\psi\in\{\xi,\eta,\xi^2,\xi\eta,\eta^2\}$:
	\begin{equation}
		\bigl(\widehat{\bm\phi}_1,\ldots,\widehat{\bm\phi}_5\bigr)=
		\begin{pmatrix}
			0 & 1 & 0     & \xi   & 2\eta\\
			-1 & 0 & -2\xi & -\eta & 0
		\end{pmatrix}.
		\label{eq:reference-basis}
	\end{equation}
	Each $\widehat{\bm\phi}_\ell$ is thus divergence-free by construction. This is the
	polynomial kernel underlying $H(\mathrm{div})$ approximation according to 
	\cite{BrezziDouglasMarini1985}. The basis is mapped to the physical cell $\omega_c$ by the contravariant Piola transform of the
	affine map $T_c:\widehat\omega\to\omega_c$, with Jacobian
	$\mathds J_c=\operatorname{diag}(\Delta x,\Delta y)$ and $J_c=\det\mathds J_c$:
	\begin{equation}
		\bm\phi_{c\ell}(\xvec)
		=J_c^{-1}\mathds J_c\,\widehat{\bm\phi}_\ell\bigl(T_c^{-1}(\xvec)\bigr)
		=\begin{pmatrix} 1/\Delta y & 0\\ 0 & 1/\Delta x\end{pmatrix}
		\widehat{\bm\phi}_\ell\bigl(T_c^{-1}(\xvec)\bigr).
		\label{eq:piola}
	\end{equation}
	The magnetic field is thus approximated cellwise by 
	\begin{equation}
		\B_h|_{\omega_c}(\xvec,t)=\sum_{\ell=1}^{5}a_{c\ell}(t)\,\bm\phi_{c\ell}(\xvec),
		\qquad
		\widehat\B_h=
		\begin{pmatrix}
			a_{c2}+a_{c4}\,\xi+2a_{c5}\,\eta\\
			-a_{c1}-2a_{c3}\,\xi-a_{c4}\,\eta
		\end{pmatrix},
		\label{eq:local-field}
	\end{equation}
	where $\widehat\B_h$ is its pull-back to $\widehat\omega$ and $a_{c\ell}$ are the unknown expansion coefficients. Since
	$\widehat\nabla\cdot\widehat\B_h=a_{c4}-a_{c4}=0$ and the Piola transform gives
	$\nabla\cdot\B_h|_{\omega_c}=J_c^{-1}\,\widehat\nabla\cdot\widehat\B_h$, the
	elementwise divergence $\frakD_c(\B_h):=\nabla\cdot(\B_h|_{\omega_c})$ vanishes
	identically:
	\begin{equation}
		\boldsymbol{\frakD}_h(\B_h)
		:=\bigl(\frakD_c(\B_h)\bigr)_{c=1}^{N_c}=\bm 0
		\qquad\text{for every coefficient vector } \bm a_c\in\R^5 .
		\label{eq:gothic-local-kernel}
	\end{equation}
	Local solenoidality is thus built into the space. It holds for any time
	integrator and any stabilization acting within this space. Averaging \eqref{eq:local-field} over the cell gives the magnetic mean that
	enters $\Q_c$:
	\begin{equation}
		\overline\B_c=
		\begin{pmatrix}
			(a_{c2}+a_{c4}/2+a_{c5})/\Delta y\\
			(-a_{c1}-a_{c3}-a_{c4}/2)/\Delta x
		\end{pmatrix},
		\label{eq:magnetic-mean-map}
	\end{equation}
	which is also used in the pressure definition \eqref{eq:cell-fv-pressure}.
	The induction and magnetic-stress face fluxes use the polynomial traces
	$\B_h^\mp$ instead of $\overline\B_c$.
	
	The element mass matrix and its global block-diagonal assembly are
	\begin{equation}
		(\mathds M_c)_{k\ell}=\int_{\omega_c}\bm\phi_{ck}\cdot\bm\phi_{c\ell}\dd\xvec,
		\qquad
		\Mh=\operatorname{diag}_c\mathds M_c.
		\label{eq:mass-matrix}
	\end{equation}
	The cell-averaged magnetic pressure of the polynomial field is
	\begin{equation}
		\overline{m}_{B,c}(\bm a_c)
		=\frac{1}{8\pi|\omega_c|}\int_{\omega_c}|\B_h|^2\dd\xvec
		=\frac{\bm a_c^\top\mathds M_c\,\bm a_c}{8\pi|\omega_c|}.
		\label{eq:cell-magnetic-energy}
	\end{equation}
	The pressure predictor, the enthalpy and the cell diagnostics use this energy,
	\begin{equation}
		\overline p_c=(\gamma-1)
		\left(E_c-\frac{|\mvec_c|^2}{2\rho_c}-\overline{m}_{B,c}(\bm a_c)\right),
		\label{eq:cell-polynomial-pressure}
	\end{equation}
	whereas the finite volume pressure \eqref{eq:cell-fv-pressure} uses
	$|\overline\B_c|^2/(8\pi)$. Splitting $\B_h=\overline\B_c+(\B_h-\overline\B_c)$
	into its mean and a zero-mean part, which are $L^2$-orthogonal, gives
	\begin{equation}
		p_{\mathrm{FV},c}-\overline p_c
		=\frac{\gamma-1}{8\pi|\omega_c|}
		\int_{\omega_c}|\B_h-\overline\B_c|^2\dd\xvec\;\geq\;0 .
		\label{eq:pressure-representation-gap}
	\end{equation}

	\paragraph{Discrete state}
	Collecting $\bm a=(\bm a_c)_c$, the discrete state is
	$\bm U=\bigl((\Q_c)_c,\bm a\bigr)$. Its independent unknowns are therefore the fluid
	averages, from the finite volume part, and the coefficients $\bm a$, from the DG
	part. The magnetic entries of $\Q_c$ are given by \eqref{eq:magnetic-mean-map}.

	\subsection{Convective sub-system}
	The convective sub-system \eqref{eqn.sd-conv} is solved with the aid of an explicit finite volume scheme 
	\begin{equation}
		\Q_c^{*} = \Q_c^{n} -\frac{\dt}{|\omega_c|}\sum_{f\in\mathcal F(c)}|f|
		\sum_q\varpi_q\widehat\F_{f,q}^{\mathrm c},
		\label{eqn.fvscheme}
	\end{equation}
	where the face integral has been approximated by quadrature rules with nodes $\xvec_{fq}$ and normalized weights $\varpi_q$. On $[0,1]$, the three Gauss nodes are $s_q=\tfrac12(1-\sqrt{3/5}),\tfrac12,\tfrac12(1+\sqrt{3/5})$, with weights $\varpi_q=5/18,4/9,5/18$. Volume integrals use their $3\times3$ tensor product. The magnetic field is evaluated directly from its five coefficients at these nodes. The numerical flux function is of Rusanov-type \cite{Rusanov1962} and is given by
	\begin{equation}
		\widehat\F_{f,q}^{\mathrm c} = \frac12\Big(\F^{\mathrm c}(\Q_{f,q}^{-}) + \F^{\mathrm c}(\Q_{f,q}^{+})\Big)\,\nvec_{cf}
		-\frac12\, s_{f,q}^{\max}\,\Big(\Q_{f,q}^{+} - \Q_{f,q}^{-}\Big),
		\label{eqn.rusanov}
	\end{equation}
	where the reconstructed values $\Q_{f,q}^{\pm}$ are obtained from the reconstruction detailed in \ref{app:reconstruction-coefficients} and the numerical dissipation is only proportional to the material speed, hence
	\begin{equation}
		s_{f,q}^{\max} = 2 \max(|\uvec^\pm\cdot\nvec_{cf}|).
		\label{eqn.rusanov_speed}
	\end{equation}

	\subsection{Semi-implicit magnetic and pressure sub-systems}
	\label{sec:semiimplicit}
	
	In the following subsection, we describe the fully discrete counterpart of the implicit steps of the semi-discrete scheme presented in \Cref{sec:sdscheme}. The convective step \eqref{eqn.sd-conv} is solved with the finite volume scheme \eqref{eqn.fvscheme}, which provides the intermediate state $\Q^{*}$ and in particular the new density $\rho^{n+1}=\rho^{*}$. The remaining part of the time step from $t^n$ to $t^{n+1}$ is made of the following stages:
	\begin{enumerate}
		\item a pressure predictor for the velocity entering the induction equation;
		\item a linear implicit solve for the magnetic field, which is the discrete counterpart of \eqref{eqn.sd-waveEqn_B};
		\item a local reconstruction of the magnetic field in troubled cells;
		\item the magnetic update of momentum and energy, according to \eqref{eqn.sd-magnetic-rhou} and \eqref{eqn.sd-magnetic-E};
		\item a linear implicit solve of the pressure sub-system \eqref{eqn.sd-pressure}.
	\end{enumerate}
	
	\paragraph{Transport and pressure predictors} We make use of two centered operators, namely the gradient of a cell quantity $z$ and the divergence of a vector or tensor flux $\mathds F$, which are given by
	\begin{align}
		(\mathds G_h z)_c&=\frac1{|\omega_c|} \sum_{f\in\mathcal F(c)}|f|\avg{z}_f\nvec_{cf}, \\
		(\mathcal D_h\mathds F)_c
		& =\frac1{|\omega_c|}\sum_{f\in\mathcal F(c)}|f|
		\sum_q\varpi_q\avg{\mathds F}_{f,q}\nvec_{cf}.
		\label{eq:centered-divergence}
	\end{align}
Let us notice that, at this stage, the magnetic components of $\Q^{*}$ are not evolved by the convective step, but they are given by the coefficients $\bm a^{n}$. 
	
\begin{figure}[!htb]
		\centering
		\begin{tikzpicture}[
  font=\small,
  cell/.style={draw=black!75,line width=.8pt},
  subedge/.style={draw=black!35,line width=.45pt},
  normal/.style={-Latex,draw=teal!75!black,line width=.8pt},
  lab/.style={font=\footnotesize,inner sep=2pt}
]
\pgfmathsetmacro{\gaussleft}{(1-sqrt(3/5))/2}
\pgfmathsetmacro{\gaussright}{(1+sqrt(3/5))/2}
\pgfmathsetmacro{\subleft}{5/18}
\pgfmathsetmacro{\subright}{13/18}
\begin{scope}[x=4.8cm,y=4.8cm]
  \node[anchor=south] at (.5,1.10) {(a) Local stress balances};
  \fill[blue!2] (0,0) rectangle (1,1);
  \fill[teal!12] (\subleft,\subleft) rectangle (\subright,\subright);
  \foreach \s in {\subleft,\subright} {
    \draw[subedge] (\s,0) -- (\s,1);
    \draw[subedge] (0,\s) -- (1,\s);
  }
  \draw[cell] (0,0) rectangle (1,1);
  \foreach \s in {0,\subleft,\subright,1} {
    \foreach \t in {0,\subleft,\subright,1} {
      \fill[blue!75!black] (\s,\t) circle[radius=1.7pt];
    }
  }
  \foreach \s in {\gaussleft,.5,\gaussright} {
    \foreach \t in {\gaussleft,.5,\gaussright} {
      \fill[orange!85!black] (\s,\t) circle[radius=1.8pt];
    }
  }
  \draw[normal] (.5,\subleft) -- (.5,.16);
  \draw[normal] (\subright,.5) -- (.84,.5);
  \draw[normal] (.5,\subright) -- (.5,.84);
  \draw[normal] (\subleft,.5) -- (.16,.5);
  \node[lab] at (.5,.63) {$\widehat\omega^{\mathrm{sg}}_{22}$};
  \node[lab,anchor=north] at (0,-.02) {$0$};
  \node[lab,anchor=north] at (\subleft,-.02) {$\frac5{18}$};
  \node[lab,anchor=north] at (\subright,-.02) {$\frac{13}{18}$};
  \node[lab,anchor=north] at (1,-.02) {$1$};
  \node[lab,anchor=east] at (-.01,1) {$1$};
  \node[lab,anchor=west] at (1.035,0) {$\xi$};
  \node[lab,anchor=south] at (0,1.035) {$\eta$};
  \node[lab] at (.5,-.19) {
    \textcolor{blue!75!black}{$\bullet$}\; stress vertices
    \quad\textcolor{orange!85!black}{$\bullet$}\; Gauss nodes};
  \node[lab] at (.5,-.30) {$\overline{\bm g}_{c,ij}\;\longmapsto\;
    \bm g_B(\bm Z;\xvec^{\mathrm G}_{c,ij})$};
\end{scope}
\begin{scope}[x=3cm,y=3cm,xshift=6.1cm,yshift=.9cm]
  \node[anchor=south] at (1,1.46) {(b) Interior-face coupling};
  \fill[teal!12] (\subright,\subleft) rectangle (1,\subright);
  \fill[teal!12] (1,\subleft) rectangle ({1+\subleft},\subright);
  \foreach \a in {0,1} {
    \foreach \s in {\subleft,\subright} {
      \draw[subedge] ({\a+\s},0) -- ({\a+\s},1);
      \draw[subedge] (\a,\s) -- ({\a+1},\s);
    }
    \draw[cell] (\a,0) rectangle ({\a+1},1);
  }
  \draw[blue!70!black,line width=1pt] (1,0) -- (1,1);
  \node at (.5,1.13) {$\omega_c$};
  \node at (1.5,1.13) {$\omega_d$};
  \foreach \s in {\gaussleft,.5,\gaussright} {
    \fill[orange!85!black] (1,\s) circle[radius=2pt];
  }
  \draw[normal] (.78,.5) -- (.975,.5);
  \draw[normal] (1.22,.5) -- (1.025,.5);
  \node[lab,anchor=east] at (.97,.67) {$\overline{\bm g}_{c,3q}$};
  \node[lab,anchor=west] at (1.03,.67) {$\overline{\bm g}_{d,1q}$};
  \node[lab,anchor=west] at (1.06,.32) {$\xvec_{fq}$};
  \draw[black!50,line width=.4pt] (1.06,.35) -- (1.02,.47);
  \node[lab] at (1,-.33) {Common force at the face quadrature node};
  \node[lab] at (1,-.59) {$\displaystyle
    \bm g_B(\bm Z;\xvec_{fq})=
    \frac12\bigl(\overline{\bm g}_{c,3q}+\overline{\bm g}_{d,1q}\bigr)$};
\end{scope}
\end{tikzpicture}
		\caption{Subcell evaluation of the magnetic force. The stress is sampled at the sixteen vertices of the local partition. Boundary balances supply nine subcell averages, assigned to the volume quadrature nodes. At an interior parent-cell face, the two adjacent subcell values are averaged. All magnetic samples depend on the same five local divergence-free basis functions.}
		\label{fig:subgrid-geometry}
	\end{figure}
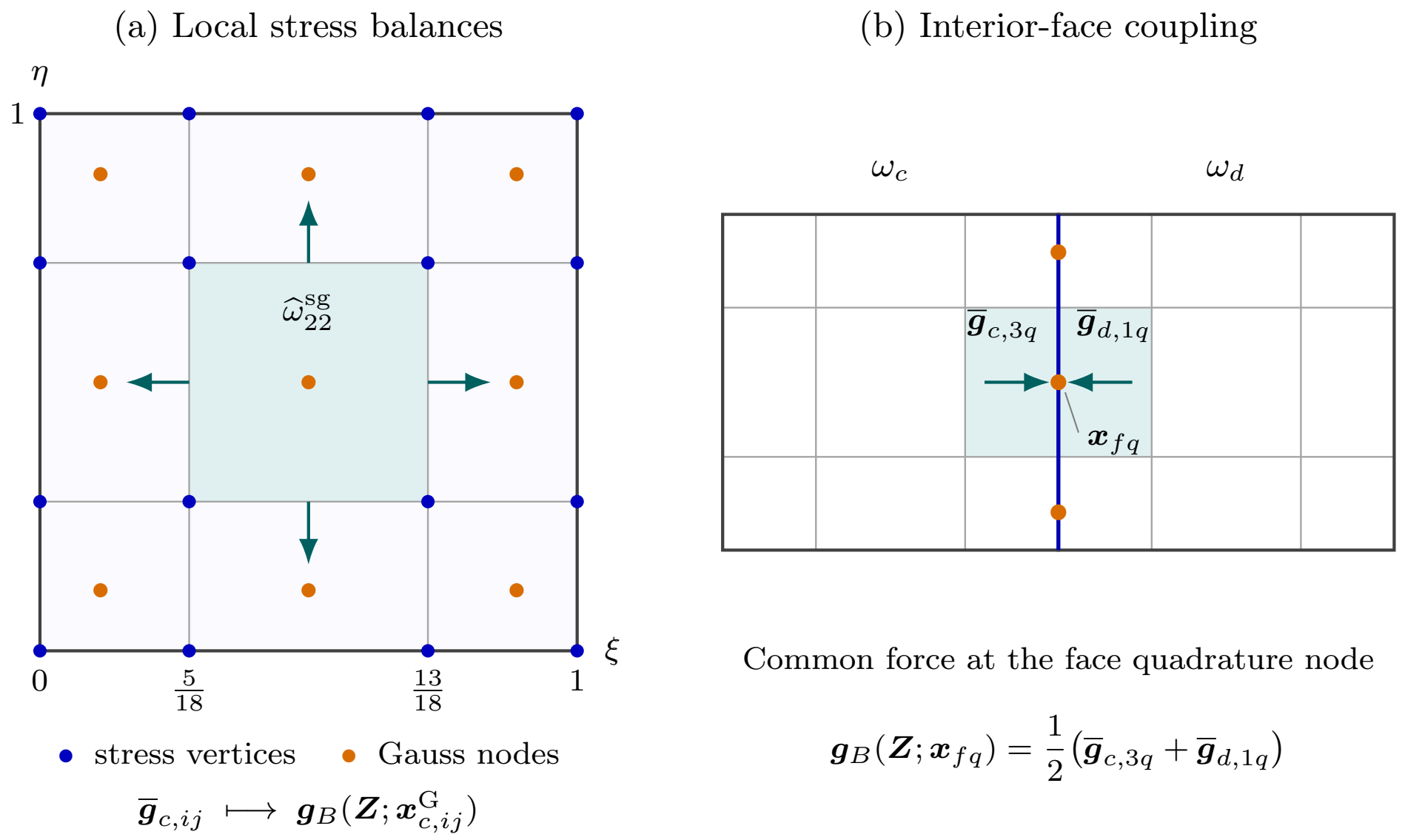

	\paragraph{Subcell evaluation of the magnetic force}
	As in \eqref{eqn.sd-waveEqn_B}, inserting the momentum equation into the induction equation gives rise to the divergence of a magnetic stress, which is linearized by taking one of its two magnetic factors at time $t^n$. For a trial field $\bm Z$ in the five-dimensional solenoidal space, the linearized stress tensor is
	\begin{equation}
		\mathds T_B(\B^{n},\bm Z)
		=\frac{\B^{n}\cdot\bm Z}{8\pi}\I
		-\frac{\bm Z\otimes\B^{n}}{4\pi}.
		\label{eq:bilinear-magnetic-stress}
	\end{equation}
	We compute its divergence on a $3\times3$ subgrid built in each cell based on the tensor product of one-dimensional Gauss-Lobatto points, see \Cref{fig:subgrid-geometry}. With $z_0=0$, $z_1=5/18$, $z_2=13/18$ and $z_3=1$, the subcells are defined by
	\begin{equation}
		\widehat\omega^{\mathrm{sg}}_{ij}
		=[z_{i-1},z_i]\times[z_{j-1},z_j],\qquad
		\omega^{\mathrm{sg}}_{c,ij}=T_c(\widehat\omega^{\mathrm{sg}}_{ij}),
		\quad i,j\in\llbracket 1,3\rrbracket.
		\label{eq:subgrid-partition}
	\end{equation}
	The subcell widths coincide with the Gauss weights, i.e. $z_i-z_{i-1}=\varpi_i$, thus $|\omega^{\mathrm{sg}}_{c,ij}|=|\omega_c|\varpi_i\varpi_j$, and each subcell contains exactly one volume Gauss node $\xvec^{\mathrm G}_{c,ij}=T_c(s_i,s_j)$. The two magnetic fields $\B_h^{n}$ and $\bm Z$ are evaluated at the sixteen subgrid vertices $\xvec^{\mathrm v}_{c,ij}=T_c(z_i,z_j)$ from their five basis coefficients, e.g.
	\begin{equation}
		\bm Z(\xvec^{\mathrm v}_{c,ij})
		=\sum_{\ell=1}^5 a^{Z}_{c\ell}\bm\phi_{c\ell}(\xvec^{\mathrm v}_{c,ij}),
		\quad i,j\in\llbracket 0,3\rrbracket.
		\label{eq:subgrid-evaluation}
	\end{equation}
	The force in each subcell is then obtained by applying the trapezoidal rule along its four edges $e\in\mathcal F^{\mathrm{sg}}_{c,ij}$, with endpoints $\xvec_{e,1},\xvec_{e,2}$ and outward normal $\nvec_e$:
	\begin{equation}
		\overline{\bm g}_{c,ij}(\bm Z)
		=\frac{1}{|\omega^{\mathrm{sg}}_{c,ij}|}
		\sum_{e\in\mathcal F^{\mathrm{sg}}_{c,ij}}\frac{|e|}{2}
		\sum_{r=1}^2
		\mathds T_B(\B_h^{n},\bm Z)(\xvec_{e,r})\nvec_e.
		\label{eq:subgrid-stress-balance}
	\end{equation}
	Since the stress is quadratic in each cell, the trapezoidal rule is not exact on a single edge. However, opposite edges share the same tangential quadratic coefficient, so the errors cancel and \eqref{eq:subgrid-stress-balance} gives exactly the subcell average of $\nabla\cdot\mathds T_B$. This average is finally assigned to the Gauss node of the subcell:
	\begin{equation}
		\bm g_B(\bm Z;\xvec^{\mathrm G}_{c,ij})
		=\overline{\bm g}_{c,ij}(\bm Z).
		\label{eq:subgrid-volume-force}
	\end{equation}
	We remark that \eqref{eq:subgrid-volume-force} is an approximation, because the Gauss node does not coincide in general with the subcell centroid.
	
	On the cell boundary faces, the force is taken as the mean of the two adjacent subcell values, see Figure \ref{fig:subgrid-geometry}. At quadrature node $q$, this gives $\tfrac12(\overline{\bm g}_{c,3q}+\overline{\bm g}_{d,1q})$ on a vertical face between a left cell $c$ and a right cell $d$, and $\tfrac12(\overline{\bm g}_{c,q3}+\overline{\bm g}_{d,q1})$ on a horizontal face between a lower cell $c$ and an upper cell $d$. At physical boundaries the one-sided value is used. This common force then enters the two induction traces, each with its own density and magnetic field, and the induction flux is averaged as in \eqref{eq:weak-induction}.
	
	Let us notice that the subgrid adds neither new unknowns nor new equations: it only serves to evaluate the force inside the linear magnetic operator. Moreover, it is compatible with the divergence-free constraint. Indeed, the normal trace of $\B_h$ is affine on each edge, so the trapezoidal rule is exact and
	\begin{equation}
		\frac{1}{|\omega^{\mathrm{sg}}_{c,ij}|}
		\sum_{e\in\mathcal F^{\mathrm{sg}}_{c,ij}}\frac{|e|}{2}
		\sum_{r=1}^2\B_h(\xvec_{e,r})\cdot\nvec_e
		=\frac{1}{|\omega^{\mathrm{sg}}_{c,ij}|}
		\int_{\omega^{\mathrm{sg}}_{c,ij}}\frakD_c(\B_h)\dd\xvec=0.
		\label{eq:subgrid-divergence}
	\end{equation}
	The fluxes on internal subcell edges cancel since they are traces of the same polynomial, while across cell faces the cancellation is ensured by the normal-trace constraints of \Cref{sec:global-div}. In other words, the basis gives the local divergence-free property, and the face constraints give global compatibility.
	
	\paragraph{Induction in the solenoidal basis}
	Let $\mathds F_B(\uvec,\B)=\B\otimes\uvec-\uvec\otimes\B$ denote the induction flux. Following the semi-discrete scheme, the velocity $\uvec^{**}$ given by \eqref{eqn.sd-magnetic-rhou} is approximated by
	\begin{eqnarray}
		\uvec^{**} = \frac{(\mvec^{*}-\dt\,\bm g_B(\B_h^{n+1}))}{\rho^{*}},
	\end{eqnarray}
	and the other magnetic factor in \eqref{eqn.sd-magnetic-B} is taken at time $t^n$. This splits the flux into an explicit and an implicit part:
	\begin{equation}
		\mathds F_B^{n}
		=\mathds F_B(\mvec^{*}/\rho^*,\B_h^{n}),
		\qquad
		\mathds F_B^{n+1}
		=-\dt\,\mathds F_B\bigl(\bm g_B(\B_h^{n+1})/\rho^{*},\B_h^{n}\bigr).
		\label{eq:linearized-induction-flux}
	\end{equation}
	As for the convective sub-system, the explicit part is discretized with a Rusanov-type numerical flux \cite{Rusanov1962}, which at the quadrature point $\xvec_{fq}$ writes
	\begin{align}
		\widehat{\bm f}_{B,f,q}^{n}
		&=\frac12\Big(\mathds F_B(\uvec^{*,-}_{f,q},\B_{h,f,q}^{n,-})+\mathds F_B(\uvec^{*,+}_{f,q},\B_{h,f,q}^{n,+})\Big)\nvec_{cf} \nonumber\\
		&-\frac12\Big[(1-\sigma)\,s_{f,q}^{\max}\Big(\B_{h,f,q}^{n,+}-\B_{h,f,q}^{n,-}\Big)
		+\sigma\, s_{f,q}^{\mathrm f}\Big((\mathcal R_h^B\overline\B)^{+}_{f,q}-(\mathcal R_h^B\overline\B)^{-}_{f,q}\Big)\Big],
		\label{eq:explicit-induction-flux}
	\end{align}
	where $\B_{h,f,q}^{n,\pm}$ are the polynomial traces of the divergence-free magnetic field and $\widetilde\uvec^{*,\pm}_{f,q}$ are the reconstructed velocities at the two sides of the face. The parameter $\sigma\in\{0,1\}$ is a global switch: $\sigma=1$ if at least one cell is marked by the troubled-cell indicator \eqref{eq:troubled-cell-indicator}, and $\sigma=0$ otherwise. For $\sigma=0$, the numerical dissipation is proportional to the material speed \eqref{eqn.rusanov_speed} and acts on the polynomial traces. For $\sigma=1$, it is proportional to the fast magnetosonic speed \eqref{eq:fast-magnetosonic-speed} of the explicit finite volume traces,
	\begin{equation}
		s_{f,q}^{\mathrm f}=\max_\pm\left(|\uvec^{\pm}_{f,q}\cdot\nvec_{cf}|+c_f(\Q^{\pm}_{f,q};\nvec_{cf})\right),
		\label{eq:explicit-induction-fast-speed}
	\end{equation}
	and it acts on the magnetic components $(\mathcal R_h^B\overline\B)^{\pm}_{f,q}$ of the finite volume reconstruction of \ref{app:reconstruction-coefficients}. Furthermore, for $\sigma=0$ the explicit predictor uses an unlimited centered reconstruction, while for $\sigma=1$ the limited reconstruction with the positivity scaling detailed in \ref{app:reconstruction-coefficients} is applied on the whole mesh. The subcell force is the same in both cases. We remark that the finite volume reconstruction of the magnetic field only enters the dissipative term, whereas the central part of the flux always uses the locally solenoidal field $\B_h$.
	
	We adopt a discontinuous Galerkin scheme for solving the magnetic sub-system. For each test function $\bm\phi_{k}$, the discrete form of \eqref{eqn.sd-waveEqn_B} reads
	\begin{equation}
		\begin{split}
			\int_{\omega_c}\bm\phi_{k}\cdot
			\frac{\B_h^{n+1}-\B_h^{n}}{\dt}\dd\xvec
			-\int_{\omega_c}\nabla\bm\phi_{k}:
			\mathds F_B^{n+1}\dd\xvec
			+\int_{\partial\omega_c}\bm\phi_{k}\cdot
			\avg{\mathds F_B^{n+1}}
			\nvec_{cf}\dd s
			=\mathcal R_{B,ck}^{n},
			\label{eq:weak-induction}
		\end{split}
	\end{equation}
	where $\mathcal R_{B,ck}^{n}$ is the weak residual of $-\nabla\cdot\mathds F_B^{n}$ computed with the numerical flux \eqref{eq:explicit-induction-flux}, and it corresponds to the right hand side \eqref{eqn.sd-rhs_B}. After assembly, and including the stabilization term introduced below, we obtain the linear system
	\begin{equation}
		(\Mh+\dt\,\Ah_B+\dt\,\Mh\mathds D_\nu)\, \bm a^{n+1}
		=\bm r_B,\qquad
		\bm r_B=\Mh\bm a^{n}+\dt\,\bm{\mathcal R}_B^{n},
		\label{eq:magnetic-linear-system}
	\end{equation}
	where $\Ah_B$ collects the two implicit flux terms of \eqref{eq:weak-induction} and $\mathds D_\nu=\operatorname{diag}_c\mathds D_{\nu,c}$.	To stabilize the magnetic solve, we add an implicit damping of the zero-mean part of the magnetic field in each cell. With the coefficient ordering of \eqref{eq:reference-basis}, this term is defined by
	\begin{equation}
		\mathds P_0\bm a_c=
		\begin{pmatrix}
			-a_{c3}-a_{c4}/2\\
			-a_{c4}/2-a_{c5}\\
			a_{c3}\\a_{c4}\\a_{c5}
		\end{pmatrix},
		\qquad \mathds D_{\nu,c}=\nu_c\mathds P_0 .
		\label{eq:zero-mean-modal-operator}
	\end{equation}
	The rate $\nu_c\geq0$ is evaluated at time $t^n$ and measures how far the directional CWENO weights are from their linear values (see \ref{app:reconstruction-coefficients}). From \eqref{eq:magnetic-mean-map}, the field associated with $\mathds P_0\bm a_c$ is exactly $\B_h-\overline\B_c$. Thus, the damping neither modifies the cell mean nor introduces a divergence. Furthermore, the constant and zero-mean parts are orthogonal, hence the term is dissipative:
	\begin{equation}
		\bm a_c^\top\mathds M_c\mathds D_{\nu,c}\bm a_c
		=\nu_c\|\B_h-\overline\B_c\|_{L^2(\omega_c)}^2\geq0 .
		\label{eq:modal-dissipation}
	\end{equation}
	
	\paragraph{Locally solenoidal troubled-cell reconstruction}
	\label{sec:postsolve-magnetic-reconstruction}
	After the magnetic solve, the field is reconstructed in the cells marked by \eqref{eq:troubled-cell-indicator} and in their face neighbors. Let $\bm a_c^{n+1}$ be the coefficients given by the linear system \eqref{eq:magnetic-linear-system} and $\overline\B_c$ the associated cell mean. The reconstruction is carried out in three steps.
	\begin{enumerate}
		\item Compute componentwise slopes from the neighboring cell means following \cite{vanLeer1977},
		\begin{equation}
			\bm s_{d,c}=\operatorname{minmod}\!\left(
			\frac{\overline\B_{c+e_d}-\overline\B_{c-e_d}}{2},
			2(\overline\B_c-\overline\B_{c-e_d}),
			2(\overline\B_{c+e_d}-\overline\B_c)\right),
			\qquad d\in\{x,y\},
			\label{eq:limited-magnetic-slopes}
		\end{equation}
		which define the target field
		\begin{equation}
			\B_c^{\mathrm{tar}}(\xi,\eta)=\overline\B_c
			+\bm s_{x,c}(\xi-\tfrac12)+\bm s_{y,c}(\eta-\tfrac12).
			\label{eq:magnetic-target-polynomial}
		\end{equation}
		\item Since componentwise slopes are in general not solenoidal, project $\B_c^{\mathrm{tar}}$ in the $L^2$ norm onto $V^{\mathrm{div}}_1(\omega_c)$, and correct the two constant coefficients so that the projection $\bm a_c^{\mathrm{tar}}$ has mean $\overline\B_c$.
		\item Blend the target with the constant field $\bm a_c^0$ of minimum energy and mean $\overline\B_c$:
		\begin{equation}
			\bm a_c^{\mathrm{lim}}=\bm a_c^0+
			\theta_c(\bm a_c^{\mathrm{tar}}-\bm a_c^0),
			\qquad 0\leq\theta_c\leq1.
			\label{eq:magnetic-mode-reconstruction}
		\end{equation}
		We set $\theta_c=1$, unless the target has more magnetic energy than the field $\bm a_c^{n+1}$ returned by the linear solve. In that case
		\begin{equation}
			\theta_c=\min\!\left\{1,
			\left[\frac{\varepsilon_{B,c}(\bm a_c^{n+1})-
				\varepsilon_{B,c}(\bm a_c^0)}
			{\varepsilon_{B,c}(\bm a_c^{\mathrm{tar}})-
				\varepsilon_{B,c}(\bm a_c^0)}\right]_+^{1/2}\right\}.
			\label{eq:magnetic-reconstruction-energy-cap}
		\end{equation}
	\end{enumerate}
	In the marked cells, $\bm a_c^{n+1}$ is then replaced by $\bm a_c^{\mathrm{lim}}$. Both $\bm a_c^0$ and $\bm a_c^{\mathrm{tar}}$ belong to $V^{\mathrm{div}}_1(\omega_c)$ and have the same mean, therefore the reconstruction is divergence-free and preserves the magnetic mean. Moreover, the magnetic energy is quadratic in $\theta_c$, so that \eqref{eq:magnetic-reconstruction-energy-cap} gives the largest factor which does not increase the energy of the cell. However, we remark that the projection used to restore the continuity of the normal traces, introduced in Section \ref{sec:global-div}, may violate this cellwise energy bound.
	
	\paragraph{Magnetic momentum and energy update}
	Once the new magnetic field $\B_h^{n+1}$ is known, the magnetic contributions to momentum and energy in \eqref{eqn.sd-magnetic-rhou} and \eqref{eqn.sd-magnetic-E} can be computed explicitly. In the fully discrete scheme, we also move here the explicit part of the pressure. Indeed, writing the new pressure as $p^{n+1}=(\gamma-1)(E^{n+1}-(\rho k)^n-m_B^{n+1})$, only the term $(\gamma-1)E^{n+1}\I$ is left to the pressure sub-system, while the remaining pressure and magnetic terms give the stress
	\begin{equation}
		\mathds T_m
		=\left[-(\gamma-1)\rho^{n}k^{n}
		+(2-\gamma)m_B^{n+1}\right]\I
		-\frac{\B_h^{n+1}\otimes\B_h^{n+1}}{4\pi},
		\label{eq:magnetic-stress}
	\end{equation}
	where $\B_h^{n+1}$ already includes the troubled-cell reconstruction. In unmarked cells, the magnetic energy flux is
	\begin{equation}
		\bm f_{m,c}^E=\frac{|\overline\B_c^{n+1}|^2}{4\pi}
		\uvec_c^{n}
		-\frac{(\uvec_c^{n}\cdot\overline\B_c^{n+1})
			\overline\B_c^{n+1}}{4\pi},
		\qquad
		\widehat{\bm f}_{m,f}^E=\avg{\bm f_{m,c}^E},
		\label{eq:compatible-magnetic-energy-flux}
	\end{equation}
	while at a detected discontinuity the same expression is evaluated on each side of the face with the limited FV velocity and the polynomial magnetic trace, before taking the central flux. 
	The intermediate momentum and energy are then
	\begin{equation}
		\mvec^{**}=\mvec^{*}-\dt\,\mathcal D_h\mathds T_m,
		\qquad
		E^{**}=E^{*}-\dt\,\mathcal D_h\widehat{\bm f}_m^E ,
		\label{eq:magnetic-update}
	\end{equation}
	while the density remains $\rho^{n+1}=\rho^{*}$. 
	
	\paragraph{Scalar energy equation}
	The last stage is the implicit solve of the pressure sub-system \eqref{eqn.sd-pressure}. 
	Taking into account the magnetic stress update \eqref{eq:magnetic-stress} into $\mvec^{**}$, the momentum equation \eqref{eqn.sd-pressure-rhou} is left with
	\begin{equation}
		\mvec^{n+1}
		=\mvec^{**}-\dt\,(\gamma-1)\mathds G_h E^{n+1}.
		\label{eq:pressure-momentum}
	\end{equation}
	Inserting \eqref{eq:pressure-momentum} into the energy equation \eqref{eqn.sd-pressure-E}, and using on periodic meshes the same centered operators as above, we define
	\begin{equation}
		\mathcal L_H z=\mathcal D_h\!\left(h^{n}\mathds G_hz\right),
		\label{eq:enthalpy-laplacian}
	\end{equation}
	and we obtain the discrete counterpart of \eqref{eqn.sd-waveEq-p}:
	\begin{equation}
		\mathds P_H(\dt)E^{n+1}
		=E^{**}-\dt\,\mathcal D_h(h^{n}\mvec^{**}),
		\qquad
		\mathds P_H(\dt)=\I-(\gamma-1)\dt^2\mathcal L_H ,
		\label{eq:pressure-system}
	\end{equation}
	where the enthalpy flux on the right hand side is $\avg{h^{n}\mvec^{**}}$. The cell enthalpy is given by $h_c^{n}=e_c^{n}+\overline p_c^{\,n}/\rho_c^{n}=\gamma\,\overline p_c^{\,n}/\bigl((\gamma-1)\rho_c^{n}\bigr)$. In the volume-weighted inner product $\langle z,w\rangle_V=\sum_{\omega_c\in\Th}|\omega_c|z_c\cdot w_c$, the periodic operators satisfy $\mathcal D_h=-\mathds G_h^*$. Thus, for positive enthalpy, the system is symmetric positive definite:
	\begin{equation}
		\langle z,\mathds P_H z\rangle_V
		=\sum_{\omega_c\in\Th}|\omega_c|z_c^2
		+(\gamma-1)\dt^2
		\langle\mathds G_hz,h^{n}\mathds G_hz\rangle_V>0
		\quad(z\ne0).
		\label{eq:pressure-spd-identity}
	\end{equation}
	This guarantees solvability, but not positivity of the pressure. With physical boundaries, the scalar equation uses instead a compact symmetric two-point closure, and the prescribed states enter the enthalpy flux on the right hand side. Finally, the momentum $\mvec^{n+1}$ is updated with \eqref{eq:pressure-momentum} and the cell pressure $\overline p^{\,n+1}$ is recovered from \eqref{eq:cell-polynomial-pressure}.
	
	%
	%
	
	\section{Globally divergence-free magnetic field}
	\label{sec:global-div}
	
	The basis \eqref{eq:reference-basis} guarantees that the magnetic field is divergence-free inside each cell. In this section we show how to make it divergence-free on the whole domain, which only requires the continuity of the normal component across the cell faces.
	
	For a piecewise polynomial field, the divergence in the sense of distributions writes
	\begin{equation}
		\langle\nabla\cdot\B_h,\psi\rangle
		=\sum_{\omega_c\in\Th}\int_{\omega_c}(\nabla\cdot\B_h)\psi\dd\xvec
		+\sum_{f\in\Ehi}\int_f\jump{\B_h\cdot\nvec_f}\psi\dd s,
		\qquad\psi\in C_c^\infty(\Omega).
		\label{eq:distributional-div}
	\end{equation}
	The volume term vanishes by construction, see \eqref{eq:gothic-local-kernel}. Therefore, the field is globally divergence-free if its normal component is continuous across all interior faces, while jumps of the tangential component are still allowed.
	
	Since $\B_h\cdot\nvec_f\in\mathbb P^1(f)$, two moments are sufficient on each interior face. With the face parametrized by $s\in[0,1]$, and $q_0(s)=1$, $q_1(s)=s-1/2$, we define the interface involutions
	\begin{equation}
		\frakJ_{f,r}(\B_h)=
		\int_f\jump{\B_h\cdot\nvec_f}q_r\dd s,\qquad r=0,1 .
		\label{eq:gothic-face-involution}
	\end{equation}
	Collecting all the moments gives $\boldsymbol{\frakJ}_h(\B_h)=\Ch\bm a$, where the matrix $\Ch$ couples the coefficients of neighboring cells; its blocks are reported in \ref{app:constraint-blocks}. Let $V_h^{\mathrm{comp}}$ be the space of piecewise $V^{\mathrm{div}}_1(\omega_c)$ fields that satisfy $\Ch\bm a=0$. The complete discrete involution then reads
	\begin{equation}
		\boldsymbol{\frakL}_h(\B_h)=
		\begin{pmatrix}\boldsymbol{\frakD}_h(\B_h)\\
			\boldsymbol{\frakJ}_h(\B_h)\end{pmatrix}=\bm0 .
		\label{eq:complete-gothic-involution}
	\end{equation}
	If $\Ch\bm a=0$, the normal jump on each face is a linear function orthogonal to the whole $\mathbb P^1(f)$, including itself, thus it vanishes. Together with \eqref{eq:gothic-local-kernel}, relation \eqref{eq:distributional-div} gives a globally divergence-free field:
	\begin{equation}
		\B_h\in H(\operatorname{div};\Omega),
		\qquad \nabla\cdot\B_h=0.
		\label{eq:global-solenoidal-result}
	\end{equation}
	
	\subsection{Projection onto the compatible space}
	The magnetic solve and the troubled-cell reconstruction of \Cref{sec:semiimplicit} do not enforce the face constraints. Therefore, at the end of each time step, the coefficients $\bm a^{\mathrm{pre}}$ are replaced by their projection onto $V_h^{\mathrm{comp}}$ in the mass-matrix norm:
	\begin{equation}
		\bm a=\arg\min_{\Ch\bm z=0}
		\frac12(\bm z-\bm a^{\mathrm{pre}})^\top
		\Mh(\bm z-\bm a^{\mathrm{pre}}).
		\label{eq:constrained-projection}
	\end{equation}
	With the Schur complement $\mathds S=\Ch\Mh^{-1}\Ch^\top$, the projection is computed as
	\begin{equation}
		\mathds S\bm\lambda=\Ch\bm a^{\mathrm{pre}},
		\qquad
		\bm a=\bm a^{\mathrm{pre}}-\Mh^{-1}\Ch^\top\bm\lambda ,
		\label{eq:schur-projection}
	\end{equation}
	and we denote this linear projection by $\mathds P_B$. Only interior faces, including periodic ones, enter $\Ch$, thus no normal boundary condition is imposed on physical boundaries.
	
	\paragraph{Properties of the projection}
	Let us notice that the minimization problem \eqref{eq:constrained-projection} is strictly convex, so the projected field is unique. If some constraints are redundant, $\bm\lambda$ is not unique, but any two solutions differ by an element of $\ker\Ch^\top$ and give the same field. Furthermore, the projection is orthogonal in the mass-matrix norm, hence
	\begin{equation*}
		\|\bm a^{\mathrm{pre}}\|_{\Mh}^2=\|\bm a\|_{\Mh}^2+\|\bm a-\bm a^{\mathrm{pre}}\|_{\Mh}^2,
	\end{equation*}
	which means that the magnetic energy can only decrease. Since constant fields belong to $\ker\Ch$, the correction is orthogonal to them, and the integral of each magnetic component is preserved:
	\begin{equation}
		\int_\Omega\B_h\dd\xvec
		=\int_\Omega\B_h^{\mathrm{pre}}\dd\xvec.
		\label{eq:projection-global-mean}
	\end{equation}
	We remark that this only holds globally: the cell means and the normal trace on physical boundaries can be in general modified.
	
	
	\subsection{Energy correction}
	The projection modifies both the magnetic means and the magnetic energy \eqref{eq:cell-magnetic-energy}. Let $\delta_c=\overline{m}_{B,c}(\bm a_c^{\mathrm{pre}}) -\overline{m}_{B,c}(\bm a_c)$ be the energy change in cell $\omega_c$. Thanks to \eqref{eq:constrained-projection}, the global loss $\Delta \overline{m}_B=\sum_{\omega_c\in\Th}|\omega_c|\delta_c$ is nonnegative. The total energy is then corrected as
	\begin{equation}
		E_c=E_c^{\mathrm{pre}}-\delta_c
		+w_c\Delta \overline{m}_B,
		\qquad
		w_c=\frac{[\delta_c]_+}
		{\sum_v|\omega_v|[\delta_v]_+},
		\label{eq:projection-energy-redistribution}
	\end{equation}
	where the last term is set to zero if the denominator vanishes. The term $-\delta_c$ keeps the pressure $\overline p_c$ unchanged for fixed density and momentum, while the last term redistributes the magnetic energy lost by the projection to the cells where it was lost, so that the total energy is conserved. Consequently, the pressure satisfies
	\begin{equation}
		\overline p_c^{\mathrm{post}}
		=\overline p_c^{\mathrm{pre}}
		+(\gamma-1)w_c\Delta \overline{m}_B,
		\label{eq:projection-pressure-balance}
	\end{equation}
	thus, for positive initial pressure, the correction does not decrease the cell entropy $\log(\overline p_c/\rho_c^\gamma)$ at fixed density and momentum. We remark that this only concerns the projection: an entropy inequality for the whole scheme would require compatible energy transfers as in \cite{BustoDumbser2023}.

	Since the finite volume fluxes cancel across interior faces, mass and momentum are conserved up to boundary fluxes. Thanks to the redistribution in \eqref{eq:projection-energy-redistribution}, the projection does not modify the total energy, which is therefore conserved as well.
	
	\section{Time integration}
	\label{sec:time}
	
	\paragraph{First order scheme}
	Let $\Phi_{\dt}$ denote one complete time step of the first order semi-implicit scheme described in the previous sections, which is made of: the explicit convective step, the magnetic solve with the troubled-cell reconstruction, the magnetic update of momentum and energy, the pressure solve, and finally the projection onto the divergence-free space with its energy correction. The first order scheme simply writes
	\begin{equation}
		\Q^{n+1}=\Phi_{\dt}(\Q^n).
		\label{eq:first-order-time}
	\end{equation}
	
	\paragraph{Second order extrapolation}
	To achieve second order of accuracy in time, we apply a Richardson extrapolation to the complete map $\Phi$. We perform one step of size $\dt$ and two steps of size $\dt/2$, and then combine the two solutions:
\begin{subequations}\label{eq:richardson-composition}
	\begin{align}
	\Q_F&=\Phi_{\dt}(\Q^n), \\
	\Q_H&=\Phi_{\dt/2}(\Q^n), \\
	\Q_{HH}&=\Phi_{\dt/2}(\Q_H), \\
	\Q^{n+1}&=2\Q_{HH}-\Q_F.
	\label{eq:richardson-final-projection}
	\end{align}
\end{subequations}
Let us notice that both $\Q_F$ and $\Q_{HH}$ already satisfy the local and the interface involutions, and so does their linear combination $\Q^{n+1}$. 
	
	To verify the order of accuracy, let us consider a smooth solution, for which the limiters are not active, and assume that one step of the scheme can be expanded as $\Phi_\tau(\Q)=\Q+\tau\mathcal L_h(\Q)+\tau^2\mathcal E_h(\Q)+O(\tau^3)$. Expanding the two half steps yields
	\begin{equation}
		2\Phi_{\dt/2}(\Phi_{\dt/2}(\Q))-\Phi_{\dt}(\Q)
		=\Q+\dt\,\mathcal L_h(\Q)
		+\frac{\dt^2}{2}\mathcal L_h'(\Q)\mathcal L_h(\Q)
		+O(\dt^3),
		\label{eq:richardson-consistency}
	\end{equation}
	which is the Taylor expansion of the exact solution up to second order. Thus, the second order error of the first order map cancels, including the contributions of the projection and of the energy correction. A numerical convergence study in time is reported in \ref{app:coupled-time-order}. We remark that the Riemann problems and the blast wave presented in \Cref{sec:numerics} are computed with the first order scheme \eqref{eq:first-order-time}. Other semi-implicit time discretizations can be found e.g. in \cite{AvgerinosEtAl2019,PareschiRusso2005,BoscheriPareschi2021,BoscheriTavelli2022,BirkeBoscheriKlingenberg2024}.
	
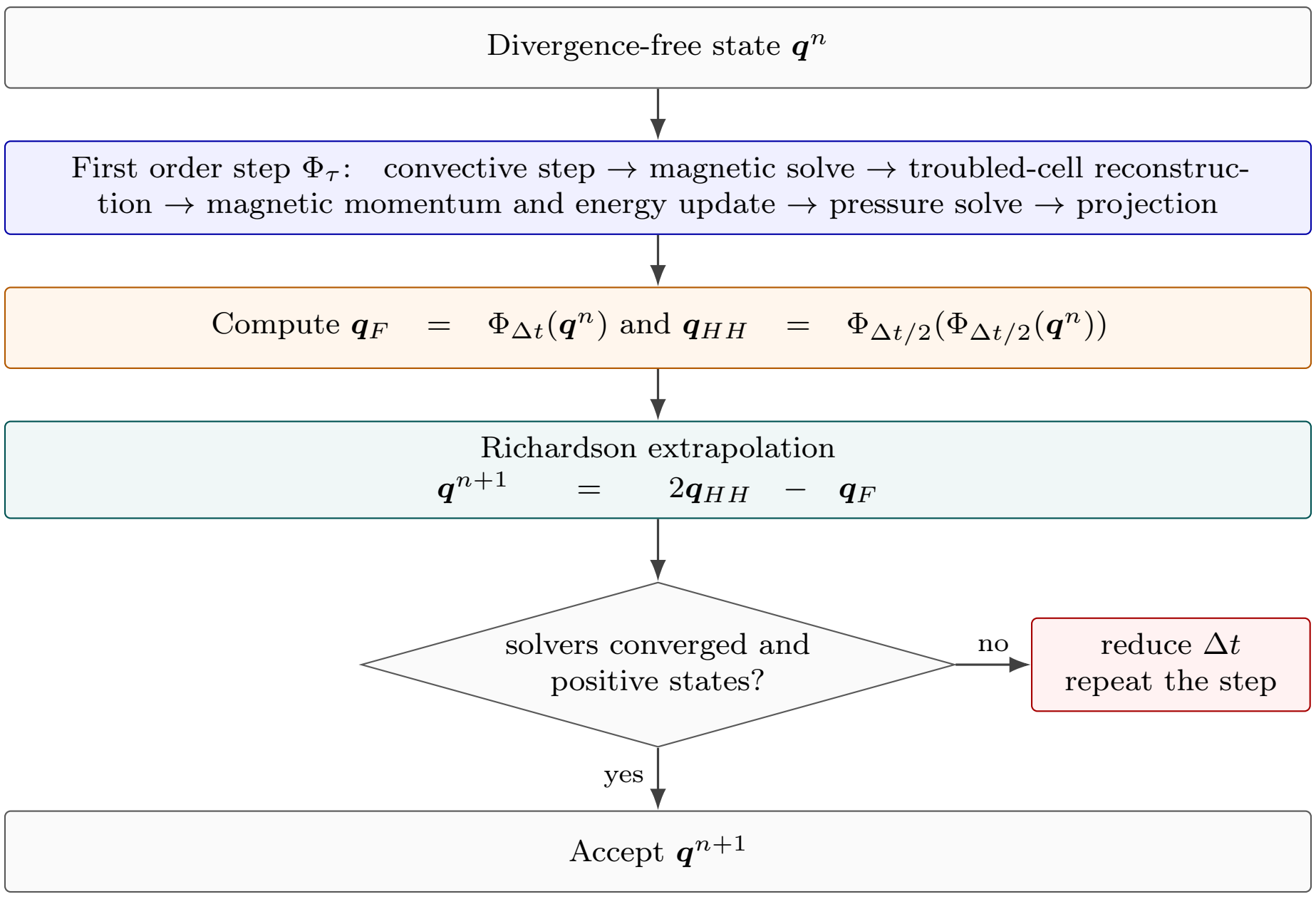
\begin{figure}[!htb]
		\centering
		\begin{tikzpicture}[
			font=\footnotesize,
			node distance=4.8mm,
			flow/.style={-Latex,semithick,draw=black!75},
			state/.style={draw=black!65,rounded corners=1.5pt,align=center,
				text width=.72\textwidth,minimum height=7.5mm,inner sep=4pt,
				fill=black!2},
			explicit/.style={state,draw=orange!70!black,fill=orange!7},
			implicit/.style={state,draw=blue!65!black,fill=blue!6},
			compatible/.style={state,draw=teal!65!black,fill=teal!6},
			decision/.style={diamond,aspect=3.6,draw=black!65,align=center,
				inner sep=1.8pt,fill=black!2}
			]
			\node[state] (input) {Divergence-free state $\Q^n$};
			\node[implicit,below=of input] (map) {First order step $\Phi_\tau$:\quad
				convective step $\rightarrow$ magnetic solve
				$\rightarrow$ troubled-cell reconstruction $\rightarrow$ magnetic momentum and energy
				update $\rightarrow$ pressure solve $\rightarrow$ projection};
			\node[explicit,below=of map] (paths) {Compute
				$\Q_F=\Phi_{\dt}(\Q^n)$ and
				$\Q_{HH}=\Phi_{\dt/2}(\Phi_{\dt/2}(\Q^n))$};
			\node[compatible,below=of paths] (combine) {Richardson extrapolation\\[1pt]
				$\Q^{n+1}=2\Q_{HH}-\Q_F$};
			\node[decision,below=5.8mm of combine] (admissible) {solvers converged and\\positive states?};
			\node[state,below=5.8mm of admissible] (output) {Accept $\Q^{n+1}$};
			\node[state,draw=red!65!black,fill=red!5,right=7mm of admissible,
			text width=.14\textwidth] (retry) {reduce $\dt$\\repeat the step};
			
			\draw[flow] (input) -- (map);
			\draw[flow] (map) -- (paths);
			\draw[flow] (paths) -- (combine);
			\draw[flow] (combine) -- (admissible);
			\draw[flow] (admissible) -- node[left,font=\scriptsize]{yes} (output);
			\draw[flow] (admissible) -- node[above,font=\scriptsize]{no} (retry);
		\end{tikzpicture}
		\caption{Flowchart of the second order semi-implicit scheme.}
		\label{fig:semi-implicit-flowchart}
	\end{figure}

	\paragraph{Time step rejection}
	If one of the linear solvers does not converge, or if a cell state is not admissible, e.g. with negative density or pressure, the whole time step is rejected and repeated with a smaller time step $\dt$, according to the criteria given in \ref{app:acceptance-criteria}. The complete algorithm is summarized in \Cref{fig:semi-implicit-flowchart}.
	
	\section{Numerical results}
	\label{sec:numerics}
	
	In this section we assess the new scheme on the test cases proposed in \cite{BoscheriThomann2024}: a smooth traveling vortex at different Mach and Alfv\'en numbers, three planar Riemann problems, a cylindrical blast wave, the MHD rotor, the Orszag--Tang vortex and the low Mach advection of a magnetic field loop. All tests use the same reconstruction, numerical dissipation and acceptance criteria. The meshes and the time integrators are listed in \Cref{tab:numerical-protocol}, where SI2 denotes the second order Richardson extrapolation \eqref{eq:richardson-composition}.
	
	\begin{table}[!htbp]
		\centering
		\footnotesize
		\setlength{\tabcolsep}{3pt}
		\caption{Meshes and time integrators.}
		\label{tab:numerical-protocol}
		\begin{tabular*}{\linewidth}{@{\extracolsep{\fill}}lcc@{}}
			\toprule
			Test & mesh & time integrator\\
			\midrule
			Smooth vortex & $64^2$--$256^2$ & SI2\\
			RP1, RP2-2D, RP3-2D & $2000\times4$ & first order\\
			Cylindrical blast & $1024^2$ & first order\\
			MHD rotor & $1024^2$ & SI2\\
			Orszag--Tang vortex & $1024^2$ & SI2\\
			Field-loop advection & $512\times256$ & SI2\\
			\bottomrule
		\end{tabular*}
	\end{table}
	
	The time step depends on the switch $\sigma$ introduced in \eqref{eq:explicit-induction-flux}. With the material speed $s^{\mathrm{mat}}_{f,q}=\max_\pm|\uvec^\pm_{f,q}\cdot\nvec_{cf}|$ at the face quadrature points, and the directional fast speed $\lambda_{d,c}=|u_{d,c}|+c_f(\Q_c;e_d)$ given by \eqref{eq:fast-magnetosonic-speed}, we set
	\begin{equation}
		\dt=\min(\dt_{\mathrm{CFL}},0.9\,\varrho_{\min}),\qquad
		\dt_{\mathrm{CFL}}=
		\begin{cases}
			\displaystyle\frac{0.9\,\varrho_{\min}}{2\max_{f,q}s^{\mathrm{mat}}_{f,q}}, & \sigma=0,\\[3mm]
			\displaystyle\frac{0.9}{\max_c\left(\lambda_{x,c}/\Delta x+\lambda_{y,c}/\Delta y\right)}, & \sigma=1,
		\end{cases}
		\label{eq:numerical-time-step}
	\end{equation}
	where the factor two is the number of space dimensions and the bound $0.9\,\varrho_{\min}$ avoids an unbounded time step when the velocity vanishes. Thus, for $\sigma=0$ the time step is only restricted by the material speed, while for $\sigma=1$ the fast magnetosonic speed is used. The first time step is always computed with the fast speed.
	
	\paragraph{Diagnostics}
	Let $M(t)=\sum_{\omega_c\in\Th}|\omega_c|\rho_c(t)$ and $\mathcal E(t)=\sum_{\omega_c\in\Th}|\omega_c|E_c(t)$ be the total mass and energy. Besides the minimum density and pressure, we monitor the relative conservation errors
	\begin{equation}
		\varepsilon_M(t)=\frac{|M(t)-M(0)|}{\max(1,|M(0)|)},
		\qquad
		\varepsilon_E(t)=\frac{|\mathcal E(t)-\mathcal E(0)|}
		{\max(1,|\mathcal E(0)|)},
		\label{eq:conservation-errors}
	\end{equation}
	and the two parts of the magnetic involution, namely the local divergence and the normal jumps across the faces:
	\begin{subequations}\label{eq:div-diagnostics}
		\begin{align}
			\frakD_\infty&=\max_c\|\frakD_c(\B_h)\|_{L^\infty(\omega_c)},
			\label{eq:local-div-diagnostic}\\
			\frakJ_{\infty,q}&=\max_{f\in\Ehi}\max_q
			|\jump{\B_h\cdot\nvec_f}(\xvec_{fq})|.
			\label{eq:jump-infinity}
		\end{align}
	\end{subequations}
	The jumps are directly evaluated from the magnetic polynomials at the face quadrature points, thus independently of the residual of the projection \eqref{eq:schur-projection}. The extrema over the whole simulation are collected in \Cref{tab:quality-indicators}. All computations keep positive density and pressure, and the local divergence is exactly zero. Among these six computations, only the blast wave requires rejected time steps (57 in total). We remark that, with physical boundaries, the total mass and energy also change because of the boundary fluxes.
	
	\begin{table}[!htbp]
		\centering
		\footnotesize
		\setlength{\tabcolsep}{3pt}
		\caption{Minimum density and pressure, and maximum of the magnetic involution diagnostics \eqref{eq:div-diagnostics}, over the whole simulation including the initial data.}
		\label{tab:quality-indicators}
		\sisetup{table-format=1.3e-1}
		\begin{tabular*}{\linewidth}{@{\extracolsep{\fill}}l c S S S S@{}}
			\toprule
			Test & mesh & {$\min\rho_h$} & {$\min p_h$} & {$\frakD_\infty$}
			& {$\frakJ_{\infty,q}$}\\
			\midrule
			RP1       & $2000\times4$ & 1.155e-1 & 8.773e-2 & 0.000e0 & 2.555e-9\\
			RP2-2D    & $2000\times4$ & 9.891e-1 & 9.500e-1 & 0.000e0 & 3.108e-9\\
			RP3-2D    & $2000\times4$ & 2.000e-1 & 1.976e-1 & 0.000e0 & 3.162e-9\\
			Blast wave& $1024^2$     & 1.165e-1 & 1.031e-2 & 0.000e0 & 3.823e-12\\
			Rotor     & $1024^2$     & 6.318e-1 & 6.212e-2 & 0.000e0 & 1.377e-13\\
			Orszag--Tang & $1024^2$  & 4.150e-2 & 8.164e-3 & 0.000e0 & 4.188e-14\\
			\bottomrule
		\end{tabular*}
	\end{table}
	
	\subsection{Convergence study at different Mach and Alfv\'en numbers}
	\label{sec:num-convergence}
	
	We consider the smooth vortex introduced in \cite{BoscheriThomann2024}, which is transported on the periodic domain $\Omega=[-5,5]^2$ up to $t_f=1$. With $\widetilde x=x-t$, $\widetilde y=y-t$, $r^2=\widetilde x^2+\widetilde y^2$, and
	\begin{equation}
		f=\exp\left(\frac{1-r^2}{2}\right),\qquad
		\epsilon=\frac{\sqrt{2\pi}}{2\pi},\qquad
		\mu=\frac{\sqrt{4\pi}}{2\pi},
		\label{eq:vortex-parameters}
	\end{equation}
	the exact solution writes
	\begin{subequations}\label{eq:vortex-exact}
		\begin{align}
			\rho&=\rho_0,
			&\uvec&=(1-\epsilon\widetilde y f,
			1+\epsilon\widetilde x f)^\top,\\
			\B&=(-\mu\widetilde y f,\mu\widetilde x f)^\top,
			&p&=1+e^{1-r^2}\left[
			\frac{\mu^2(1-r^2)}{8\pi}-\frac{\rho_0\epsilon^2}{2}\right].
		\end{align}
	\end{subequations}
	By choosing the background density $\rho_0=10^{-k}$, $k=0,\ldots,5$, we modify both the Mach and the Alfv\'en numbers. The maximum Mach number $\mathcal M$ and Alfv\'en speed $c_A$ of the initial data are listed in \Cref{tab:vortex-regimes}. We report the Alfv\'en speed instead of the Alfv\'en Mach number, since the latter is singular where the magnetic field vanishes. We set $\gamma=5/3$ and use uniform meshes with $N\in\{64,128,256\}$ cells per direction.
	
	\begin{table}[!htbp]
		\centering
		\footnotesize
		\setlength{\tabcolsep}{3pt}
		\caption{Maximum Mach number and Alfv\'en speed of the initial data of the smooth vortex.}
		\label{tab:vortex-regimes}
		\begin{tabular*}{\linewidth}{@{\extracolsep{\fill}}lrrrrrr@{}}
			\toprule
			$\rho_0$ & $10^0$ & $10^{-1}$ & $10^{-2}$ & $10^{-3}$ & $10^{-4}$ & $10^{-5}$\\
			\midrule
			$\max \mathcal M$ & $1.473\mathrm{e}{+0}$ & $4.459\mathrm{e}{-1}$ & $1.405\mathrm{e}{-1}$
			& $4.442\mathrm{e}{-2}$ & $1.405\mathrm{e}{-2}$ & $4.442\mathrm{e}{-3}$\\
			$\max c_A$ & $1.592\mathrm{e}{-1}$ & $5.033\mathrm{e}{-1}$ & $1.592\mathrm{e}{+0}$
			& $5.033\mathrm{e}{+0}$ & $1.592\mathrm{e}{+1}$ & $5.033\mathrm{e}{+1}$\\
			\bottomrule
		\end{tabular*}
	\end{table}
	
	The $L^2$ error of a variable $w$ and the corresponding order of convergence are computed as
	\begin{equation}
		e_2(w;h)=\left[\frac1{|\Omega|}\sum_{\omega_c\in\Th}|\omega_c|
		\bigl(w_c(t_f)-w(\xvec_c,t_f)\bigr)^2\right]^{1/2},
		\qquad
		\mathcal O_2=\log_2\frac{e_2(w;2h)}{e_2(w;h)}.
		\label{eq:observed-order}
	\end{equation}
	Let us notice that the exact solution is evaluated at the cell centers, which differs from the exact cell average by $O(h^2)$, and that the Gaussian tails of the vortex are not exactly periodic on the truncated domain. No troubled cells are detected in this test, thus $\sigma=0$ throughout the simulation and the time step \eqref{eq:numerical-time-step} only depends on the material speed, independently of the acoustic and Alfv\'en speeds.
	
	\begin{table}[!htbp]
		\centering
		\footnotesize
		\setlength{\tabcolsep}{3pt}
		\caption{Smooth vortex: $L^2$ errors and orders of convergence at $t_f=1$. The pressure is computed with \eqref{eq:cell-polynomial-pressure}, and the magnetic errors refer to the cell averages.}
		\label{tab:convergence-summary}
		\begin{tabular*}{\linewidth}{@{\extracolsep{\fill}}cc*{5}{S[table-format=1.3e-1]S[table-format=-1.3]}@{}}
			\toprule
			& & \multicolumn{2}{c}{$\rho$} & \multicolumn{2}{c}{$u$}
			& \multicolumn{2}{c}{$p$} & \multicolumn{2}{c}{$B_x$}
			& \multicolumn{2}{c}{$B_y$}\\
			\cmidrule(lr){3-4}\cmidrule(lr){5-6}\cmidrule(lr){7-8}
			\cmidrule(lr){9-10}\cmidrule(lr){11-12}
			{$\rho_0$} & {$N$} & {error} & {$\mathcal O_2$} & {error} & {$\mathcal O_2$}
			& {error} & {$\mathcal O_2$} & {error} & {$\mathcal O_2$} & {error} & {$\mathcal O_2$}\\
			\midrule
			\multirow{3}{*}{$10^{0}$}
			& 64 & 1.624e-03 & {--} & 1.309e-03 & {--} & 1.602e-03 & {--} & 2.086e-03 & {--} & 1.948e-03 & {--}\\
			& 128 & 3.946e-04 & 2.041 & 2.744e-04 & 2.254 & 4.067e-04 & 1.978 & 5.765e-04 & 1.856 & 5.431e-04 & 1.843\\
			& 256 & 9.449e-05 & 2.062 & 7.553e-05 & 1.861 & 1.026e-04 & 1.987 & 1.860e-04 & 1.632 & 1.761e-04 & 1.625\\
			\midrule
			\multirow{3}{*}{$10^{-1}$}
			& 64 & 6.033e-05 & {--} & 1.883e-03 & {--} & 1.865e-04 & {--} & 2.980e-03 & {--} & 2.907e-03 & {--}\\
			& 128 & 1.605e-05 & 1.910 & 3.252e-04 & 2.534 & 4.609e-05 & 2.017 & 6.705e-04 & 2.152 & 6.376e-04 & 2.189\\
			& 256 & 4.063e-06 & 1.982 & 8.675e-05 & 1.906 & 1.189e-05 & 1.954 & 1.913e-04 & 1.809 & 1.818e-04 & 1.811\\
			\midrule
			\multirow{3}{*}{$10^{-2}$}
			& 64 & 5.720e-06 & {--} & 4.474e-03 & {--} & 2.831e-04 & {--} & 6.007e-03 & {--} & 6.010e-03 & {--}\\
			& 128 & 1.583e-06 & 1.853 & 7.685e-04 & 2.541 & 5.467e-05 & 2.373 & 1.119e-03 & 2.424 & 1.111e-03 & 2.435\\
			& 256 & 3.964e-07 & 1.998 & 1.825e-04 & 2.074 & 1.221e-05 & 2.163 & 2.529e-04 & 2.146 & 2.489e-04 & 2.158\\
			\midrule
			\multirow{3}{*}{$10^{-3}$}
			& 64 & 5.454e-07 & {--} & 1.412e-02 & {--} & 7.130e-04 & {--} & 1.475e-02 & {--} & 1.474e-02 & {--}\\
			& 128 & 1.849e-07 & 1.561 & 2.319e-03 & 2.606 & 1.469e-04 & 2.279 & 2.796e-03 & 2.399 & 2.796e-03 & 2.398\\
			& 256 & 4.272e-08 & 2.113 & 5.894e-04 & 1.976 & 3.001e-05 & 2.292 & 5.414e-04 & 2.369 & 5.414e-04 & 2.369\\
			\midrule
			\multirow{3}{*}{$10^{-4}$}
			& 64 & 5.304e-08 & {--} & 7.786e-02 & {--} & 1.496e-03 & {--} & 3.287e-02 & {--} & 3.288e-02 & {--}\\
			& 128 & 2.718e-08 & 0.965 & 1.272e-02 & 2.614 & 4.077e-04 & 1.875 & 7.773e-03 & 2.080 & 7.773e-03 & 2.081\\
			& 256 & 6.061e-09 & 2.165 & 2.411e-03 & 2.399 & 8.393e-05 & 2.280 & 1.485e-03 & 2.388 & 1.485e-03 & 2.388\\
			\midrule
			\multirow{3}{*}{$10^{-5}$}
			& 64 & 1.128e-08 & {--} & 3.481e-01 & {--} & 2.334e-03 & {--} & 5.851e-02 & {--} & 5.851e-02 & {--}\\
			& 128 & 4.121e-09 & 1.452 & 7.285e-02 & 2.257 & 9.835e-04 & 1.247 & 1.995e-02 & 1.553 & 1.994e-02 & 1.553\\
			& 256 & 1.207e-09 & 1.772 & 1.123e-02 & 2.698 & 2.383e-04 & 2.045 & 4.280e-03 & 2.220 & 4.280e-03 & 2.220\\
			\bottomrule
		\end{tabular*}
	\end{table}
	
	The results are reported in \Cref{tab:convergence-summary}. For all background densities, i.e. for Mach numbers ranging from $1.5$ down to $4.4\times10^{-3}$, the errors decrease under mesh refinement and second order of accuracy is retrieved. On the two finest meshes, the orders of convergence range from $1.954$ to $2.292$ for the pressure, from $1.861$ to $2.698$ for the velocity and from $1.625$ to $2.388$ for the magnetic field. The density rates are more irregular on the coarsest meshes, where the density errors are very small (below $10^{-6}$ for $\rho_0\leq10^{-3}$). We recall that the mesh size and the time step are refined simultaneously.
	
	\subsection{Planar Riemann problems}
	\label{sec:num-riemann}
	
	We consider three Riemann problems based on Table 3 of \cite{BoscheriThomann2024}, solved on a $2000\times4$ strip with $\Delta x=\Delta y$ and periodic boundaries in the $y$-direction. RP1 uses the original initial data, while RP2-2D and RP3-2D only retain the in-plane components with $w=B_z=0$, thus they differ from the full MHD problems RP2 and RP3 of \cite{BoscheriThomann2024}. The left and right states are reported in \Cref{tab:riemann-states}. The computational domain is $[-0.5,0.5]$ in the $x$-direction, the initial discontinuity is located at $x=0$ for RP1 and at $x=-0.1$ for RP2-2D and RP3-2D, and the left and right states are imposed at the corresponding boundaries. The final times are $t_f=0.10$, $0.20$ and $0.15$, respectively, and $\gamma=5/3$.
	
	\begin{table}[!htbp]
		\centering
		\footnotesize
		\setlength{\tabcolsep}{3pt}
		\caption{Left and right primitive states $(\rho,u,v,p,B_x,B_y)$ of the planar Riemann problems.}
		\label{tab:riemann-states}
		\begin{tabular*}{\linewidth}{@{\extracolsep{\fill}}lcc@{}}
			\toprule
			Case & Left state & Right state\\
			\midrule
			RP1 & $(1,0,0,1,0.75\sqrt{4\pi},\sqrt{4\pi})$
			& $(0.125,0,0,0.1,0.75\sqrt{4\pi},-\sqrt{4\pi})$\\
			RP2-2D & $(1.08,1.2,0.01,0.95,2,3.6)$
			& $(0.9891,-0.0131,0.0269,0.97159,2,4.0244)$\\
			RP3-2D & $(1.7,0,0,1.7,3.899398,3.544908)$
			& $(0.2,0,0,0.2,3.899398,2.785898)$\\
			\bottomrule
		\end{tabular*}
	\end{table}
	
\begin{figure}[!htb]
		\centering
		\figurepair{0.47}{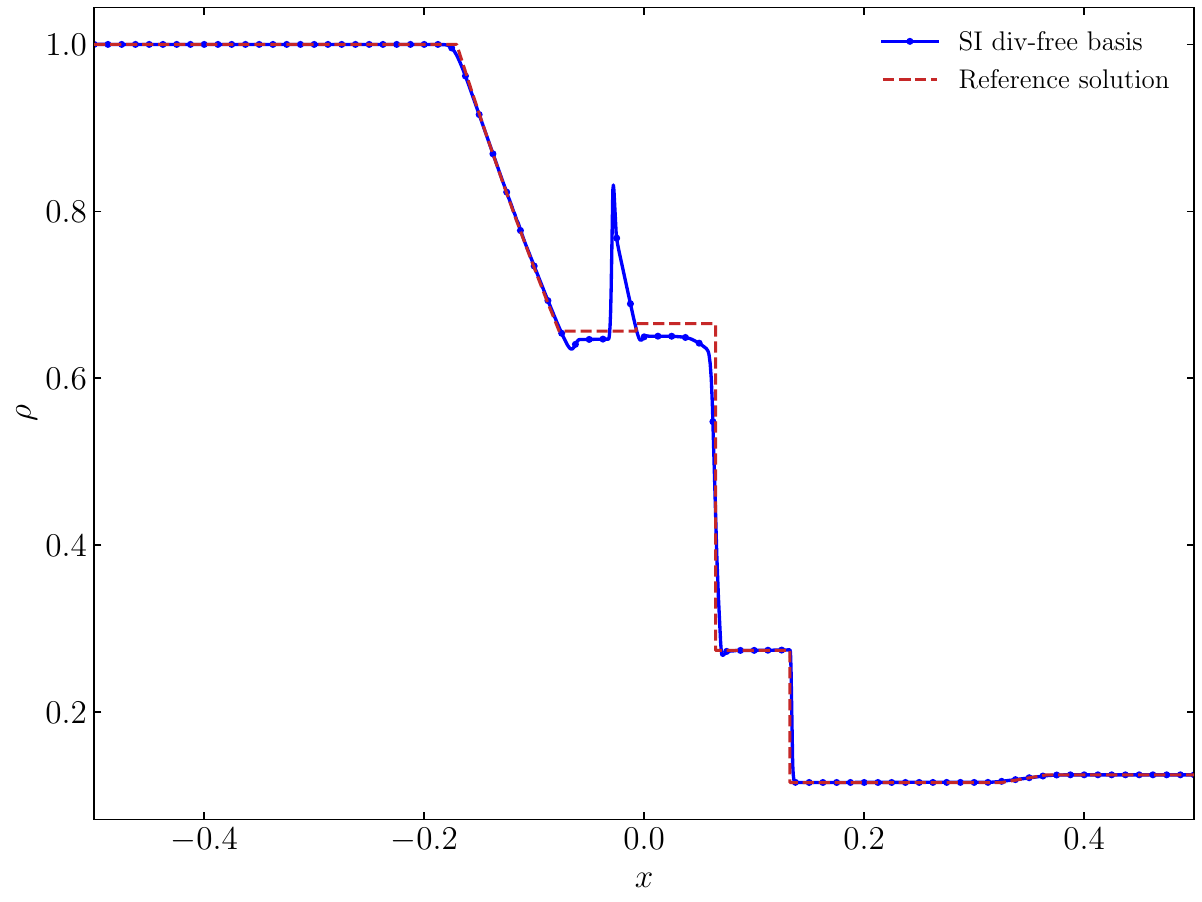}{(a) Density $\rho$}
		{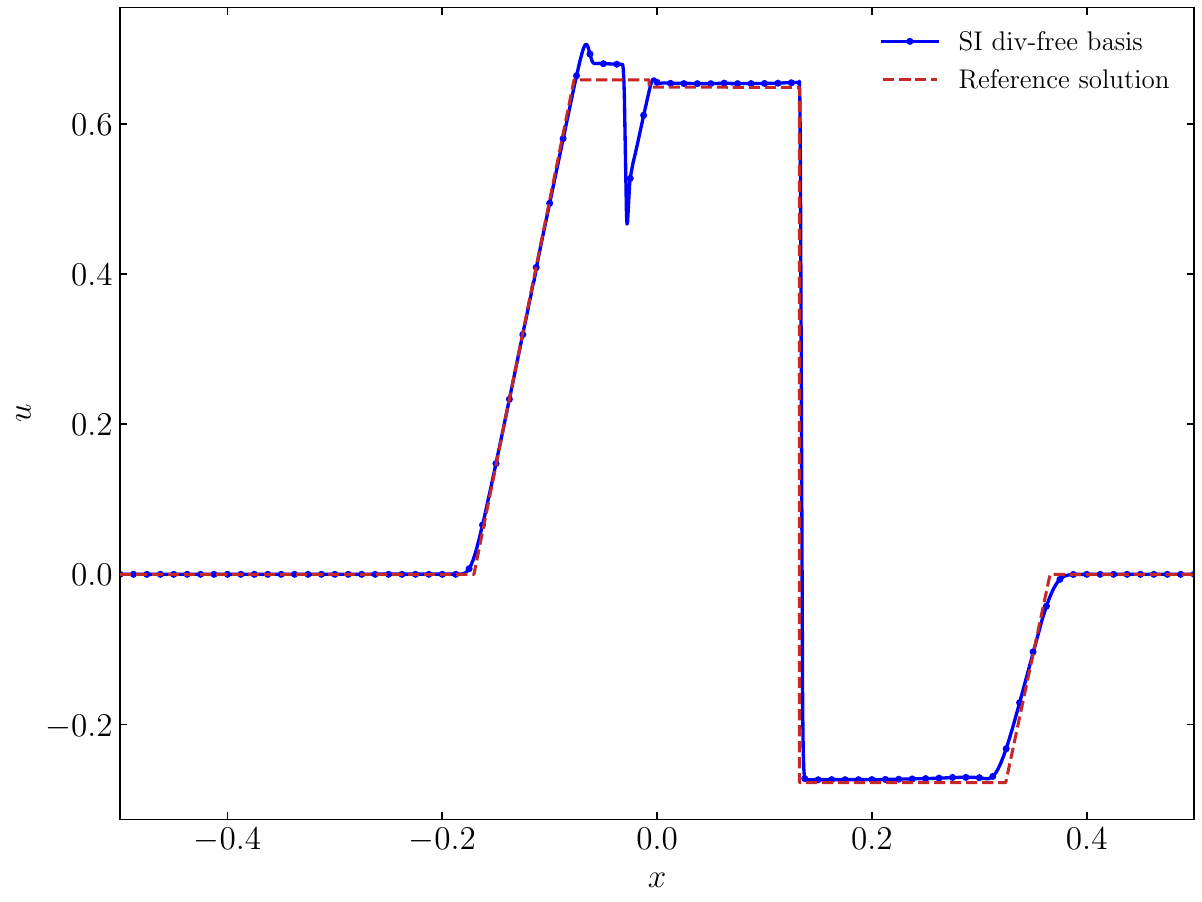}{(b) Normal velocity $u$}\par\smallskip
		\figurepair{0.47}{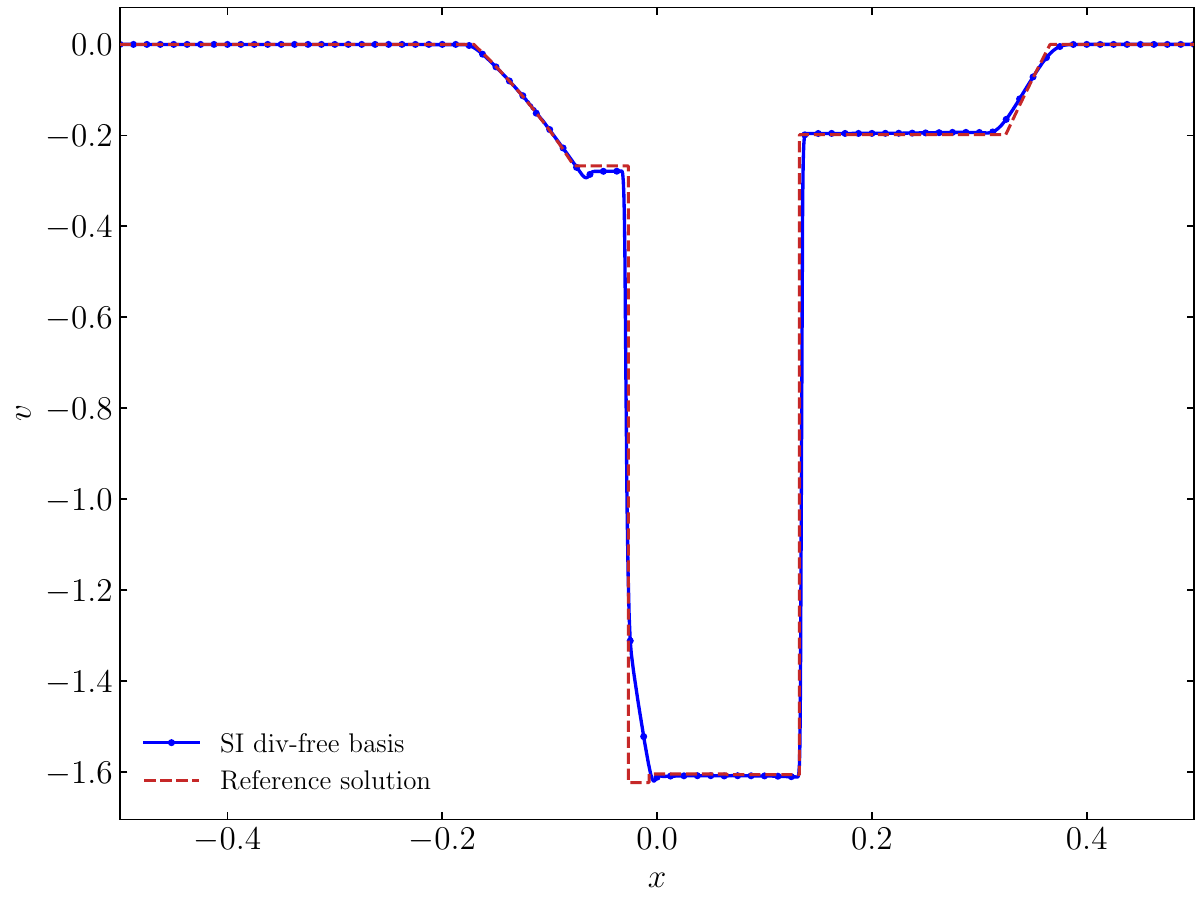}{(c) Transverse velocity $v$}
		{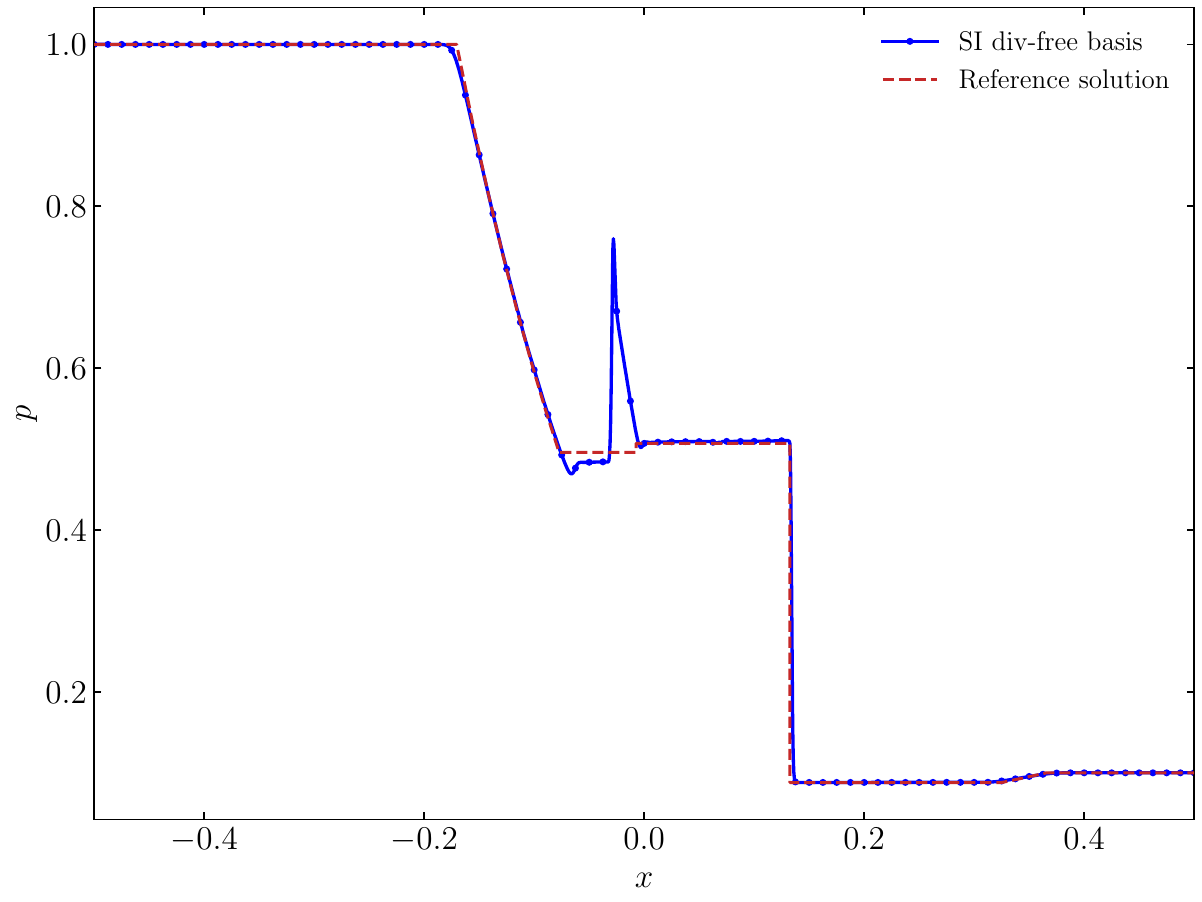}{(d) Pressure $p$}\par\smallskip
		\figurepair{0.47}{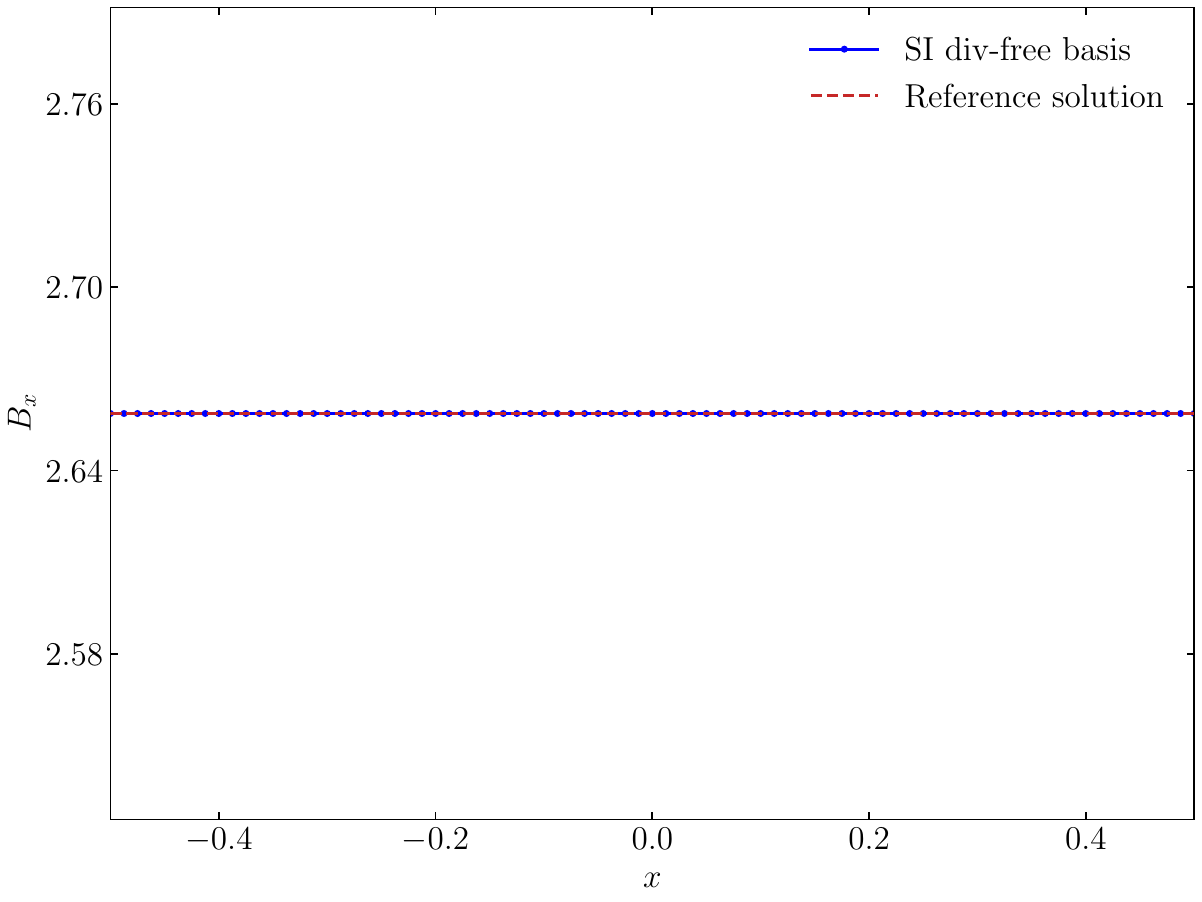}{(e) Normal field $B_x$}
		{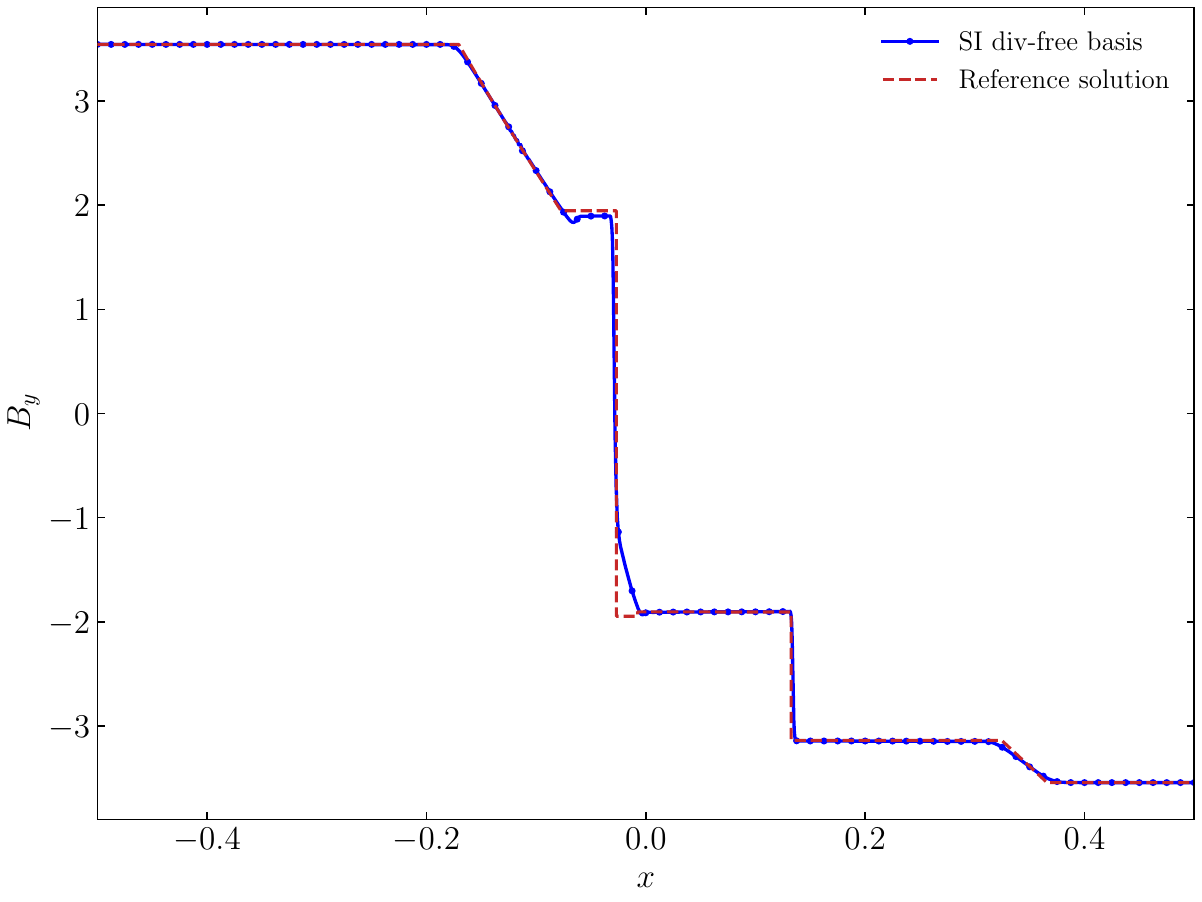}{(f) Transverse field $B_y$}
		\caption{RP1 at $t_f=0.10$ on $2000\times4$ cells: profiles along $y=0$ compared with the exact solution of \cite{Falle2002,FalleKomissarov2001}.}
		\label{fig:rp1-slices}
	\end{figure}

\begin{figure}[!htb]
		\centering
		\figurepair{0.47}{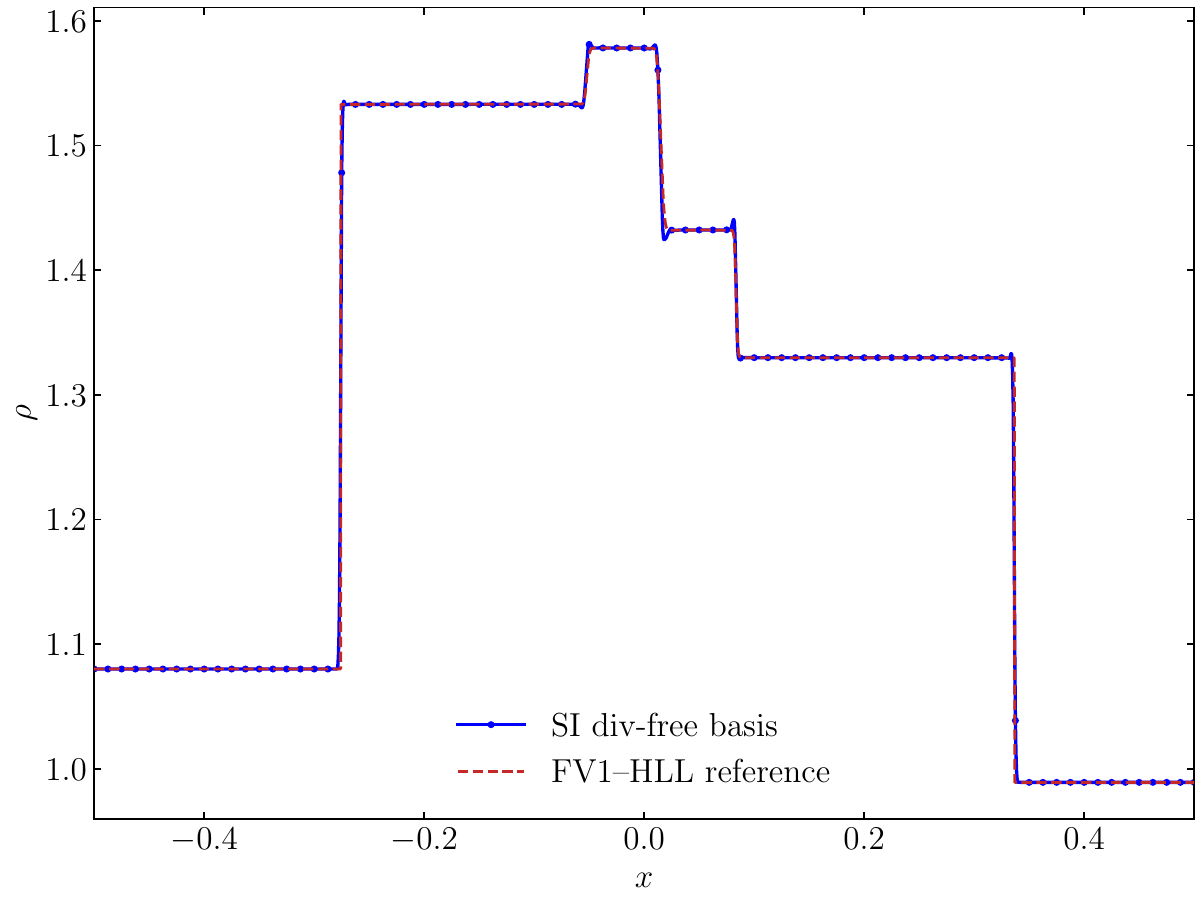}{(a) Density $\rho$}
		{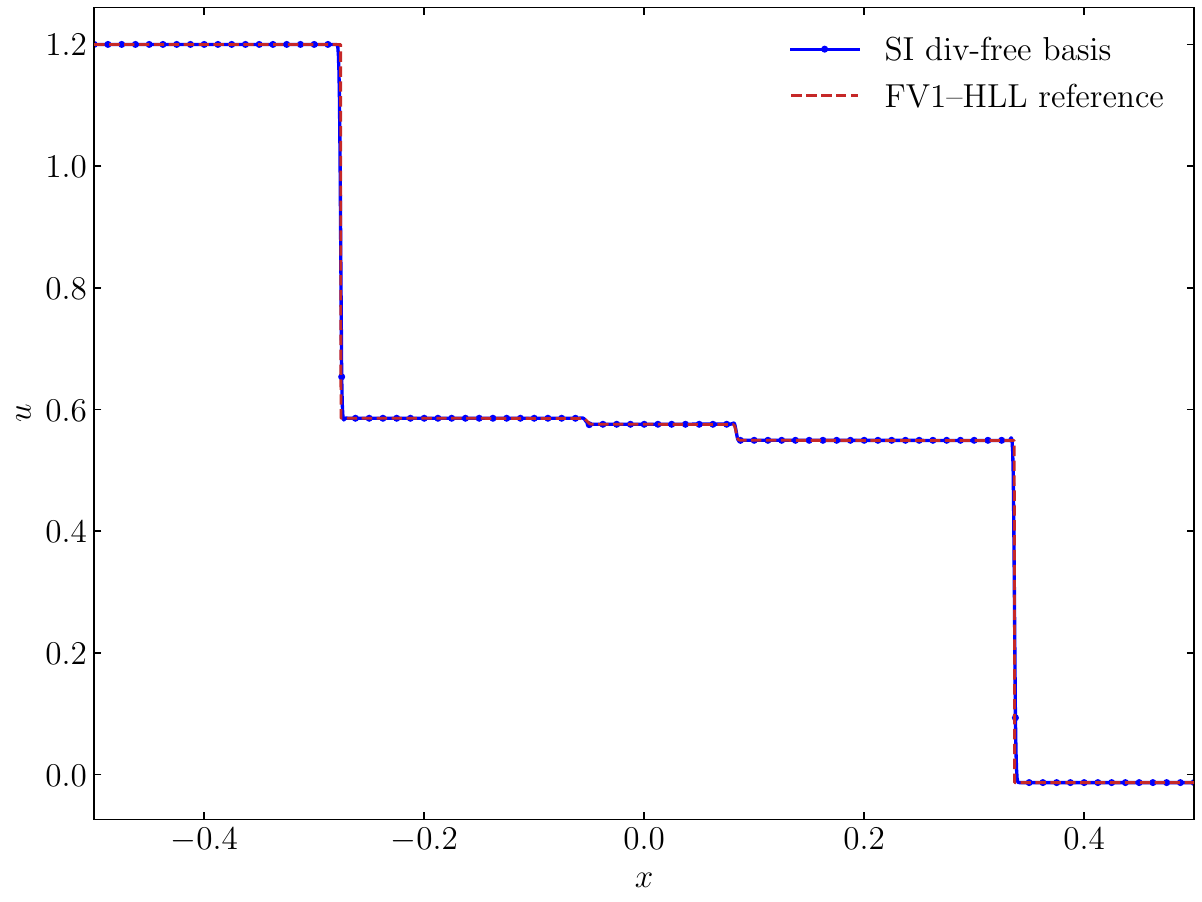}{(b) Normal velocity $u$}\par\smallskip
		\figurepair{0.47}{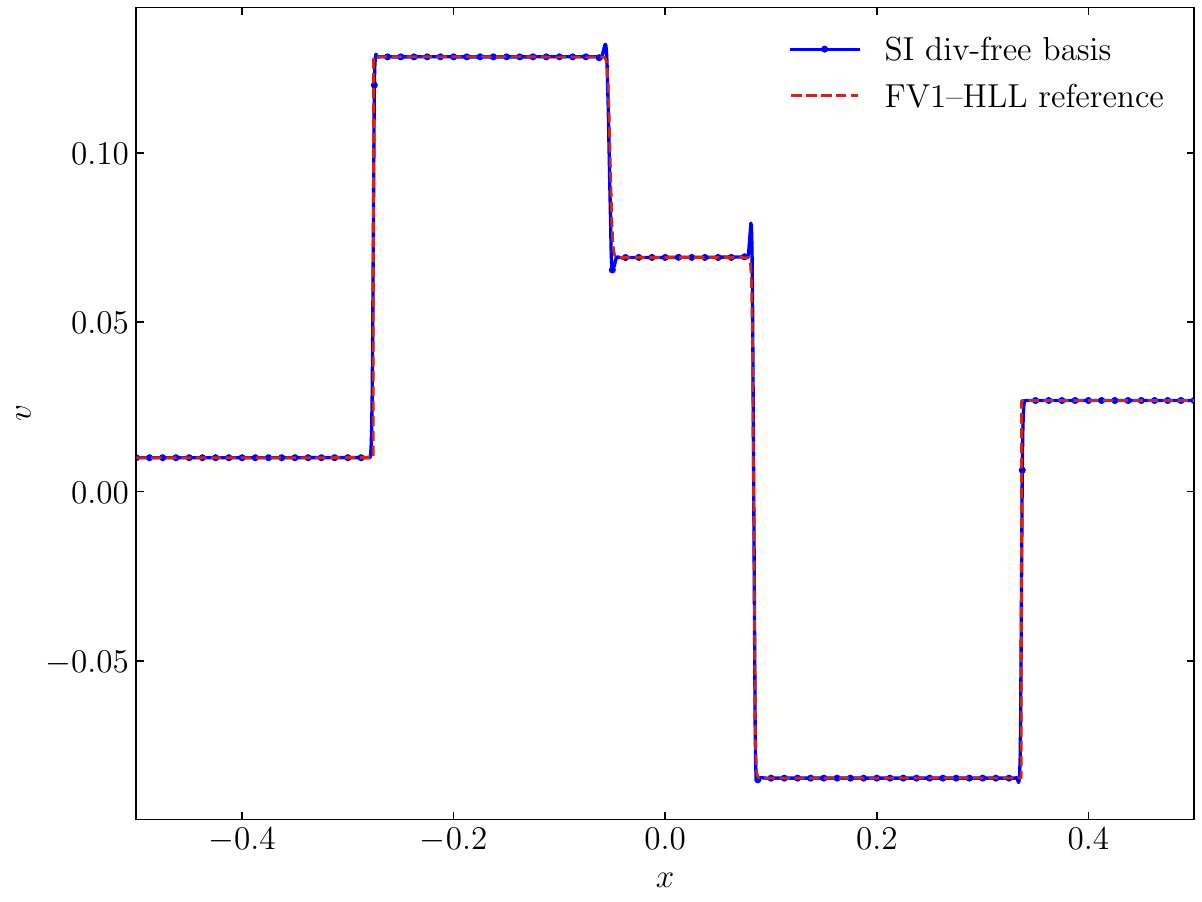}{(c) Transverse velocity $v$}
		{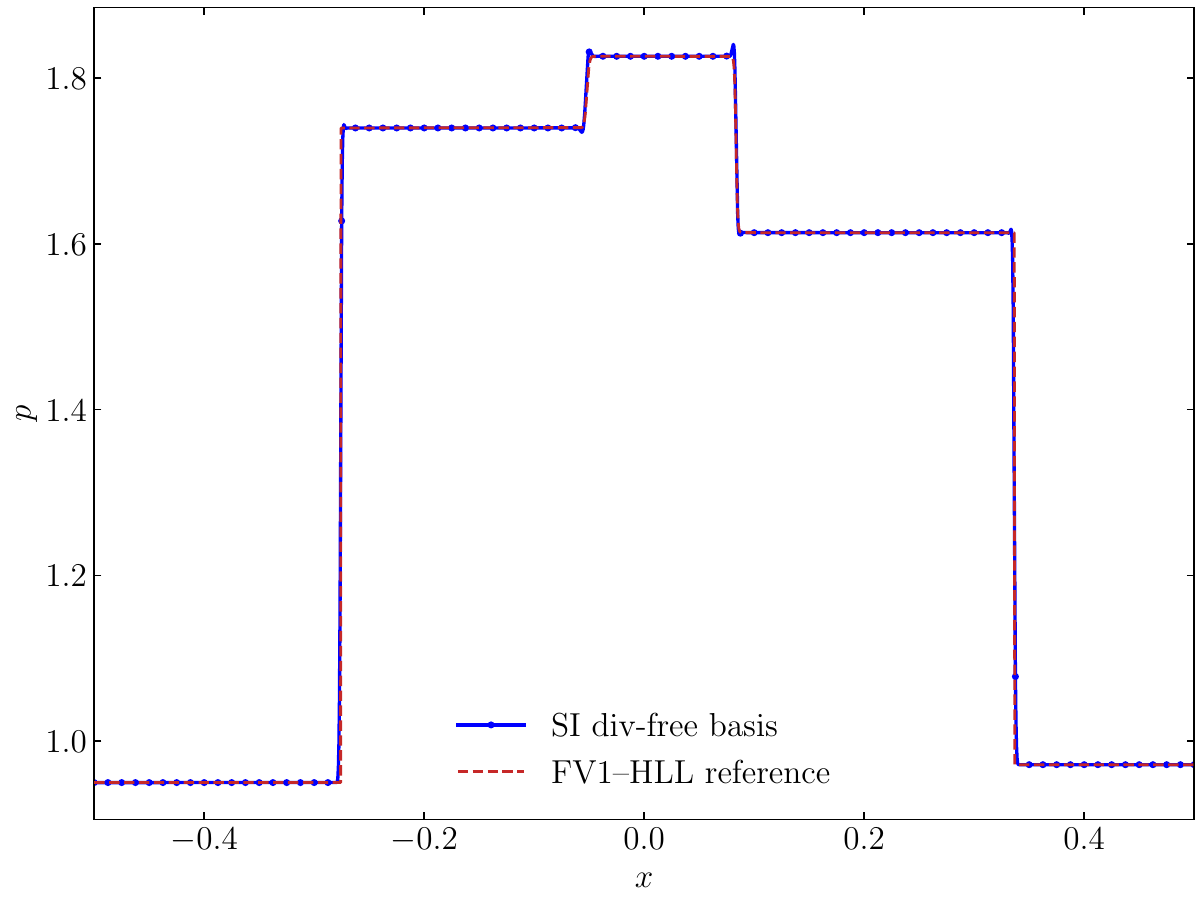}{(d) Pressure $p$}\par\smallskip
		\figurepair{0.47}{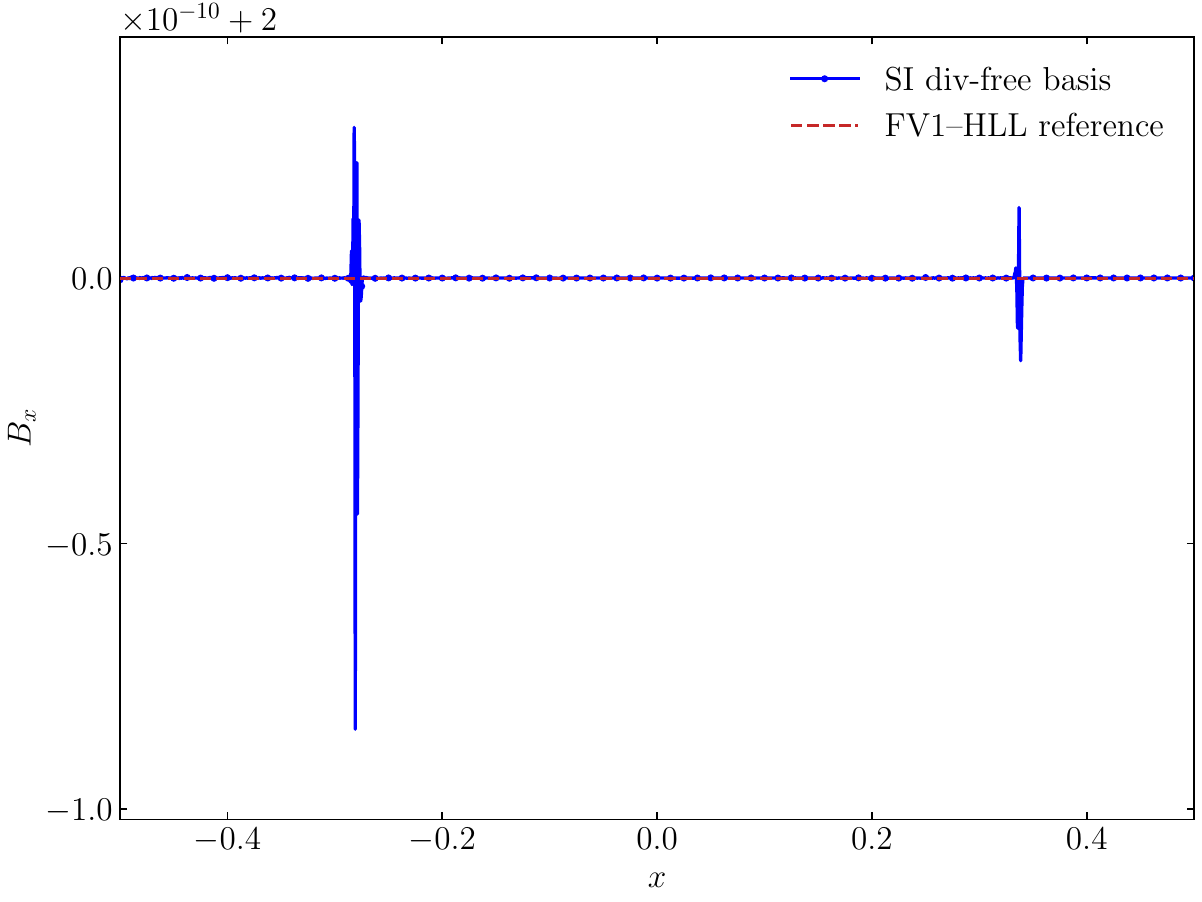}{(e) Normal field $B_x$}
		{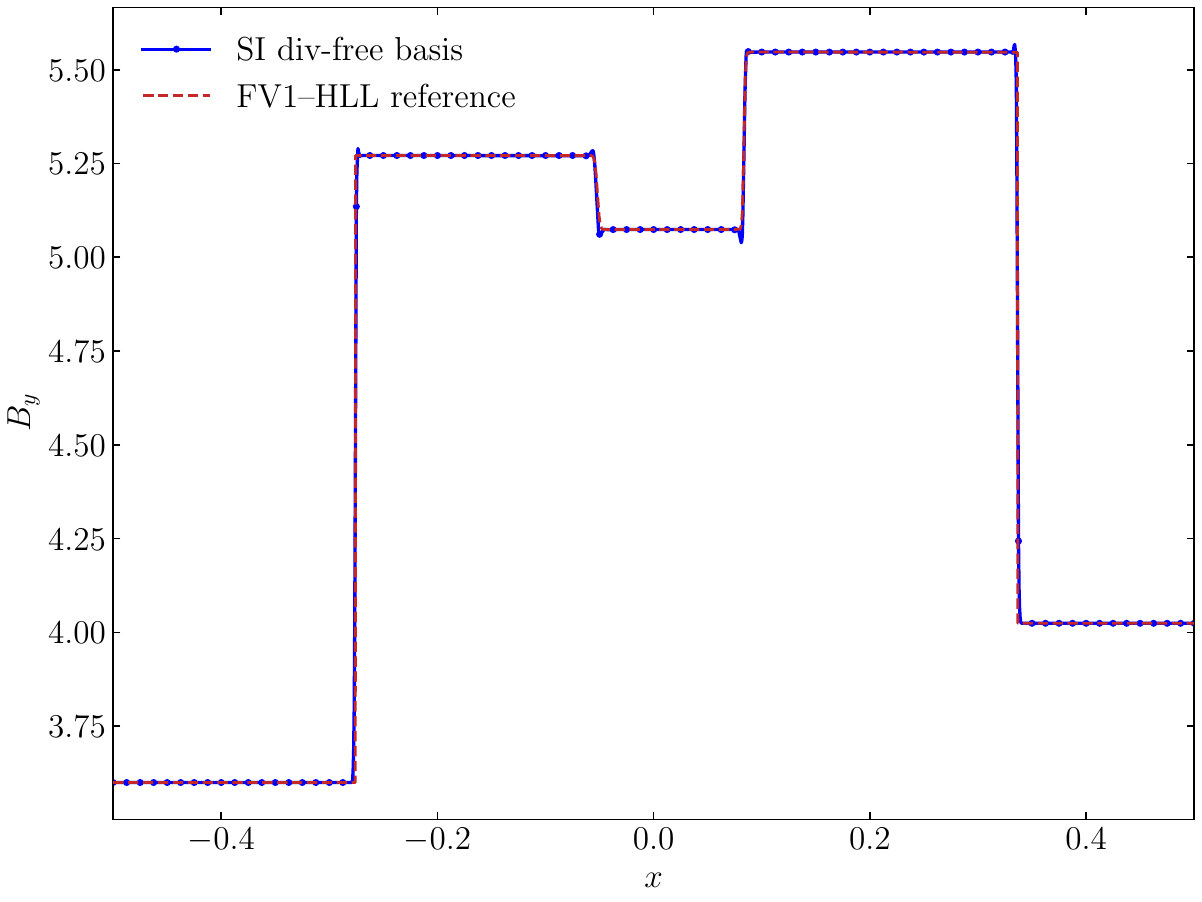}{(f) Transverse field $B_y$}
		\caption{RP2-2D at $t_f=0.20$ on $2000\times4$ cells: profiles along $y=0$ compared with a first order HLL reference solution on $64000$ cells.}
		\label{fig:rp2-slices}
	\end{figure}

\begin{figure}[!htb]
		\centering
		\figurepair{0.47}{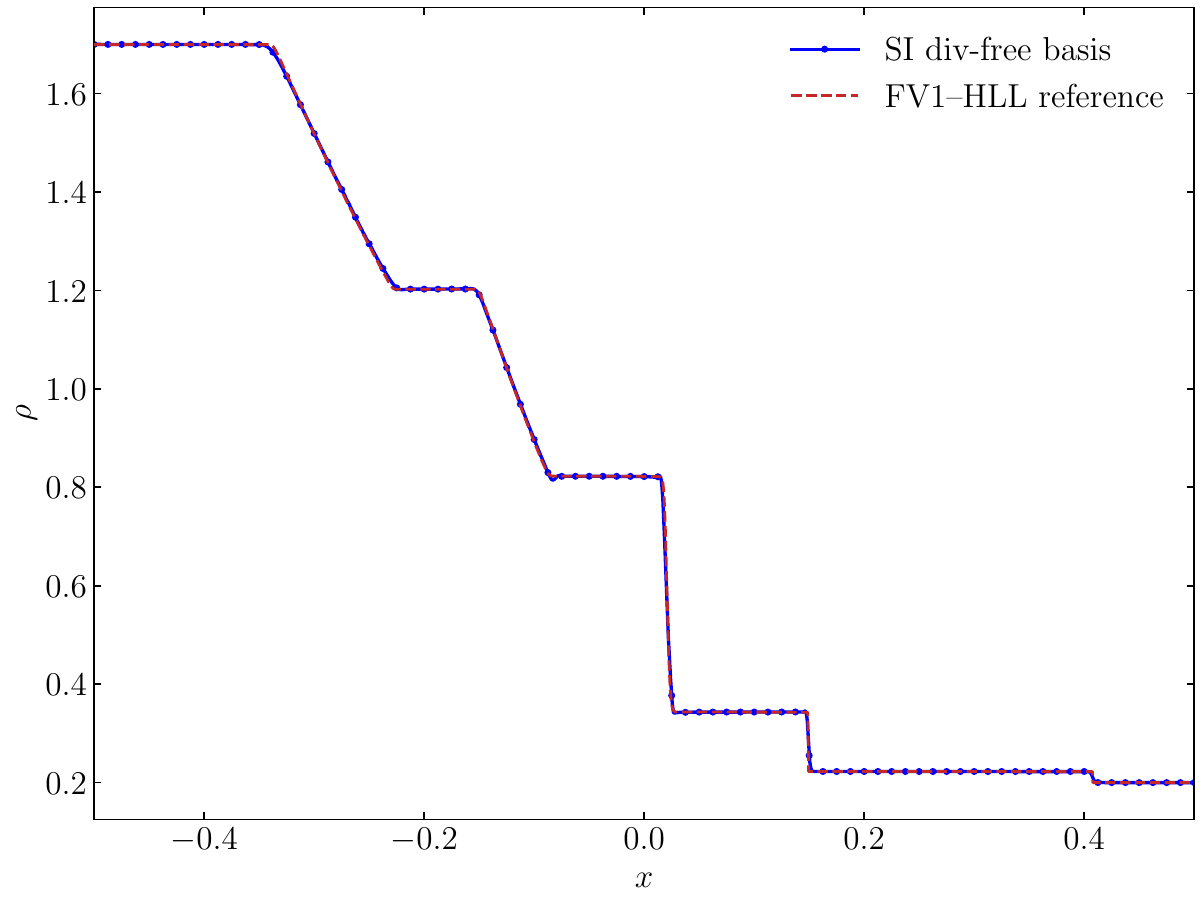}{(a) Density $\rho$}
		{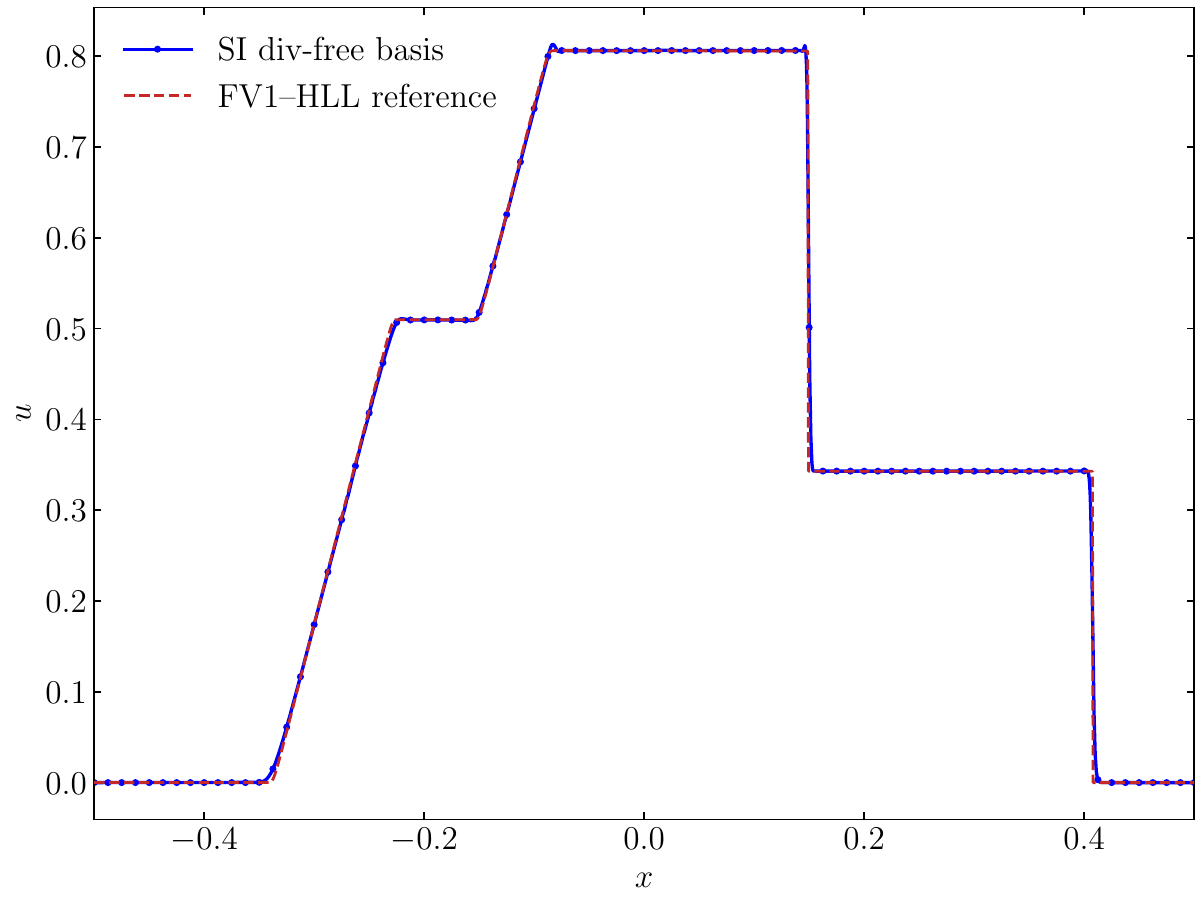}{(b) Normal velocity $u$}\par\smallskip
		\figurepair{0.47}{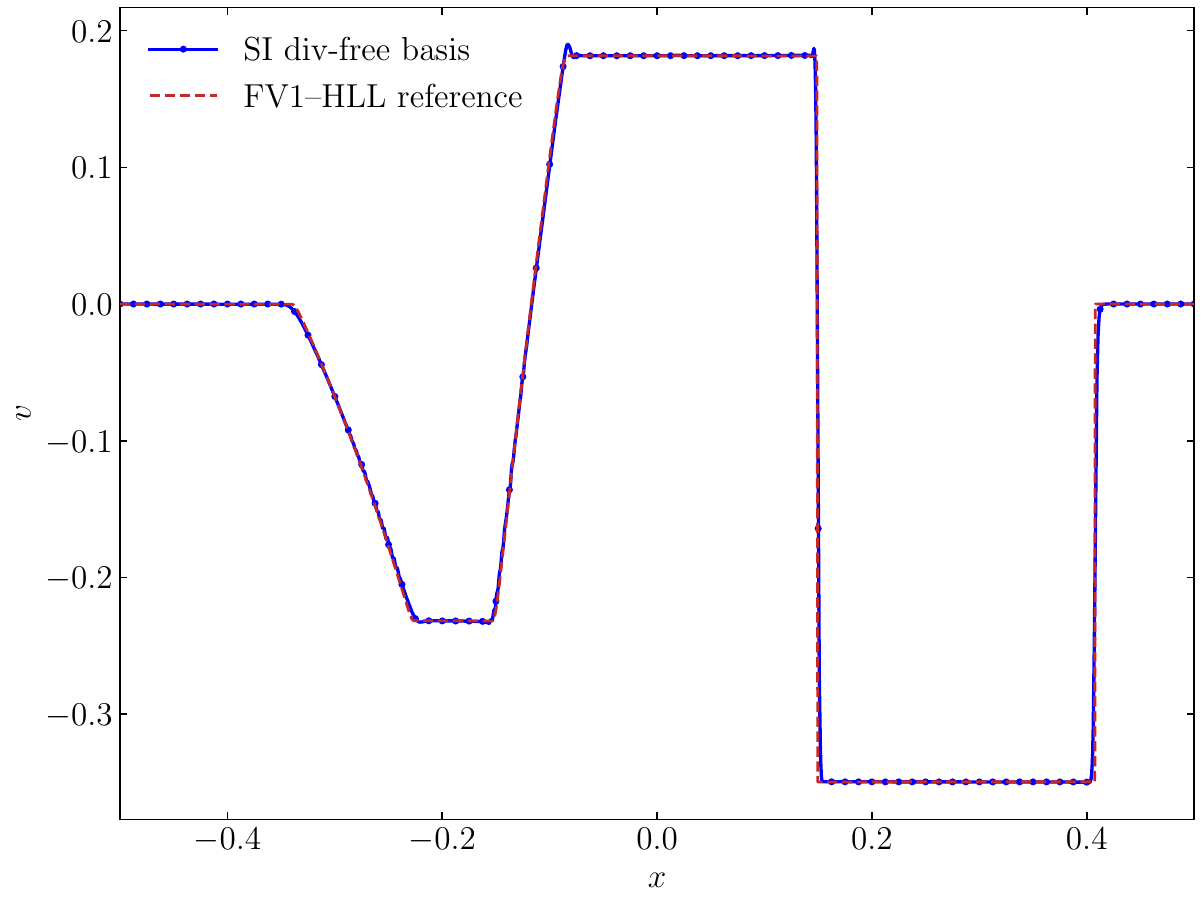}{(c) Transverse velocity $v$}
		{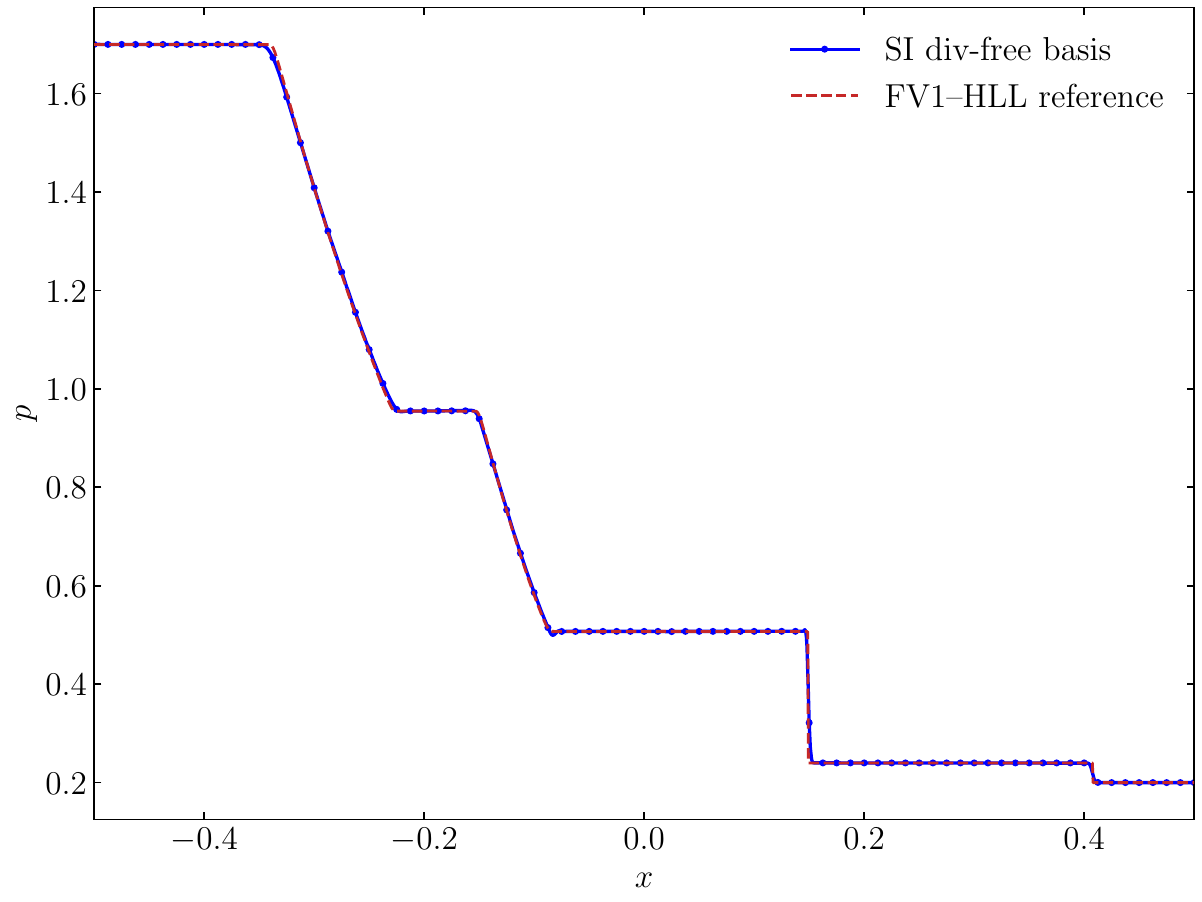}{(d) Pressure $p$}\par\smallskip
		\figurepair{0.47}{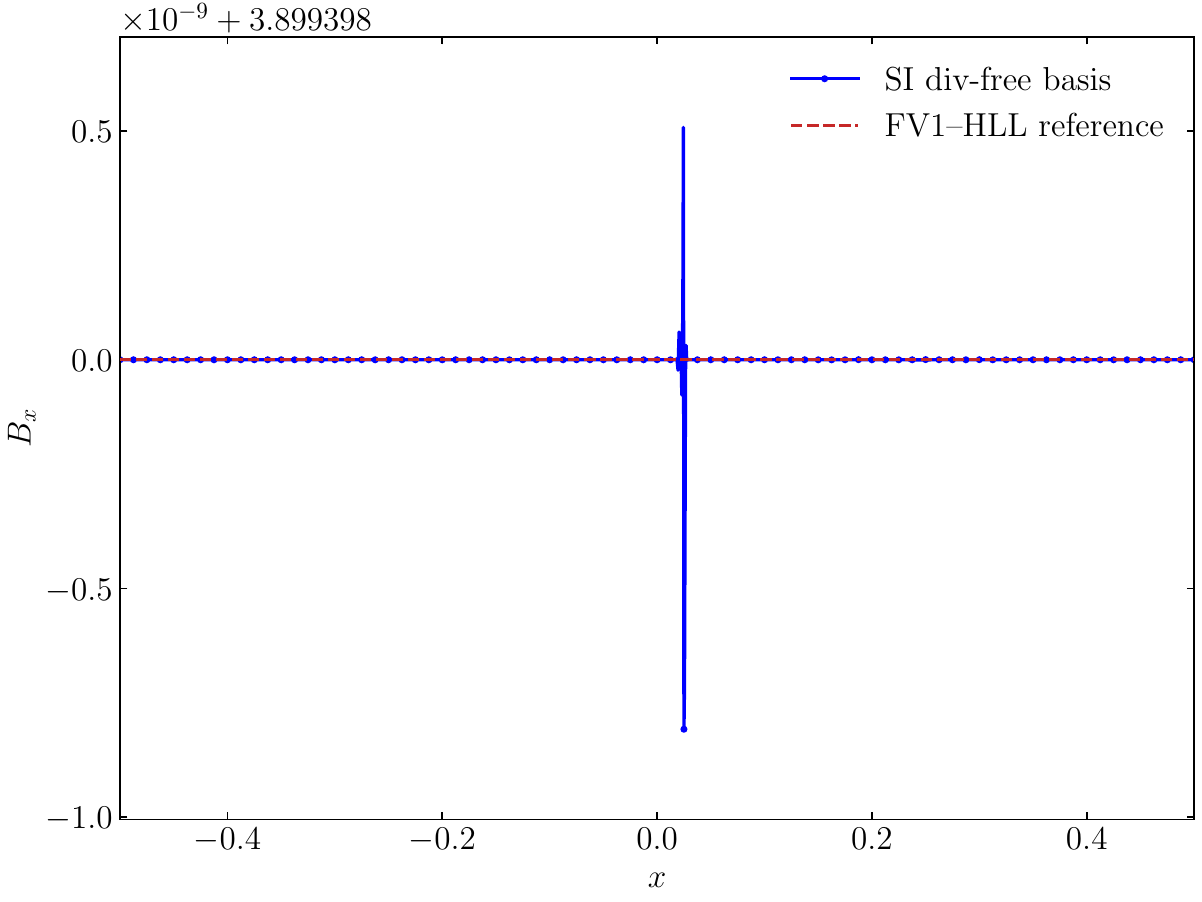}{(e) Normal field $B_x$}
		{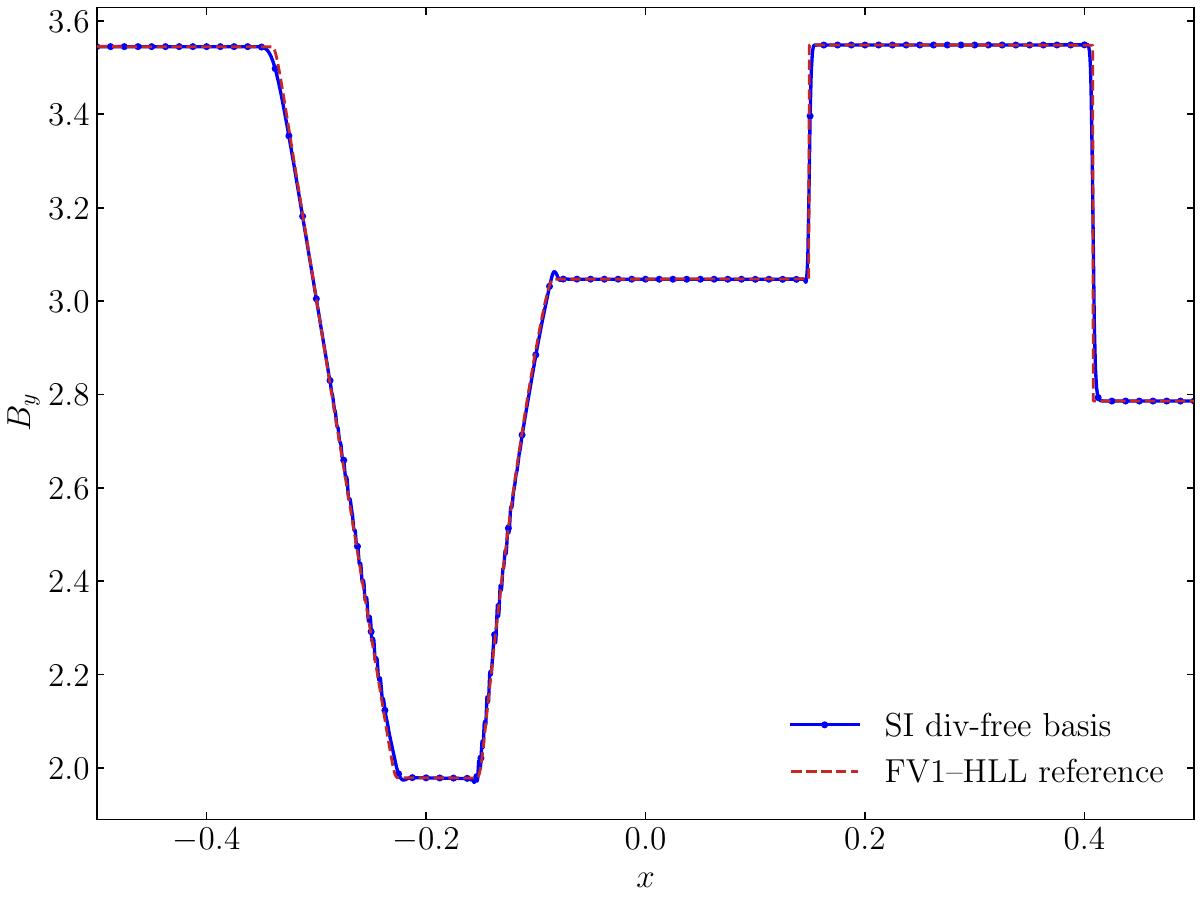}{(f) Transverse field $B_y$}
		\caption{RP3-2D at $t_f=0.15$ on $2000\times4$ cells: profiles along $y=0$ compared with a first order HLL reference solution on $64000$ cells.}
		\label{fig:rp3-slices}
	\end{figure}

	For RP1, the reference solution has been provided by Falle and Komissarov and it is obtained with their exact MHD Riemann solver \cite{Falle2002,FalleKomissarov2001}. For RP2-2D and RP3-2D, the reference solution is computed with a first order finite volume scheme on $64000$ cells, using the HLL flux \cite{HartenLaxVanLeer1983} and a forward Euler time stepping with CFL $0.8$.	The numerical results are depicted in \Cref{fig:rp1-slices,fig:rp2-slices,fig:rp3-slices}. For RP1, the sequence of waves and the intermediate states agree with the exact solution. Small overshoots in density and pressure, as well as an undershoot in velocity, are visible near the compound wave, which have also been reported for several approximate MHD solvers in \cite{BoscheriThomann2024}. For RP2-2D and RP3-2D, the position of the waves and the intermediate states agree well with the reference solution, including the sharp transitions in the transverse velocity and magnetic field, and the deviations are localized at the fronts. In all three cases, the solution remains invariant in the transverse direction, the local divergence is zero and the normal jumps stay below $3.2\times10^{-9}$, see \Cref{tab:quality-indicators}. Therefore, the involution constraints are preserved also in the presence of shocks and limiting.
	
	\subsection{Cylindrical blast wave}
	\label{sec:num-blast}
	
\begin{figure}[!htb]
		\centering
		\includegraphics[width=0.94\textwidth]{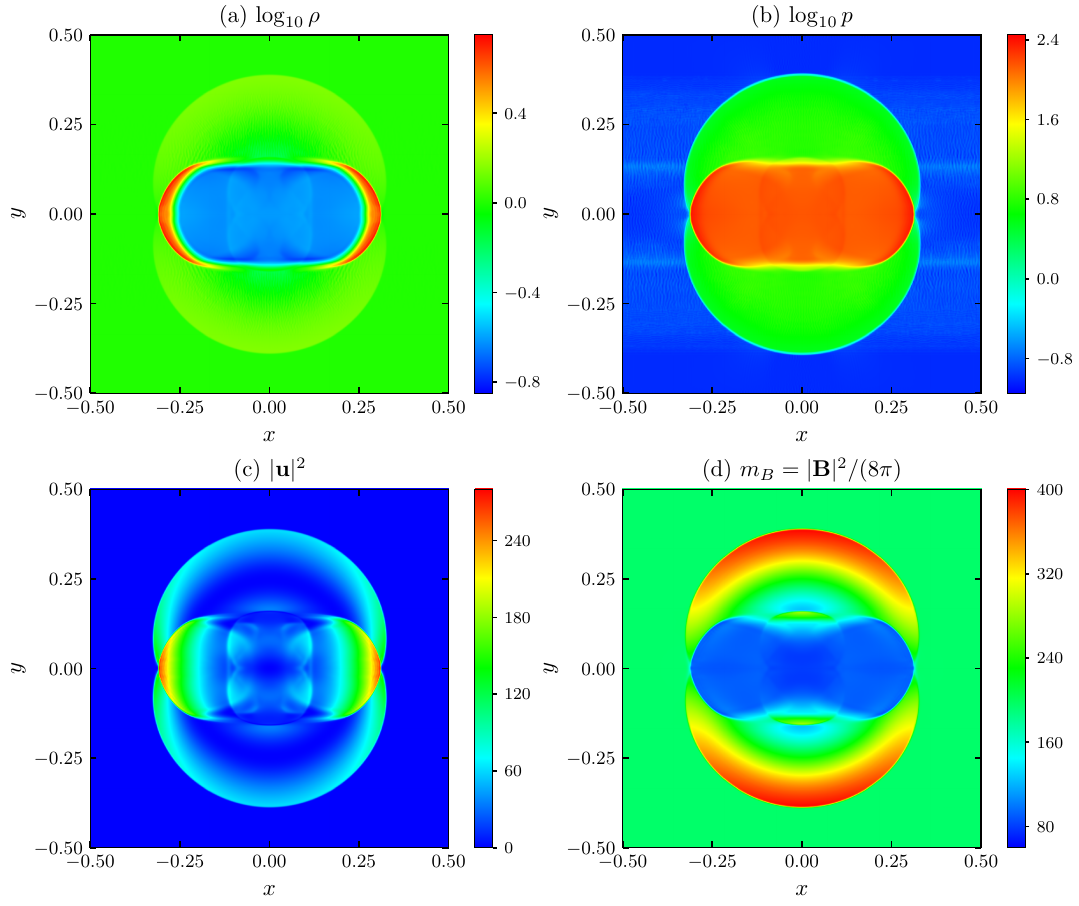}
		\caption{Cylindrical blast wave on $1024^2$ cells at $t_f=0.01$.}
		\label{fig:blast-validation}
	\end{figure}

	The cylindrical blast wave describes an expanding shock in a strong horizontal magnetic field. On $\Omega=[-0.5,0.5]^2$, the initial data are
	\begin{equation}
		\rho=1,\qquad \uvec=\bm0,\qquad \B=(70,0)^\top,\qquad
		p=\begin{cases}1000,&r<0.1,\\0.1,&r\geq0.1,\end{cases}
		\label{eq:blast-data}
	\end{equation}
	with $r=(x^2+y^2)^{1/2}$ and $\gamma=1.4$. As done in \cite{BoscheriThomann2024}, we use $1024^2$ cells, Dirichlet boundary conditions given by the initial exterior state and the final time $t_f=0.01$. The scheme is second order in space and first order in time. The results are shown in \Cref{fig:blast-validation}. The horizontal magnetic field confines the expansion: the central cavity is elongated along the field, with a strong compression at its lateral ends and an increase of the magnetic energy above and below it. The outer front and the symmetries of the solution agree qualitatively with Figure 6 of \cite{BoscheriThomann2024}, while some differences remain in the extrema and in the internal wave structure. We recall that the two computations use different orders of accuracy in time. The simulation reaches $t_f$ in 1797 accepted time steps, despite the initial pressure ratio of $10^4$. Small-scale oscillations remain in the compressed region, but density and pressure stay positive and the normal jumps remain below $3.823\times10^{-12}$.
	
	\subsection{MHD rotor}
	\label{sec:num-rotor}
	
\begin{figure}[!htb]
		\centering
		\includegraphics[width=0.94\textwidth]{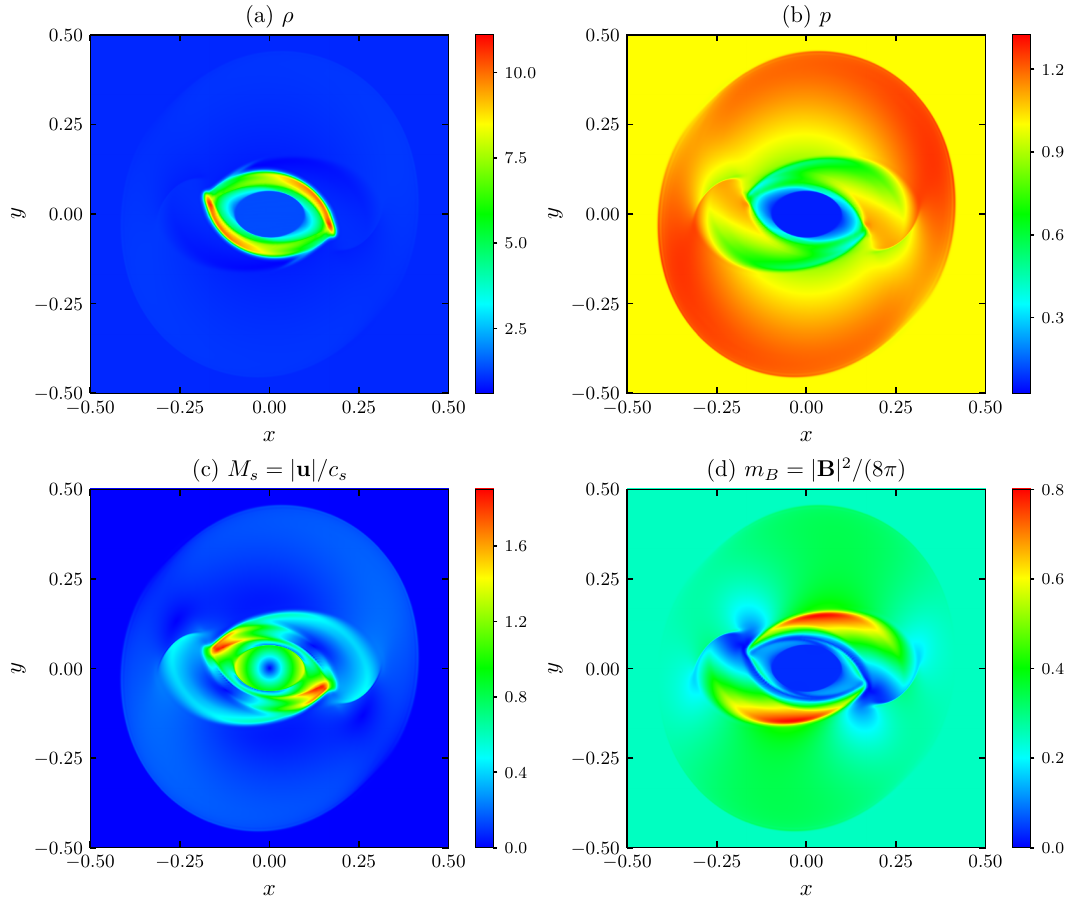}
		\caption{MHD rotor on $1024^2$ cells at $t_f=0.25$.}
		\label{fig:rotor-validation}
	\end{figure}

	A rapidly rotating dense disk is embedded in a light fluid at rest. On $\Omega=[-0.5,0.5]^2$, the pressure and the magnetic field are uniform with $p=1$ and $\B=(2.5,0)^\top$, $\gamma=1.4$, while the initial density and velocity are
	\begin{equation}
		(\rho,\uvec)=
		\begin{cases}
			(10,10(-y,x)^\top),&r\leq0.1,\\
			(1,\bm0),&r>0.1.
		\end{cases}
		\label{eq:rotor-data}
	\end{equation}
	Following \cite{BoscheriThomann2024}, we use $1024^2$ cells, the initial exterior state at the boundaries and the final time $t_f=0.25$. The rotating core bends the magnetic field lines and generates torsional Alfv\'en waves. The elliptic front, the tilted dense annulus, the low pressure in the center and the maxima of the magnetic energy shown in \Cref{fig:rotor-validation} have the same orientation and extent as in Figure 7 of \cite{BoscheriThomann2024}, with some differences in the central profiles and in the peak values. Small oscillations remain near some fronts. The simulation is completed in 1222 time steps without rejection, and the final normal jump is $5.773\times10^{-15}$.
	
	\subsection{Orszag--Tang vortex}
	\label{sec:num-orszag}
	
	The Orszag--Tang vortex \cite{OrszagTang,DahlburgPicone,PiconeDahlburg} starts from smooth periodic data and develops a complex network of interacting fast and slow shocks. On the periodic domain $\Omega=[0,1]^2$, with $\gamma=5/3$, the initial data are
	\begin{equation}
		\rho=\frac{25}{36\pi},\quad
		p=\frac{5}{12\pi},\quad
		\uvec=\binom{-\sin(2\pi y)}{\sin(2\pi x)},\quad
		\B=\frac1{\sqrt{4\pi}}
		\binom{-\sin(2\pi y)}{\sin(4\pi x)}.
		\label{eq:orszag-data}
	\end{equation}
	The simulation is run on $1024^2$ cells up to $t_f=4$.
	\begin{figure}[!htb]
		\centering
		\includegraphics[width=\textwidth]{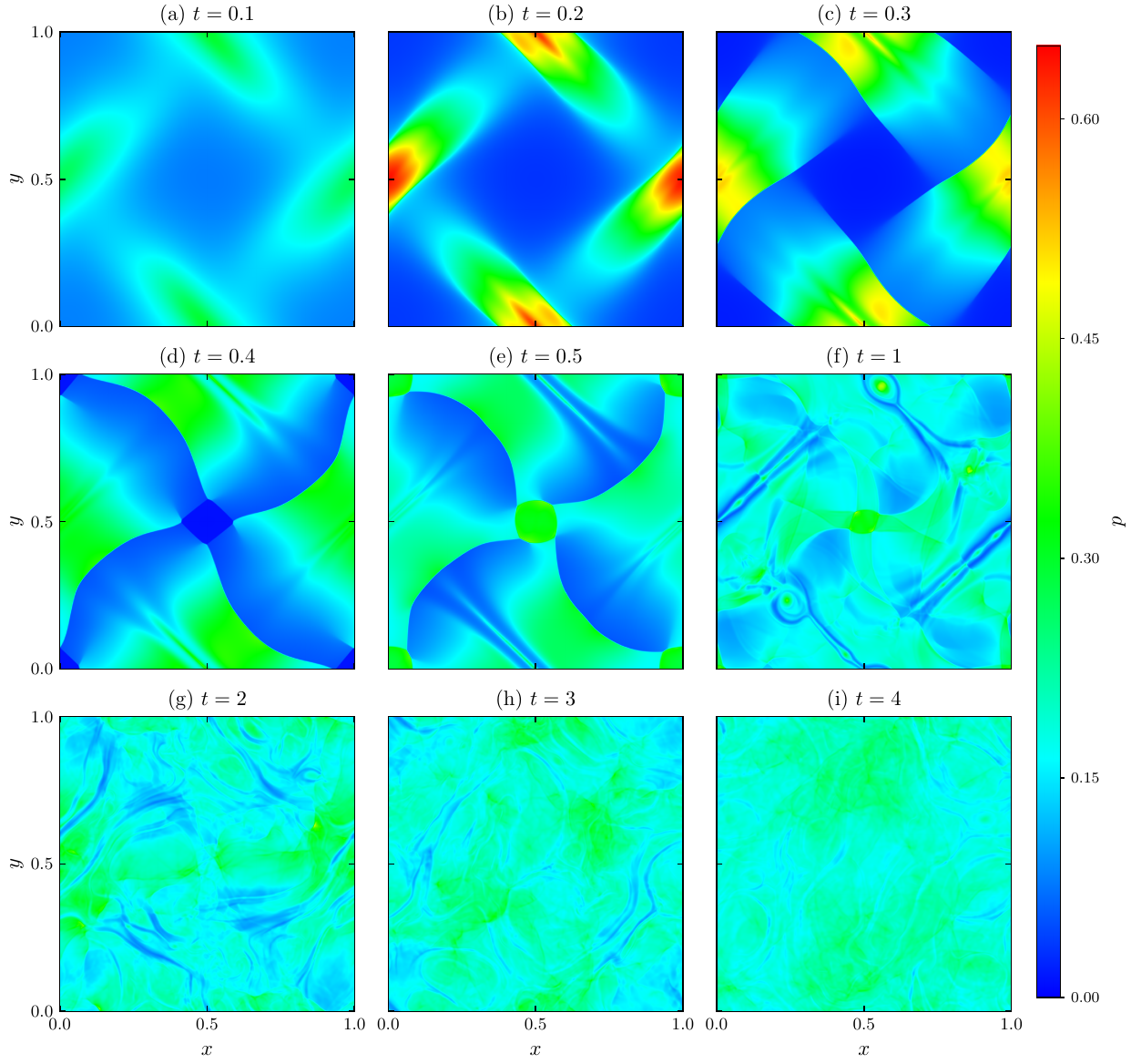}
		\caption{Orszag--Tang vortex on $1024^2$ cells: pressure at different times.}
		\label{fig:orszag-pressure-history}
	\end{figure}
	The pressure in \Cref{fig:orszag-pressure-history} shows how the oblique compressive waves steepen into a network of shocks, which then generate finer vortical structures through repeated interactions. The solution at $t=0.5$, $2$, $3$ and $4$ agrees qualitatively with \cite{BoscheriThomann2024}. The final solution in \Cref{fig:orszag-final-fields} retains the large-scale structure and the periodic symmetries. Over $20366$ time steps, the relative mass and total energy errors \eqref{eq:conservation-errors} remain below $1.897\times10^{-11}$ and $2.135\times10^{-13}$, respectively, and the final normal jump is $1.333\times10^{-15}$.
	
	\begin{figure}[!t]
		\centering
		\includegraphics[width=\textwidth]{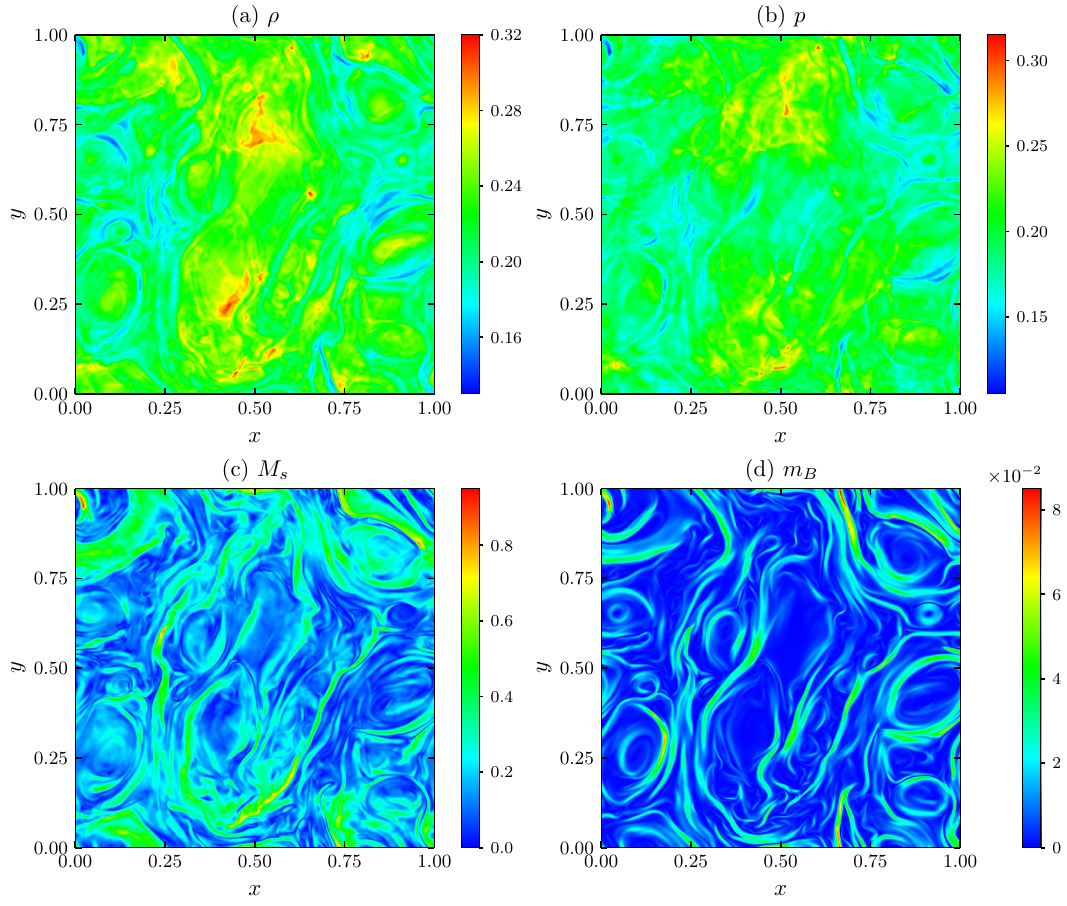}
		\caption{Orszag--Tang vortex on $1024^2$ cells at $t_f=4$.}
		\label{fig:orszag-final-fields}
	\end{figure}
	
	\subsection{Low Mach magnetic field loop advection}
	\label{sec:num-field-loop}
	
	The last test consists in the diagonal transport of a weak magnetic loop through the periodic domain $\Omega=[-1,1]\times[-0.5,0.5]$. The fluid state is uniform,
	\begin{equation}
		\rho=1,\qquad \uvec=(2,1)^\top,\qquad p=10^5,
		\label{eq:field-loop-fluid}
	\end{equation}
	and the magnetic field is given by the curl of the compactly supported potential
	\begin{equation}
		A_z=A_0\max(R-r,0),\qquad R=0.3,
		\qquad \B=\nabla\times(0,0,A_z)^\top,
		\label{eq:field-loop-potential}
	\end{equation}
	that is $\B=A_0(-y/r,x/r)^\top$ for $0<r\leq R$ and $\B=0$ outside the loop. We consider four amplitudes $A_0\in\{10^{-3},10^{-2},10^{-1},1\}$, on $512\times256$ cells up to $t_f=1$. At this time, the exact solution coincides with the initial data. To measure the magnetic dissipation, we compute the ratio between the final and the initial magnetic energy:
	\begin{equation}
		\mathcal R_B(t_f)=
		\frac{\int_\Omega|\B_h(\xvec,t_f)|^2\dd\xvec}
		{\int_\Omega|\B_h(\xvec,0)|^2\dd\xvec}.
		\label{eq:magnetic-retention}
	\end{equation}
	
	\begin{figure}[!t]
		\centering
		\includegraphics[width=\textwidth]{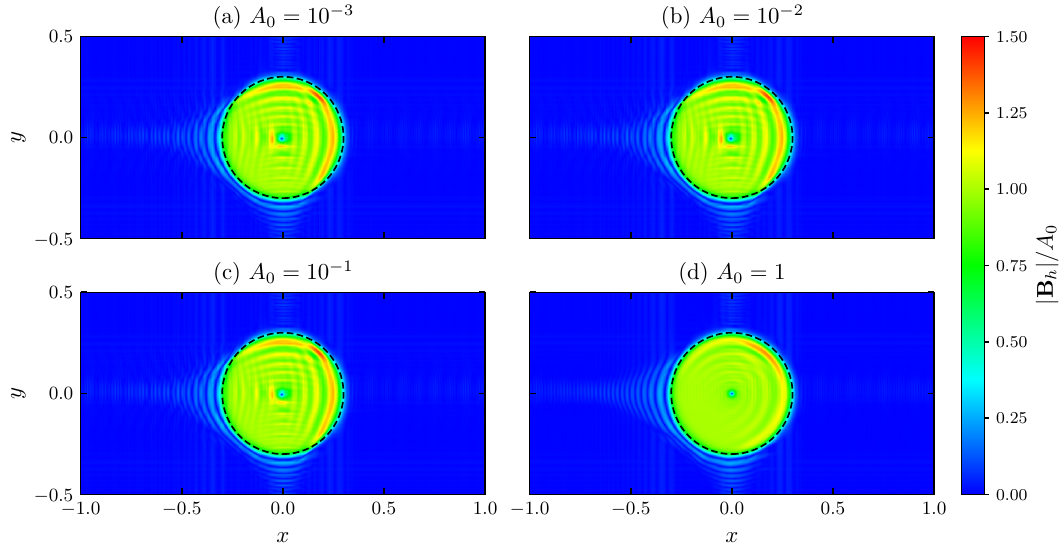}
		\caption{Field loop advection on $512\times256$ cells at $t_f=1$.}
		\label{fig:field-loop}
	\end{figure}
	
	The normalized solutions in \Cref{fig:field-loop} are almost identical for the three smallest amplitudes, and the same circular structure is kept for $A_0=1$. For the four runs, the Alfv\'en Mach number ranges from approximately $7.9\times10^3$ to $7.9$, the ratio $\mathcal R_B$ lies in $[0.9656,0.9682]$ and the number of time steps lies in $[2277,2335]$. Small dispersive ripples are visible near the boundary and the center of the loop, but they do not grow in normalized amplitude. For $A_0=1$, the final density lies in $[0.9940,1.0068]$ and the relative pressure perturbation stays below $2.33\times10^{-6}$. Therefore, both the accuracy and the computational cost depend only weakly on the magnetic amplitude in the tested range.
	
	\subsection{Preservation of the involution constraints}
	\label{sec:num-cross-test}
	
\begin{figure}[!htb]
		\centering
		\includegraphics[width=\textwidth]{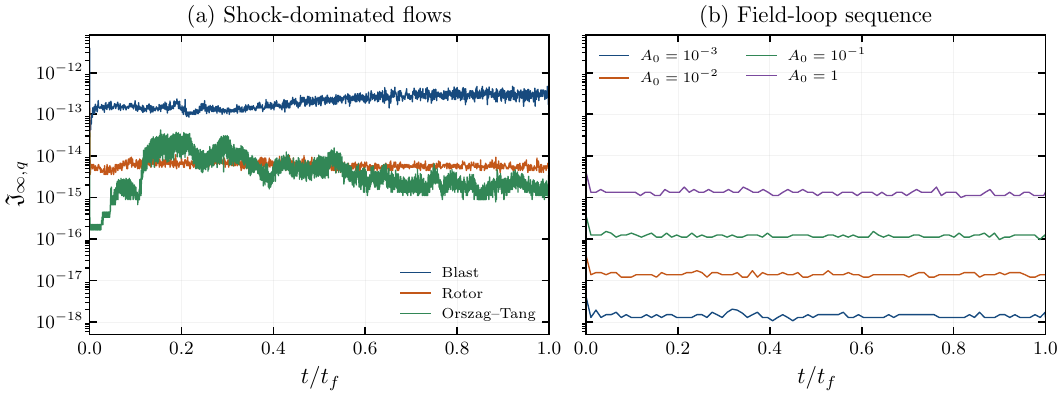}
		\caption{Time evolution of the normal jumps $\frakJ_{\infty,q}$ for the shock-dominated tests and the field loops.}
		\label{fig:diagnostic-histories}
	\end{figure}

	Finally, we summarize the behavior of the two magnetic involutions in all tests. The local divergence $\frakD_\infty$ is always exactly zero, while the normal jumps stay at the accuracy of the projection. As shown in \Cref{fig:diagnostic-histories}, $\frakJ_{\infty,q}$ remains below $3.83\times10^{-12}$ for the blast wave, the rotor, the Orszag--Tang vortex and the four field loops. The maximum value is reached in the blast wave, and it is smaller than $5.5\times10^{-14}$ relative to the magnitude of the initial magnetic field. For the field loops, the maximum normal jump ranges from $4.140\times10^{-18}$ for $A_0=10^{-3}$ to $4.052\times10^{-15}$ for $A_0=1$, that is approximately $4\times10^{-15}$ after normalization by $A_0$. These results confirm that both the local and the global divergence-free properties are preserved in the presence of shocks, magnetic stabilization and low Mach transport.
	
	\section{Conclusions}
	\label{sec:conclusions}

	In this work, we have presented a semi-implicit FV/DG scheme for ideal MHD. The magnetic field is represented by divergence-free basis functions, which preserve the local involution throughout the update. A constrained projection enforces the continuity of its normal component across cell faces. The two constructions provide local solenoidality and global compatibility within the same magnetic approximation. Material transport is treated explicitly, while the implicit magnetic and pressure sub-systems lead to sequential linear solves. The scheme is thus suitable for low Mach and Alfvén number flows, with a time step restriction based on the material speed of the fluid. Subcell stress balances supply the force entering the magnetic operator, and Richardson extrapolation gives second order temporal accuracy for smooth solutions.

	The numerical tests cover smooth vortices, Riemann problems, strongly nonlinear two-dimensional flows and low Mach magnetic-loop advection. The reported mesh-refinement studies show decreasing cell-average errors across the tested Mach and Alfv\'en regimes. The local divergence vanishes by construction, and the normal jumps remain controlled by the projection throughout the simulations, including those with shocks and limiting.

	A direction for future work is to combine this approach with the involution-preserving ALE framework based on constrained mesh motion developed in \cite{BernardelliBoscheriCardonnaALE}. The aim is to use locally divergence-free polynomial bases of arbitrary degree in a high-order DG discretization and determine the mesh motion so as to enforce continuity of the normal magnetic trace across interfaces. This would couple the local algebraic constraint with a geometric treatment of the interface constraint.

	\section*{Acknowledgments}
	The authors would like to thank the framework offered by the SUN HYPE 2026 event organized in Crete in June 2026, where the work presented in this paper was conceived and its initial development took place. In particular, we thank the Technical University of Crete for hosting the event and the financial support received from the European Union through the ERC Starting Grant ALcHyMiA (grant agreement No. 101114995).
	
	WB also acknowledges funding by the Institute de Mathématiques pour la Planète Terre (IMPT - France) under the project AAP2025-ERUPTA. WB is member of the Gruppo Nazionale per il Calcolo Scientifico dell'Istituto Nazionale di Alta Matematica (GNCS-INdAM).
	
	LG also acknowledges funding by the Spanish Ministry of Science, Innovation and Universities (MICIU), the Spanish AEI (MCIN/AEI/10.13039/501100011033) and European Social Fund Plus under the project RYC2022-036355-I, by Xunta de Galicia under Axudas Pro-ERC Oportunius and by the European Research Council under ERC 2025 Starting Grant - 101219299 SUPREMUM.
	
	Views and opinions expressed are however those of the author(s) only and do not necessarily reflect those of the European Union or the European Research Council Executive Agency. Neither the European Union nor the granting authority can be held responsible for them.
	
	\section*{Data availability}
	Data and material are available upon reasonable request addressed
	to the corresponding author.
	
	\section*{Declaration of competing interest}
	The authors declare that they have no known competing financial interests or personal relationships that could have appeared to influence the work reported in this paper.

	\appendix
	\renewcommand{\theHtable}{\thesection.\arabic{table}}
	
	\section{Compact finite volume reconstruction}
	\label{app:reconstruction-coefficients}
	
	\paragraph{CWENO reconstruction}
	For each conservative variable $w$, the reconstruction in cell $\omega_c$ is a linear polynomial of the form
	\begin{equation}
		w_c(\xvec)=w_c+\bm g_c\cdot(\xvec-\xvec_c),
		\label{eq:linear-reconstruction}
	\end{equation}
	which is obtained with a CWENO strategy \cite{dumbser2017central,cravero2018cweno}. Let $w_c^{\mathrm{opt}}$ be the centered linear polynomial, and $w_{c,s}$, $s=1,\ldots,4$, the four linear candidates given by the vertex patch gradients of \Cref{app:periodic-reconstruction}. All candidates have mean $w_c$. The central candidate $w_{c,0}$ is defined such that
	\begin{equation}
		w_c^{\mathrm{opt}}=\sum_{s=0}^4\lambda_s w_{c,s},
		\label{eq:cweno-optimal-recovery}
	\end{equation}
	with positive normalized linear weights $\lambda_s$. The reconstruction is then given by the nonlinear combination $w_c(\xvec)=\sum_s\varpi_{c,s}w_{c,s}(\xvec)$, with
	\begin{equation}
		\varpi_{c,s}=\frac{\alpha_{c,s}}{\sum_j\alpha_{c,j}},
		\qquad
		\alpha_{c,s}=\frac{\lambda_s}{(I_{c,s}+\epsilon)^r},
		\label{eq:cweno-weights}
	\end{equation}
	where $I_{c,s}=|\partial_\xi w_{c,s}|^2+|\partial_\eta w_{c,s}|^2$, $r=4$ and $\epsilon=10^{-12}$. The unnormalized linear weights are $10^5$ for the central candidate and $1$ for the four peripheral ones, and they are the same in all tests.
	
	\paragraph{Troubled cells and positivity}
	Let $\eta_h=(\varrho_{\min}/L_{\max})^{1/2}$, and let $w_{\min}$ be the floor value of $w\in\{\rho,p\}$. A cell is marked as troubled if
	\begin{equation}
		\max_{w\in\{\rho,p\}}\max_{v\in\mathcal N(c)}
		\frac{|w_v-w_c|}{\max(w_{\min},|w_v|+|w_c|)}
		>\eta_h ,
		\label{eq:troubled-cell-indicator}
	\end{equation}
	and its face neighbors are marked as well. In the marked cells, the CWENO gradient is replaced by the minmod slopes \cite{Harten1983}, where minmod returns the argument of smallest magnitude if both have the same sign and zero otherwise:
	\begin{equation}
		(\bm g_c^{\mathrm{lim}})_d=
		\frac{1}{\varrho_d}\operatorname{minmod}
		(w_{c+e_d}-w_c,\;w_c-w_{c-e_d}),
		\qquad d\in\{x,y\},
		\label{eq:tvd-gradient}
	\end{equation}
	where $c\pm e_d$ denote the neighbors of cell $\omega_c$ in direction $d$. Finally, the reconstruction $\Q_h-\Q_c$ is scaled in each cell so that density and pressure are positive at the vertices and at the quadrature points \cite{ZhangShu2010}. We remark that this only guarantees the positivity of the reconstructed states, and not of the whole time step. The magnetic components of this reconstruction are denoted by $\mathcal R_h^B\overline\B$; they only enter the numerical dissipation of \eqref{eq:explicit-induction-flux} for $\sigma=1$, while the central part of the induction flux always uses $\B_h$.
	
	\subsection{Magnetic relaxation rate}
	
	The rate $\nu_c$ of \eqref{eq:zero-mean-modal-operator} is computed from the cell means of the magnetic field at time $t^n$. For each component $\overline B_{j,c}$ and direction $d$, let $\delta^-=\overline B_{j,c}-\overline B_{j,c-e_d}$, $\delta^+=\overline B_{j,c+e_d}-\overline B_{j,c}$ and $\delta^0=(\delta^-+\delta^+)/2$. The smoothness indicators are
	\begin{equation}
		\beta_0=(\delta^+-\delta^-)^2+(\delta^0)^2+\epsilon_B,\qquad
		\beta_\pm=(\delta^\pm)^2+\epsilon_B,
		\label{eq:magnetic-cweno-indicators}
	\end{equation}
	with $\epsilon_B=10^{-12}$. With the linear weights $\lambda_0=1/2$ and $\lambda_-=\lambda_+=1/4$, the nonlinear weights are
	\begin{equation}
		\alpha_s=\lambda_s\left[1+
		\left(\frac{|\beta_--\beta_+|}{\beta_s}\right)^2\right],
		\qquad
		\varpi_s=\frac{\alpha_s}{\alpha_0+\alpha_-+\alpha_+},
		\quad s\in\{0,-,+\}.
		\label{eq:magnetic-cweno-weights}
	\end{equation}
	The distance from the linear weights is measured by $\zeta_{d,j,c}=\frac23\sum_s|\varpi_s-\lambda_s|$, and we set $\zeta_{d,c}=\max_{j\in\{x,y\}}\zeta_{d,j,c}$. The relaxation rate is finally
	\begin{equation}
		\nu_c=\max_{d\in\{x,y\}}\frac{c_{B,d,c}}{h_d}\zeta_{d,c}^{\,2},
		\qquad
		c_{B,d,c}=\frac12\left(|u_{d,c}|+
		\sqrt{u_{d,c}^2+\frac{4|\overline\B_c|^2}{4\pi\rho_c}}\right).
		\label{eq:magnetic-modal-rate}
	\end{equation}
	If one of the two neighbors in direction $d$ is missing at a physical boundary, the contribution of that direction is set to zero.
	
	\paragraph{Remarks}
	Let us notice that $\nu_c$ depends only on the cell means, thus it does not control all polynomial modes. For instance, on an even periodic mesh with unit cells, the field $\B_h=(a_{ij}X,b_i-a_{ij}Y)^\top$, with $X=x-x_c$, $Y=y-y_c$, $a_{ij}=(-1)^{i+j}$, $(b_0,b_1,b_2,b_3)=(0,0,1,0)$ and $\uvec=0$, has continuous normal traces, but in the cells of column one it gives $\overline\B_c=0$ and $\nu_c=0$, so its zero-mean part is not damped. Furthermore, for smooth data the rate is small: a one-dimensional Taylor expansion gives $\nu=O(h^5)$ away from extrema and $\nu=O(h^3)$ at a generic extremum. However, these estimates are not uniform with respect to the position of the extremum inside the cell, and $\nu=O(h^{-1})$ may occur as long as $\epsilon_B\ll h^4|f''|^2$, before the regularization $\epsilon_B$ becomes dominant. Therefore, \eqref{eq:modal-dissipation} alone guarantees neither uniform consistency for smooth data nor coercivity on all modes.
	
	\subsection{Reference face-moment matrices}
	\label{app:constraint-blocks}
	
	With the basis ordered as in \eqref{eq:reference-basis}, the constraint blocks on the left ($L$), right ($R$), bottom ($B$) and top ($T$) faces of the reference cell, computed with the outward normals, are
	\begin{equation}
		\begin{aligned}
			\widehat{\mathds C}_{L}
			&=\begin{pmatrix}
				0&-1&0&0&-1\\ 0&0&0&0&-1/6
			\end{pmatrix},
			&\widehat{\mathds C}_{R}
			&=\begin{pmatrix}
				0&1&0&1&1\\ 0&0&0&0&1/6
			\end{pmatrix},\\[2mm]
			\widehat{\mathds C}_{B}
			&=\begin{pmatrix}
				1&0&1&0&0\\ 0&0&1/6&0&0
			\end{pmatrix},
			&\widehat{\mathds C}_{T}
			&=\begin{pmatrix}
				-1&0&-1&-1&0\\ 0&0&-1/6&0&0
			\end{pmatrix}.
		\end{aligned}
		\label{eq:reference-constraint-blocks}
	\end{equation}
	The blocks of the physical cells are obtained with the Piola transform \eqref{eq:piola}. For an interior face, the negatives of the outward-normal blocks of the two adjacent cells are inserted in the same two rows of $\Ch$, consistently with the orientation of the jump in \eqref{eq:gothic-face-involution}.
	
	\section{Boundary treatment and solution of the linear systems}
	\label{app:multidimensional-consistency}
	
	\subsection{Periodic vertex patches}
	\label{app:periodic-reconstruction}
	
	Let $\mathcal C(p)$ be the set of cells sharing the vertex $p$, including the periodic images if $p$ lies on a periodic boundary. The vertex control volume is $\omega_p=\bigcup_{c\in\mathcal C(p)}\omega_{pc}$, where $\omega_{pc}$ is the quarter of cell $\omega_c$ attached to $p$. Let $\bm\ell_{pc}$ be the sum of the outward normals of the two edges of $\omega_{pc}$ incident to $p$, each one weighted by half of the edge length. The patch gradient is
	\begin{equation}
		\bm g_p(w)=-\frac1{|\omega_p|}
		\sum_{c\in\mathcal C(p)}
		w_c\,\bm\ell_{pc}.
		\label{eq:periodic-patch-gradient}
	\end{equation}
	At physical boundaries, the corresponding boundary closure is used.
	
	\subsection{Boundary conditions}
	\label{app:boundary-traces}
	
	The prescribed boundary data $\Q_b(\xvec,t)$ are evaluated at the face quadrature points, i.e. $\Q^+_{f,q}=\Q_b(\xvec_{fq},t)$, and they enter the convective flux, the explicit induction flux, the magnetic energy flux and the enthalpy flux on the right hand side of \eqref{eq:pressure-system}. For the implicit terms, we use the values of the interior cell:
	\begin{equation}
		\widehat{\bm f}_{B,f}^{\,n+1}
		=\mathds F_{B,c}^{n+1}\nvec_{cf},\qquad
		\widehat{\mathds T}_{B,f}=\mathds T_{B,c},\qquad
		E^{n+1}_f=E_c^{n+1}.
		\label{eq:physical-correction-traces}
	\end{equation}
	The last relation is used in the pressure gradient \eqref{eq:pressure-momentum}, while the boundary contribution to $\mathcal L_H$ is zero. In the enthalpy flux, the boundary value is $\tfrac12(h_c^{n}\mvec_c^{**}+h_b\mvec_b)$, with $h_b=\gamma p_b/\bigl((\gamma-1)\rho_b\bigr)$. These closures do not affect the normal magnetic trace, which is not constrained by the projection on physical boundaries. In the Riemann problems, the left and right states are prescribed at the vertical boundaries, while for the blast wave and the rotor the initial exterior state is imposed. Periodic boundaries are handled through the neighbor map.
	
	\subsection{Linear solvers and time step acceptance}
	\label{app:acceptance-criteria}
	
	\paragraph{Linear solvers}
	The magnetic system \eqref{eq:magnetic-linear-system} is not symmetric and it is solved with preconditioned GMRES \cite{SaadSchultz1986}. The preconditioner combines the element mass matrices with the exact inverse of the stabilization term,
	\begin{equation}
		(\I+\dt\,\nu_c\mathds P_0)^{-1}
		=\I-\frac{\dt\,\nu_c}{1+\dt\,\nu_c}\mathds P_0 ,
		\label{eq:modal-block-inverse}
	\end{equation}
	which follows from $\mathds P_0^2=\mathds P_0$. The energy system \eqref{eq:pressure-system} is solved with preconditioned conjugate gradients. On Cartesian grids we use a Fourier preconditioner for the Helmholtz operator with constant enthalpy, otherwise a diagonal preconditioner. For the projection \eqref{eq:schur-projection}, the rows are normalized in the dual mass norm. On periodic Cartesian meshes, the Schur complement is block diagonal in Fourier space, with $4\times4$ Hermitian blocks for the two moments on the two face orientations, which are inverted with pseudoinverses. On non-periodic grids, an odd extension to a periodic grid of double size provides a Fourier preconditioner, while small systems can be solved with a sparse factorization.
	
	\paragraph{Tolerances and floors}
	The density and pressure floors are $10^{-10}$, and the pressure checked for acceptance is $\overline p_c$ from \eqref{eq:cell-polynomial-pressure}. The pressure and the enthalpy entering the implicit coefficients are floored at $10^{-10}$ and $10^{-12}$, respectively, while the pressure diagnostics are never replaced by the floor. The magnetic and energy solvers target a residual $10^{-14}\max(1,\|b\|_2)$, where $b$ is the right hand side and the magnetic residual is measured after left preconditioning, with a GMRES restart dimension of $30$, at most $5N_c$ restart cycles, and at most $1000$ conjugate gradient iterations. The projection targets $10^{-12}\max(1,\|b\|_2)$, with at most $\max(500,8(N_x+N_y))$ iterations.
	
	\paragraph{Time step rejection}
	A time step is rejected if a normalized residual exceeds $10^{-11}$ for the magnetic or energy solver or $10^{-10}$ for the projection, or if a cell state is not admissible. In that case, the solution is reset to $\Q^n$ and the time step is halved. This reduced value is kept as an upper bound for $\dt$, and it is increased by a factor $1.25$ after eight consecutive accepted steps. The simulation stops after twelve rejections at the same time level, one hundred rejections in total, or if $\dt$ becomes smaller than $\max(100\epsilon_{\mathrm{mach}}\max(1,|t_f|),10^{-8}h_{\min})$. Once the troubled-cell indicator has been activated, and thus $\sigma=1$, it stays active even if the time step is rejected. We remark that these criteria do not guarantee the positivity or the stability of the scheme.
	
	\section{Temporal accuracy}
	\label{app:coupled-time-order}
	\setcounter{table}{0}
	
	To verify the second order of accuracy in time, we consider a periodic $16\times16$ mesh, with background state $\rho=p=1$ and oblique magnetic field $\B=\sqrt{4\pi}(1,0.7)^\top$. The smooth perturbation
	\begin{equation}
		u(x,y,0)=10^{-2}\sin(2\pi x)\sin(2\pi y),\qquad
		v(x,y,0)=0.6\cdot10^{-2}\cos(2\pi x)\sin(2\pi y)
		\label{eq:coupled-temporal-probe}
	\end{equation}
	is evolved on the fixed mesh with a constant time step. In this test, acoustic, magnetic and transport effects are coupled, the projection is active at every time step, and neither limiting nor time step rejection occur. The orders of convergence reported in \Cref{tab:coupled-temporal-probe} are close to two for all variables.
	
	\begin{table}[!htbp]
		\centering
		\footnotesize
		\setlength{\tabcolsep}{4pt}
		\caption{Orders of convergence in time of the second order scheme \eqref{eq:richardson-composition} on a fixed mesh, computed from the root-mean-square differences of the cell averages at the final time.}
		\label{tab:coupled-temporal-probe}
		\begin{tabular*}{\linewidth}{@{\extracolsep{\fill}}ccccccc@{}}
			\toprule
			number of time steps & $\rho$ & $m_x$ & $m_y$ & $E$ & $B_x$ & $B_y$\\
			\midrule
			$20,40,80$   & 2.0022 & 1.9740 & 1.9141 & 1.9835 & 1.9960 & 1.9948\\
			$40,80,160$  & 2.0007 & 2.0347 & 2.1106 & 1.9882 & 2.0297 & 2.0215\\
			$80,160,320$ & 1.9986 & 2.0200 & 2.0250 & 1.9788 & 2.0327 & 2.0407\\
			\bottomrule
		\end{tabular*}
	\end{table}
	
	\subsection{Projection and order of accuracy}
	\label{app:stagewise-projection-order}
	
	Here we explain why the projection is combined with a Richardson extrapolation rather than with a Runge--Kutta method. Let us consider a linear magnetic model with generator $A$ and projection $\mathds P_B$. For a compatible input, one projected implicit step with step size $\tau$ writes
	\begin{equation}
		\mathds P_B(\I-\tau A)^{-1}\mathds P_B
		=\mathds P_B+\tau\mathds P_BA\mathds P_B
		+\tau^2\mathds P_BA^2\mathds P_B+O(\tau^3),
		\label{eq:projected-stage-expansion}
	\end{equation}
	and the second order term differs from the square of the projected generator by
	\begin{equation}
		\mathds P_BA^2\mathds P_B-(\mathds P_BA\mathds P_B)^2
		=\mathds P_BA(\I-\mathds P_B)A\mathds P_B.
		\label{eq:projection-commutator-defect}
	\end{equation}
	This error vanishes if $A$ maps compatible fields into compatible fields, but not in general. If the projection is applied after each stage of the LSDIRK2 scheme used in \cite{BoscheriThomann2024}, this term gives an error $\alpha\dt^2\mathds P_BA(\I-\mathds P_B)A\mathds P_B$, with $\alpha=1-1/\sqrt2$, which can reduce the order of accuracy to one. Indeed, in the test of \Cref{app:coupled-time-order} with perturbation amplitude $10^{-4}$, the finest time steps give orders $0.994$, $0.988$, $1.071$ and $1.082$ for $m_x$, $m_y$, $B_x$ and $B_y$, respectively. On the contrary, the Richardson extrapolation of the complete first order step cancels the second order error, as shown by \eqref{eq:richardson-consistency}. We remark that this problem does not occur in \cite{BoscheriThomann2024}, since the magnetic field is computed as the discrete curl of a vector potential, and no projection is needed between the stages.

	\bibliographystyle{elsarticle-num}
	\bibliography{references}

@Article{Toro_Vazquez12,
  author   = {E.F. Toro and M.E. V{\'a}zquez-Cend{\'o}n},
  journal  = {Comput. \& Fluids},
  title    = {{Flux splitting schemes for the Euler equations}},
  year     = {2012},
  pages    = {1--12},
  volume   = {70},
  fjournal = {Computers \& Fluids},
}

@Article{GH,
  author   = {H. Guillard and A. Murrone},
  journal  = {Comput. \& Fluids},
  title    = {{On the behavior of upwind schemes in the low Mach number limit : II. Godunov type schemes}},
  year     = {2004},
  pages    = {655--675},
  volume   = {33},
  fjournal = {Computers \& Fluids},
}

@Article{Dellacherie1,
  author   = {S. Dellacherie},
  journal  = {J. Comput. Phys.},
  title    = {{Analysis of Godunov type schemes applied to the compressible Euler system at low Mach number}},
  year     = {2010},
  pages    = {978--1016},
  volume   = {229},
  fjournal = {Journal of Computational Physics},
}

@Article{BosFil2016,
  author   = {S. Boscarino and F. Filbet and G. Russo},
  journal  = {J. Sci. Comput.},
  title    = {{High Order Semi-implicit Schemes for Time Dependent Partial Differential Equations}},
  year     = {2016},
  pages    = {975--1001},
  volume   = {68},
  fjournal = {Journal of Scientific Computing},
}

@Article{AscRuuSpi,
  author   = {U.~M. Ascher and S.~J. Ruuth and R.~J. Spiteri},
  journal  = {Appl. Numer. Math.},
  title    = {{Implicit-explicit {R}unge-{K}utta methods for time-dependent partial differential equations}},
  year     = {1982},
  pages    = {151--167},
  volume   = {25},
  fjournal = {Applied Numerical Mathematics},
}

@Article{Busto2021,
  author   = {S. Busto and L. Rio-Martin and M.E. Vazquez-Cendon and M. Dumbser},
  journal  = {Appl. Math. Comput.},
  title    = {A semi-implicit hybrid finite volume/finite element scheme for all Mach number flows on staggered unstructured meshes},
  year     = {2021},
  volume   = {402},
  fjournal = {Applied Mathematics and Computation},
}

@Article{Chen2022,
  author  = {W. Chen and K. Wu and T. Xiong},
  journal = {Journal Of Computational Physics},
  title   = {{High order asymptotic preserving finite difference WENO schemes with constrained transport for MHD equations in all sonic Mach numbers}},
  year    = {2023},
  pages   = {112240},
  volume  = {488},
}

@Article{yee1966,
  author  = {Yee, K.},
  journal = {{IEEE Transactions on antennas and propagation}},
  title   = {Numerical solution of initial boundary value problems involving {Maxwell}'s equations in isotropic media},
  year    = {1966},
  number  = {3},
  pages   = {302--307},
  volume  = {14},
}

@Article{dedner2002,
  author   = {Dedner, A. and Kemm, F. and Kr{\"o}ner, D. and Munz, C.-D. and Schnitzer, T. and Wesenberg, M.},
  journal  = {J. Comput. Phys.},
  title    = {{Hyperbolic divergence cleaning for the MHD equations}},
  year     = {2002},
  number   = {2},
  pages    = {645--673},
  volume   = {175},
  fjournal = {Journal of Computational Physics},
}

@Article{Godunov1959,
  author  = {S. Godunov},
  journal = {Mat. Sb.},
  title   = {{Finite difference methods for the computation of discontinuous solutions of the equations of fluid dynamics}},
  year    = {1959},
  pages   = {271--306},
  volume  = {47},
}

@article{TavelliDumbser2017AllMach,
  author = {Tavelli, Maurizio and Dumbser, Michael},
  title = {A pressure-based semi-implicit space--time discontinuous Galerkin method on staggered unstructured meshes for the solution of the compressible Navier--Stokes equations at all Mach numbers},
  journal = {Journal of Computational Physics},
  volume = {341},
  pages = {341--376},
  year = {2017},
  doi = {10.1016/j.jcp.2017.03.030}
}

@article{ZampaBustoDumbser2024,
  author = {Zampa, Enrico and Busto, Saray and Dumbser, Michael},
  title = {A divergence-free hybrid finite volume / finite element scheme for the incompressible MHD equations based on compatible finite element spaces with a posteriori limiting},
  journal = {Applied Numerical Mathematics},
  volume = {198},
  pages = {346--374},
  year = {2024},
  doi = {10.1016/j.apnum.2024.01.014}
}

@article{LiXu2012GlobalDG,
  author = {Li, Fengyan and Xu, Liwei},
  title = {Arbitrary order exactly divergence-free central discontinuous Galerkin methods for ideal MHD equations},
  journal = {Journal of Computational Physics},
  volume = {231},
  number = {6},
  pages = {2655--2675},
  year = {2012},
  doi = {10.1016/j.jcp.2011.12.016}
}

@article{LiXu2013LocalDG,
  author = {Li, Fengyan and Xu, Liwei},
  title = {Locally divergence-free central discontinuous Galerkin methods for ideal MHD equations},
  journal = {Journal of Computational Science},
  volume = {4},
  number = {1--2},
  pages = {80--91},
  year = {2013},
  doi = {10.1016/j.jocs.2012.05.002}
}

@article{BustoDumbser2023,
  author={S. Busto and M. Dumbser},
  title={A new thermodynamically compatible finite volume scheme for magnetohydrodynamics},
  journal={SIAM Journal on Numerical Analysis}, year={2023},
  doi={10.1137/22M147815X}}

@article{BoscheriThomann2024,
  author={W. Boscheri and A. Thomann},
  title={A structure-preserving semi-implicit {IMEX} finite volume scheme for ideal magnetohydrodynamics at all {Mach} and {Alfv\'en} numbers},
  journal={Journal of Scientific Computing}, volume={100}, pages={67}, year={2024},
  doi={10.1007/s10915-024-02606-1}}

@article{BoscheriDimarcoPareschi2023,
  author={W. Boscheri and G. Dimarco and L. Pareschi},
  title={Locally structure-preserving div-curl operators for high order discontinuous {Galerkin} schemes},
  journal={Journal of Computational Physics}, volume={486}, pages={112130}, year={2023},
  doi={10.1016/j.jcp.2023.112130}}

@article{DumbserBalsaraTavelliFambri2019,
  author={M. Dumbser and D. S. Balsara and M. Tavelli and F. Fambri},
  title={A divergence-free semi-implicit finite volume scheme for ideal, viscous, and resistive magnetohydrodynamics},
  journal={International Journal for Numerical Methods in Fluids}, volume={89}, pages={16--42}, year={2019},
  doi={10.1002/fld.4681}}

@article{Fambri2021,
  author={F. Fambri}, title={A novel structure preserving semi-implicit finite volume method for viscous and resistive magnetohydrodynamics},
  journal={International Journal for Numerical Methods in Fluids}, volume={93}, pages={3447--3489}, year={2021},
  doi={10.1002/fld.5041}}

@article{FambriEtAl2023,
  author={F. Fambri and E. Zampa and S. Busto and L. R\'io-Mart\'in and F. Hindenlang and E. Sonnendr\"ucker and M. Dumbser},
  title={A well-balanced and exactly divergence-free staggered semi-implicit hybrid finite volume/finite element scheme for the incompressible {MHD} equations},
  journal={Journal of Computational Physics}, volume={493}, pages={112493}, year={2023},
  doi={10.1016/j.jcp.2023.112493}}

@article{BirkeBoscheriKlingenberg2024,
  author={C. Birke and W. Boscheri and C. Klingenberg},
  title={A well-balanced semi-implicit {IMEX} finite volume scheme for ideal magnetohydrodynamics at all {Mach} numbers},
  journal={Journal of Scientific Computing}, volume={98}, pages={34}, year={2024},
  doi={10.1007/s10915-023-02422-z}}

@article{BoscheriPareschi2021,
  author={W. Boscheri and L. Pareschi}, title={High order pressure-based semi-implicit {IMEX} schemes for the 3D {Navier--Stokes} equations at all {Mach} numbers},
  journal={Journal of Computational Physics}, volume={434}, pages={110206}, year={2021},
  doi={10.1016/j.jcp.2021.110206}}

@article{BoscheriTavelli2022,
  author={W. Boscheri and M. Tavelli}, title={High order semi-implicit schemes for viscous compressible flows in 3D},
  journal={Applied Mathematics and Computation}, volume={434}, pages={127457}, year={2022},
  doi={10.1016/j.amc.2022.127457}}

@article{BrezziDouglasMarini1985,
  author={F. Brezzi and J. Douglas, Jr. and L. D. Marini},
  title={Two families of mixed finite elements for second order elliptic problems},
  journal={Numerische Mathematik}, volume={47}, pages={217--235}, year={1985},
  doi={10.1007/BF01389710}}

@article{ZhangShu2010,
  author={X. Zhang and C.-W. Shu},
  title={On positivity-preserving high order discontinuous {Galerkin} schemes for compressible {Euler} equations on rectangular meshes},
  journal={Journal of Computational Physics}, volume={229}, number={23}, pages={8918--8934}, year={2010},
  doi={10.1016/j.jcp.2010.08.016}}

@article{AvgerinosEtAl2019,
  author={S. Avgerinos and F. Bernard and A. Iollo and G. Russo},
  title={Linearly implicit all {Mach} number shock capturing schemes for the {Euler} equations},
  journal={Journal of Computational Physics}, volume={393}, pages={278--312}, year={2019},
  doi={10.1016/j.jcp.2019.04.020}}

@article{Harten1983,
  author={A. Harten},
  title={High resolution schemes for hyperbolic conservation laws},
  journal={Journal of Computational Physics}, volume={49}, pages={357--393}, year={1983},
  doi={10.1016/0021-9991(83)90136-5}}

@article{HartenLaxVanLeer1983,
  author={A. Harten and P. D. Lax and B. van Leer},
  title={On upstream differencing and {Godunov}-type schemes for hyperbolic conservation laws},
  journal={SIAM Review}, volume={25}, number={1}, pages={35--61}, year={1983},
  doi={10.1137/1025002}}

@article{BalsaraSpicer1999,
  author={D. S. Balsara and D. S. Spicer}, title={A staggered mesh algorithm using high order {Godunov} fluxes to ensure solenoidal magnetic fields in magnetohydrodynamic simulations},
  journal={Journal of Computational Physics}, volume={149}, pages={270--292}, year={1999}}

@article{Balsara2004,
  author={D. S. Balsara}, title={Second-order accurate schemes for magnetohydrodynamics with divergence-free reconstruction},
  journal={The Astrophysical Journal Supplement Series}, volume={151}, pages={149--184}, year={2004}}

@article{HymanShashkov1997,
  author={J. M. Hyman and M. Shashkov}, title={Natural discretizations for the divergence, gradient, and curl on logically rectangular grids},
  journal={Computers \& Mathematics with Applications}, volume={33}, pages={81--104}, year={1997}}

@article{dumbser2017central,
  author={M. Dumbser and W. Boscheri and M. Semplice and G. Russo},
  title={Central weighted {ENO} schemes for hyperbolic conservation laws on fixed and moving unstructured meshes},
  journal={SIAM Journal on Scientific Computing}, volume={39}, number={6}, pages={A2564--A2591}, year={2017},
  doi={10.1137/17M1111036}}

@article{cravero2018cweno,
  author={I. Cravero and G. Puppo and M. Semplice and G. Visconti},
  title={{CWENO}: uniformly accurate reconstructions for balance laws},
  journal={Mathematics of Computation}, volume={87}, number={312}, pages={1689--1719}, year={2018},
  doi={10.1090/mcom/3273}}

@article{PareschiRusso2005,
  author={L. Pareschi and G. Russo}, title={Implicit-explicit {Runge--Kutta} schemes and applications to hyperbolic systems with relaxation},
  journal={Journal of Scientific Computing}, volume={25}, pages={129--155}, year={2005}}

@article{Rusanov1962,
  author={V. V. Rusanov},
  title={Calculation of interaction of non-steady shock waves with obstacles},
  journal={USSR Computational Mathematics and Mathematical Physics},
  volume={1}, number={2}, pages={267--279}, year={1962},
  doi={10.1016/0041-5553(62)90062-9}}

@article{SaadSchultz1986,
  author={Y. Saad and M. H. Schultz}, title={{GMRES}: a generalized minimal residual algorithm for solving nonsymmetric linear systems},
  journal={SIAM Journal on Scientific and Statistical Computing}, volume={7}, pages={856--869}, year={1986},
  doi={10.1137/0907058}}

@article{Falle2002,
  author={S. A. E. G. Falle},
  title={Rarefaction shocks, shock errors, and low order of accuracy in {ZEUS}},
  journal={The Astrophysical Journal Letters}, volume={577}, number={2},
  pages={L123--L126}, year={2002}, doi={10.1086/344336}}

@article{FalleKomissarov2001,
  author={S. A. E. G. Falle and S. S. Komissarov},
  title={On the inadmissibility of non-evolutionary shocks},
  journal={Journal of Plasma Physics}, volume={65}, number={1},
  pages={29--58}, year={2001}, doi={10.1017/S0022377801008856}}

@article{OrszagTang,
  author={S. A. Orszag and C. M. Tang}, title={Small-scale structure of two-dimensional magnetohydrodynamic turbulence},
  journal={Journal of Fluid Mechanics}, volume={90}, number={1}, pages={129--143}, year={1979}}

@article{DahlburgPicone,
  author={R. B. Dahlburg and J. M. Picone}, title={Evolution of the {Orszag--Tang} vortex system in a compressible medium. I. Initial average subsonic flow},
  journal={Physics of Fluids B}, volume={1}, number={11}, pages={2153--2171}, year={1989}}

@article{PiconeDahlburg,
  author={J. M. Picone and R. B. Dahlburg}, title={Evolution of the {Orszag--Tang} vortex system in a compressible medium. II. Supersonic flow},
  journal={Physics of Fluids B}, volume={3}, number={1}, pages={29--44}, year={1991}}

@article{Dematte2024MHD,
  author  = {Dematt{\'e}, Riccardo and Farmakalides, Alexander A.
             and Millmore, Stephen and Nikiforakis, Nikos},
  title   = {An all Mach number scheme for visco-resistive
             magnetically-dominated MHD flows},
  journal = {Journal of Computational Physics},
  volume  = {514},
  pages   = {113229},
  year    = {2024},
  doi     = {10.1016/j.jcp.2024.113229}
}

@article{vanLeer1977,
  author  = {van Leer, Bram},
  title   = {Towards the ultimate conservative difference scheme. {IV}. {A} new approach to numerical convection},
  journal = {Journal of Computational Physics},
  volume  = {23},
  number  = {3},
  pages   = {276--299},
  year    = {1977},
  doi     = {10.1016/0021-9991(77)90095-X}
}

@misc{BernardelliBoscheriCardonnaALE,
      title={Involution-Preserving schemes through constrained mesh motion: a novel Arbitrary-Lagrangian-Eulerian framework},
      author={Sacha Cardonna and Elena Bernardelli and Walter Boscheri},
      year={2026},
      eprint={2611.12188},
      archivePrefix={arXiv},
      primaryClass={math.NA},
      url={https://arxiv.org/abs/2611.12188},
}
	
\end{document}